\documentclass[times,preprint,3p]{elsarticle}
\journal{Journal of Computational Physics: X}
\usepackage{framed,multirow}

\usepackage{bm}
\usepackage{multirow}
\usepackage{hyperref}
\usepackage{color}
\usepackage{amsfonts}
\usepackage{amsmath}
\usepackage{float}
\usepackage{graphicx}
\usepackage{rotating}
\usepackage{mathtools}
\usepackage[shortlabels]{enumitem}
\usepackage{siunitx} 
\usepackage{graphicx}
\usepackage{subcaption}
\usepackage{hhline}
\usepackage{hyperref}
\usepackage[capitalise]{cleveref}

\hypersetup{
    colorlinks=true,
    linkcolor=blue,
    filecolor=magenta,      
    urlcolor=cyan,
}

\usepackage{booktabs}
\usepackage{lineno}

\newcommand{\bfx}{\mathbf{x}}

\newcommand{\bfz}{\mathbf{z}}

\newcommand{\bfC}{\mathbf{C}}

\newcommand{\bfK}{\mathbf{K}}

\usepackage{comment}

\newcommand{\mm}[1]{\rm mm}

\newcommand{\beq}{\begin{equation}}
\newcommand{\eeq}{\end{equation}}
\newcommand{\bea}{\begin{eqnarray}}
\newcommand{\eea}{\end{eqnarray}}

\newcommand{\dt}{\Delta t}
\newcommand{\dx}{\Delta x}
\newcommand{\dy}{\Delta y}
\newcommand{\dz}{\Delta z}

\newcommand{\bit}{\begin{itemize}}
\newcommand{\eit}{\end{itemize}}
\newcommand{\ben}{\begin{enumerate}}
\newcommand{\een}{\end{enumerate}}

\newcommand{\bF}{\mathbf{F}}
\newcommand{\bG}{\mathbf{G}}

\newcommand{\bI}{\mathbf{I}}

\newcommand{\bK}{\mathbf{K}}
\newcommand{\bk}{\mathbf{k}}

\newcommand{\bU}{\mathbf{U}}
\newcommand{\bV}{\mathbf{V}}

\newcommand{\be}{\mathbf{e}}

\newcommand{\bq}{\mathbf{q}}

\newcommand{\btt}{\mathbf{t}}

\newcommand{\bx}{\mathbf{x}}
\newcommand{\by}{\mathbf{y}}
\newcommand{\bz}{\mathbf{z}}

\newcommand{\erf}[1]{\operatorname{erf} \left [ #1 \right ]}
\newcommand{\expo}[1]{\exp \left [ #1 \right ]}

\makeatletter
\def\@xfootnote[#1]{%
  \protected@xdef\@thefnmark{#1}%
  \@footnotemark\@footnotetext}
\makeatother

\usepackage{url}
\usepackage{xcolor}
\definecolor{newcolor}{rgb}{.8,.349,.1}

\begin{document}
\begin{frontmatter}

\title{GP-MOOD: A positive-preserving high-order finite volume method 
for hyperbolic conservation laws}

\author[1,2,3]{R\'{e}mi {Bourgeois}}
\ead{remi.bourgeois@cea.fr}
\author[1]{Dongwook {Lee}\corref{cor1}}
\cortext[cor1]{Corresponding author: Tel.: +1-831-502-7708}
\ead{dlee79@ucsc.edu}

\address[1]{Department of Applied Mathematics, 
The University of California, Santa Cruz, CA, United States}
\address[2]{Bordeaux INP, Enseirb-Matmeca, France}
\address[3]{Maison de la Simulation, CEA Saclay, Université Paris-Saclay, Gif-sur-Yvette, France}

\begin{abstract}
We present an \textit{a posteriori} shock-capturing finite volume method algorithm called GP-MOOD that solves a compressible hyperbolic conservative system at high-order solution accuracy (e.g., third-, fifth-, and seventh-order) in multiple spatial dimensions. The GP-MOOD method combines two methodologies, the polynomial-free spatial reconstruction methods of GP (Gaussian Process) and the \textit{a posteriori} detection algorithms of MOOD (Multidimensional Optimal Order Detection). The spatial approximation of our GP-MOOD method uses GP's unlimited spatial reconstruction that builds upon our previous studies on GP reported in Reyes  \textit{et al.}, Journal of Scientific Computing, 76 (2017) and Journal of Computational Physics, 381 (2019). This paper focuses on extending GP's flexible variability of spatial accuracy to an \textit{a posteriori}  detection formalism based on the MOOD approach. We show that GP's polynomial-free reconstruction provides a seamless pathway to the MOOD's order cascading formalism by utilizing GP's novel property of variable $(2R+1)$th-order spatial accuracy on a multidimensional GP stencil defined by the GP radius $R$, whose size is smaller than that of the standard polynomial MOOD methods.  The resulting GP-MOOD method is a positivity-preserving method. We
examine the numerical stability and accuracy of GP-MOOD on smooth and discontinuous flows in multiple spatial dimensions without resorting to any conventional, computationally expensive  \textit{a priori} nonlinear limiting mechanism to maintain numerical stability.

\end{abstract}

\end{frontmatter}



\section{Introduction}\label{sec:intro}
Nonlinear shock waves and discontinuities are ubiquitous in compressible
flows. Therefore, a strong gradient control algorithm is an essential building block 
in a numerical method for compressible flow simulations.
A well-suited numerical discretization should be accompanied by a robust mechanism
that can feel
such gradient flows and
evolve its solution stably and accurately, suppressing any
adverse evolution of numerical oscillations and unphysical states. 
Practitioners have devised novel ways to meet such computational needs
to maximize solution accuracy (at the cost of lowering numerical dissipation) on smooth flows while maximizing
numerical stability (or increasing numerical dissipation at the cost of lowering solution accuracy)
on strong gradient flows. Numerically, the corresponding resolution can be classified as one of the two approaches,
namely: an \textit{a priori} method and an \textit{a posteriori} method, 
depending on how a given method detects and controls smooth vs. non-smooth local flows.

A more popular and traditional path has been the \textit{a priori} approach, in which
all cells are spatially discretized by utilizing one or more nonlinear limiting procedures, 
integrated with a stable temporal update to the next time step.
Well-known examples of \textit{a priori} limiting methods include
second-order piecewise linear TVD (Total Variation Diminishing) 
methods (e.g., ~\cite{van1974towards,tadmor1988convenient,hubbard1999multidimensional,harten1997high}),
higher-order polynomial approximations such as 
piecewise parabolic method (PPM)~\cite{colella1984piecewise,mccorquodale2011high},
essentially non-oscillatory methods (ENO) (e.g., ~\cite{harten1997uniformly,shu1988efficient}),
weighted ENO (WENO) methods (e.g., ~\cite{liu1994weighted,jiang1996efficient,balsara2000monotonicity,gerolymos2009very})
and other variants such as central WENO (CWENO)~\cite{levy2000compact,ivan2014high,semplice2016adaptive,dumbser2017central},
Hermite WENO (HWENO)~\cite{qiu2004hermite},
adaptive-order WENO (AO-WENO)~\cite{balsara2007sub,balsara2016efficient},
polynomial-free Gaussian Process WENO (GP-WENO)~\cite{reyes2018new,reyes2019variable,reeves2020application},
to name a few.

In these shock-capturing methods, nonlinear mechanisms or switches are necessary to detect the magnitude of local flow gradients in an \textit{a priori} fashion to evolve numerical solutions stably while meeting underlying physical principles (e.g., positivity, conservation).
%
From the computational perspective, major drawbacks in these  \textit{a priori} 
nonlinear switches are two-fold: ``the computational expense'' due to the required cell-by-cell
calculation of nonlinear limiters and 
``the presence of unavoidable numerical dissipation'' resulting in
degraded solution accuracy. In \cite{kent2014determining_part2}, Kent \textit{et al.} studied the impact of 
monotone and quasi-monotone limiters on effective grid resolutions for finite difference (FD) and finite volume (FV)
methods. The study shows that a large increase in numerical diffusion and dispersion errors
is observed at the first time step, significantly reducing the effective grid resolution compared to the
corresponding unlimited schemes. The disparity between the limited and unlimited schemes then
grows at a slower rate as simulations progress to longer time steps.

Compared to the \textit{a priori} approach, the \textit{a posteriori} approach is a relatively newer paradigm
for compressible shock-capturing schemes.
For this approach, there have been growing interests in the so-called Multidimensional
Optimal Order Detection (MOOD) method in the past decade, which takes a different strategy in handling 
shocks and discontinuities in compressible flow simulations. Since its first introduction in
\cite{clain2011high}, the MOOD method has become a new \textit{a posteriori} detection
principle alternative to the conventional \textit{a priori} shock-capturing techniques based on
nonlinear limiting procedures. The original MOOD in \cite{clain2011high} focused on developing
high-order (up to third-order) two-dimensional polynomial approximations on unstructured grids.
This study introduced the first idea of its kind where the initial unlimited polynomial order of accuracy
cascades down to lower ones until a discrete candidate solution on each cell meets
a predefined admissibility condition called the Discrete Maximum Principle (DMP). 
In this context, the troubled cells that fail to meet the prescribed DMP condition need to be recomputed
locally, hence referred to as an \textit{a posteriori} scheme.
In \cite{diot2012improved}, the accuracy of the MOOD method 
was extended to sixth-order on unstructured grids
with an improved detection criterion called the $u$2 detection 
that relaxes the order-reduction criterion on smooth extrema.
Besides, the positivity-preserving property under an admissible 
CFL (Courant-Fredirichs-Lewy) condition was made explicit by 
imposing the Physical Admissibility Detection (PAD)
that checks the positivity of density and pressure variables on each cell.
The MOOD method was further extended to general three-dimensional unstructured meshes
with simplifications of the $u$2 detection \cite{diot2013multidimensional,diot2012methode};
to a third-order FV adaptive mesh refinement (AMR) scheme in \cite{semplice2018adaptive}.

The MOOD paradigm of ``repeat-until-valid'' has also been adopted in the Discontinuous Galerkin (DG)
method, where the first pass of discrete updates is done by a target high-order \textit{unlimited} 
DG solver. Then, the invalid cells that fail to meet prescribed admissibility conditions such as 
PAD and Numerical Admissibility Detection (NAD) are recomputed using other more stable
(but lower accurate) \textit{limited} solvers (such as WENO) at a more refined subcell resolution 
on each of those invalid cell,
e.g., see examples of the MOOD limiter in 
ADER\footnote{ADER stands for Arbitrary high-order DERivative, first proposed by Toro and his
collaborators \cite{toro2001towards,titarev2002ader,titarev2005ader}.}-DG applications
~\cite{dumbser2014posteriori,dumbser2016simple}.

On the other hand, there has been a new attempt \cite{bourriaud2020priori} 
to use a trained Neural Network (NN) with 
hidden layers and few perceptrons to design an NN-based classification model. In this NN approach,
the standard \textit{a posterior} MOOD detection strategy is replaced by an educated guess of
the appropriate order of polynomial accuracy. 
Thus, the potential imbalance in parallel efficiency caused by 
the extra recomputing need in the standard MOOD treatment could be ameliorated.
Moreover, the NN-based approach alters the original \textit{a posteriori} nature of the MOOD paradigm
to a trained \textit{a priori} version, which would be better suited for extending the MOOD method to an
implicit scheme. Although promising, the work in \cite{bourriaud2020priori} has some limitations; namely,
it only explored shallow NNs of two hidden layers; investigated 1D tests only; 
trained their NNs on a relatively moderate size of a training set; the positivity and \texttt{NAN}
checks are not in the part of NNs but 
imposed separately as part of the \textit{a posteriori} MOOD loop.

There are two core building blocks in the MOOD method. The first building block is the MOOD detection criteria
themselves to tag invalid cells. These problematic cells are subsequently recomputed using 
a more stable but less accurate discretization method, repeating the so-called MOOD iterative loop
until they become valid.
The second building block is a set of multiple choices 
(at least two, e.g., a high-order unlimited method and a first-order method) of spatial
discretization methods so that the invalid cells become valid 
physically (i.e., PAD) and numerically (i.e., NAD) at the end of the MOOD loop.
To date, the MOOD method has exclusively employed the use of multidimensional polynomial reconstruction
\cite{clain2011high,diot2012improved,diot2013multidimensional,diot2012methode,bourriaud2020priori}.
By far, piecewise polynomial reconstruction and interpolation of stencil data have been among the most
understood and well-established mathematical tools
in computing discrete approximations to underlying functions in computational fluid dynamics (CFD).

Regardless, there seem at least two challenges in function approximations 
using piecewise multidimensional polynomials.
The first issue appears to be well-known complications where the interpolated/reconstructed solutions of 
piecewise multidimensional polynomials may lead to an ill-conditioned or unsolvable linear system.
One workaround could be to utilize a tensor product of 
one-dimensional $n$th-order (or $(n-1)$th degree) polynomials
~\cite{shi_technique_2002,shu_high-order_2003,balsara2009divergence}
whose degrees of freedom (i.e., the number of unknown polynomial coefficients) 
is $\mbox{DoF}(n)=n(n+1)/2$. 
Another workaround is to solve a least-squares problem (LSP) 
to compute polynomial coefficients
with a number of stencil points larger
than $\mbox{DoF}(n)$~\cite{mccorquodale2015adaptive,zhang2012fourth}.
The LSP pathway is also taken by the existing MOOD methods
~\cite{clain2011high,diot2012improved,diot2013multidimensional,diot2012methode,
semplice2018adaptive,bourriaud2020priori},
wherein 2D at least 5 cells are required for 2nd-order MOOD, 
8 cells for 3rd-order MOOD,
16 cells for 4th-order MOOD, 20 cells for 5th-order MOOD, 
and 26 cells for 6th-order MOOD~\cite{diot2012improved}.
The second challenging issue is the computational cost.
Ideally, in the existing MOOD approach, a new least-squares solution
is required for all cells in the first pass
and for each troubled cell at every subsequent MOOD iteration to compute lower-order polynomials.
To minimize this computationally intensive polynomial reconstruction, 
an effort of `truncating' the highest degree polynomial to 
lower-order polynomials was proposed in~\cite{clain2011high},
although it was found later to be undesirable at discontinuities for 3rd-order or 
higher MOOD methods~\cite{diot2012improved}.

Here lies the main emphasis of this paper. In the context of FV MOOD, we propose a new set of
Gaussian Process (GP) reconstruction methods that are polynomial-free and computationally efficient. 
We will show that the prime advantage in our GP reconstruction 
is afforded by GP's re-usable high-order spatial reconstructors that seamlessly vary its solution accuracy. 
Thereby, we advance the MOOD's second building block 
significantly in solution accuracy, stability, and efficiency. 
Our approach builds around the recent developments of GP methods
~\cite{reyes2018new,reyes2019variable,reeves2020application}.
First introduced by Reyes \textit{et al.}~\cite{reyes2018new},
FV GP reconstruction has shown a versatile selectable-order property
in modeling the 1D Euler equations~\cite{reyes2018new}.
The FV GP method followed WENO's \textit{a priori} shock-handling control
with modified GP smoothness indicators, called GP-WENO.
An extension of GP-WENO to a full 3D FD method (FD-prim) was reported in~\cite{reyes2019variable}.
In the most recent work~\cite{reeves2020application} by Reeves \textit{et al.}, 
the GP-WENO method was extended to a 3rd-order prolongation
algorithm in FV AMR simulations using AMReX~\cite{amrex}.
In addition, we further relax the existing MOOD admissibility conditions in the current study.
This new  relaxed detection strategy is found to be crucial for simulating
key flow structures on some of the well-known benchmark test problems
(e.g., non-systematic solution behaviors in the 1D Shu-Osher test
in \cref{sec:shu_osher}; 
the diagonal jet formation in the 2D implosion test in \cref{sec:Implosion}).
%
The new relaxation improves the issue
in the existing MOOD method, in which the order decrement 
takes place on cells where the local flow is weakly compressible
but has not yet built steep gradients.
Such local flows are an important powerhouse 
to trigger a more nonlinear regime 
(e.g., formations of shocks, vortex roll-ups, reflected jets)
in the subsequent flow dynamics.
However, in the existing MOOD method, 
we observed that a low order method is activated too soon on such cells, 
which often suppresses the onset of such signature flow 
structures by the excessive numerical dissipation.

We organize the paper as follows. In \cref{sec:method} we briefly review the discretization
strategy of the existing MOOD method. We present a group of unlimited, selectable-order GP
reconstruction methods whose spatial accuracy is linearly varying 
at $(2R+1)$th-order with an integer GP radius $R$.
Demonstrated in this paper include GP methods of 3rd-order ($R=1$), 5th-order ($R=2$), 
and 7th-order ($R=3$). 
In \cref{sec:method} we describe a new relaxed MOOD strategy that combines
a multidimensional shock-detector and the MOOD method.
A description of a step-wise implementation of GP-MOOD is given in \cref{sec:stepwise}.
The test results of our new GP-MOOD method are presented in \cref{sec:results}.
We conclude our paper with a brief summary in \cref{sec:conclusion}.

\section{Finite volume discretization of GP-MOOD}\label{sec:method}
The FV spatial discretization of the GP-MOOD method is introduced in this section. 
We focus on designing a set of selectable-order GP-MOOD algorithms on Cartesian grid configurations,
which reconstruct point-wise Riemann states at cell faces 
using cell-centered volume averages. 
As will be shown, the GP reconstructor is a covariance kernel-based
posterior mean function that computes the $(2R+1)$th-order accurate
Riemann states at cell faces \textit{directly} (i.e., without the conventional Taylor series expansion)
from the cell-centered average values on the multidimensional GP stencil of radius $R$.
The Riemann states are computed at multiple Gaussian
quadrature points on each cell face to deliver the anticipated $(2R+1)$th spatial accuracy.
The spatial solutions are temporally evolved with
strong stability preserving Runge-Kutta (SSP-RK) methods.

\subsection{Governing equations and the basic finite volume discretization form}\label{sec:governing_eqns}
We are interested in solving a general hyperbolic system of conservative laws in 2D,
\begin{equation}\label{eq:gov}
    \partial_{t} \bU + \partial_{x} \bF (\bU) + \partial_{y} \bG (\bU)  = 0,
\end{equation}
where \( \bU \) is the vector of conservative variables and
$\bm{\mathcal{F}}= (\bF, \bG)$ are the flux functions in
\( x \)- and \( y \)-direction.
In the Euler equations in 2D, the conservative variables and
the flux functions are defined as,

\begin{equation}
    \bU = \begin{bmatrix}
        \rho \\
        \rho u \\
        \rho v \\
        \rho E
    \end{bmatrix},\quad
    \bF (\bU) = \begin{bmatrix}
        \rho u \\
        \rho u^{2} + p \\
        \rho u v \\
        u \left(\rho E + p \right)
    \end{bmatrix}, \quad
    \bG (\bU) = \begin{bmatrix}
        \rho v \\
        \rho u v \\
        \rho v^{2} + p \\
        v \left( \rho E + p \right)
    \end{bmatrix}, \quad
    \label{eq:euler_eqn}
\end{equation}
where $\rho$ denotes the fluid density, $u$ and $v$ represent 
the $x$ and $y$ fluid velocity respectively, 
and $\rho E$ is the total energy.
The system is closed with an ideal gas equation of state (EoS),
\(p = (\gamma - 1)\left(\rho E - \frac{1}{2}\rho({u^2 + v^2)}\right), \)
where $\gamma$ is the ratio of specific heats. 
The hyperbolic system in Eq.~\eqref{eq:euler_eqn} is physically admissible if both $p>0$ and $\rho>0$,
and the numerical method that maintains the positivity property is referred to as
a positivity-preserving method.
%

The basic form of the finite volume discretization of \cref{eq:gov} is derived by integrating
the equation over each cell $I_{ij}=[x_{i-1/2},x_{i+1/2}]\times  [y_{j-1/2},y_{j+1/2}]$ and over 
a time interval $[t^n, t^{n+1}]$, yielding
\beq\label{eq:fvm_discrete}
\overline{\bU}^{n+1}_{ij} = \overline{\bU}^{n}_{ij} - 
\dt \,\mathbb{F}_{\nabla}, 
\eeq
where $\overline{\bU}^{n}_{ij}$ is a vector of the volume-averaged conservative variables and
$\mathbb{F}_{\nabla}$ is a collection of the rest of spatial derivatives terms, including the face-averaged
and temporally-averaged flux in each spatial direction. 
A popular choice for the discrete temporal update in \cref{eq:fvm_discrete}
for high-order simulations is to use either 
a one-step ADER method~\cite{toro2001towards,titarev2002ader,titarev2005ader}
or a multi-stage SSP-RK method~\cite{gottlieb1998total,spiteri2002new}
(the latter is our choice in this paper).

It leaves us to determine how to update the spatial approximations 
of the relevant terms in $\mathbb{F}_{\nabla}$ 
to meet the expected high-order accuracy. 
To this end, we design a family of multidimensional FV reconstruction algorithms of GP with the 
following set of attentions~\cite{lee2017piecewise,diot2012improved}:
\bit
\item[(i)] make a clear distinction between the integral quantities (e.g., face-averaged, volume-averaged) and the pointwise values,
\item[(ii)] use volume-averaged conservative quantities in high-order reconstruction, viz., avoid the use of
primitive variables in reconstruction, which often are derived by inappropriate nonlinear operations such as 
the $x$-velocity computed as a division of the volume-averaged $x$-momentum by the volume-averaged density,
e.g., $u \leftarrow \overline{\rho u}/ \overline{\rho}$ or $\overline{u} \leftarrow \overline{\rho u}/ \overline{\rho}$,
\item[(iii)] avoid any nonlinear operations that use inappropriately converted primitive variables to construct conservative variables
and use them in the evolution of conservative solution updates, e.g., any nonlinear operations that
use the $x$-velocity in (ii) to replace any volume-averaged conservative variables,
$\overline{\rho u} \leftarrow (\overline{\rho}) (\overline{u}$).
\eit

In this paper, we use a $q$-point Gaussian quadrature rule to approximate 
the \textit{face-averaged} fluxes at $2q$-th order accuracy using $q$ many
\textit{pointwise} fluxes 
on each cell face.
For the sake of simplicity, we assume a uniform Cartesian grid configuration in 2D. 
This gives us to write $\mathbb{F}_{\nabla}$ as
\beq\label{eq:fvm_F}
\mathbb{F}_{\nabla} = 
\frac{1}{\dx}
\sum_{j_g=1}^{q}
\omega_{j_g}\Big( 
\bF^*_{i+1/2,j_g} - \bF^*_{i-1/2,j_g}
\Big)
+
\frac{1}{\dy}
\sum_{i_g=1}^{q}
\omega_{i_g}\Big( 
\bG^*_{i_g,j+1/2} - \bG^*_{i_g,j-1/2}
\Big),
\eeq
where $i_g$ and $j_g$ are the indices of the $q$-point Gaussian quadrature point locations 
on each $x$ and $y$ cell face; the corresponding $\omega_{i_g}$ and $\omega_{j_g}$ 
are the quadrature weights for the $2q$-th order numerical integration. 
The numerical fluxes $\bF^*$ and $\bG^*$ are \textit{pointwise} 
fluxes at each respective cell face, obtained by solving the corresponding Riemann problems
at the Gaussian quadrature points.
A pair of high-order accurate \textit{pointwise} Riemann states, $(\bU_{L},\bU_{R})$, 
are used as inputs to calculate the corresponding Riemann problems at each quadrature point.
In each pair, the left $\bU_{L}$ and the right $\bU_{R}$ states 
are computed using a $(2R+1)$th-order GP reconstruction method that takes inputs of
the \textit{volume-averaged} cell-centered conservative variables $\overline{\bU}^n_{ij}$
on a multidimensional GP stencil determined by the size of a GP radius $R$.

In the next section, we describe a family of $(2R+1)$th-order GP reconstruction methods.
More specifically, we aim to provide a family of three GP solvers, 
including a 3rd-order GP method with $R=1$ (GP-R1), 
a 5th-order GP method with $R=2$ (GP-R2),
and a 7th-order GP method with $R=3$ (GP-R3).
It should be noted that our GP algorithm is general and 
is not limited to the 7th-order GP-R3 method.
As expected, however, the computational cost increases
as $R$ increases, primarily due to 
the increase in the size of the resulting linear system to solve and 
the increase in the number of the needed Gaussian quadrature points.
For this reason, we present our GP method up to 7th-order in this paper.

\subsection{The GP stencil}\label{sec:gp_stencil}



This section presents how to configure local GP stencils
in the context of a two-dimensional FV reconstruction. 
The $(2R+1)$th-order GP reconstruction operates on:
\bit
\item[(i)] Input: a vector $\bq$ consisting of the \textit{volume-averaged} conservative 
variables (e.g., $\overline{\rho}^n_{ij}$) on a 2D local GP stencil of radius $R$, where the GP radius $R$ is a runtime parameter.
\item[(ii)] Output: an \textit{unlimited} $(2R+1)$th-order accurate conservative \textit{pointwise} Riemann state 
of the same input variable (e.g., density) at each Gaussian quadrature point 
(e.g., $\rho_{*}=\rho(\bx_*)$,  where $\bx_* = (x_{i\pm 1/2}, y_{j_g})$ or $\bx_* = (x_{i_g}, y_{j\pm 1/2})$) 
\eit

Let us first give a geometrical description of how our 
GP local stencil of size $R$ is configured. 
We begin with the simplest case with $R=1$ for
the 3rd-order GP reconstruction. The GP radius $R$ is an integer value 
in the unit of grid spacing, $\dx$ and $\dy$.
For example, in a 2D uniform grid configuration
the GP stencil of $R=1$ defines a five-point cross-shape stencil that
extends the local stencil centered at $\bx_{ij}=(x_i,y_j)$ 
to one neighboring cell in each $x$ and $y$ direction,
drawing hypothetically a blocky-diamond in 2D around $\bx_{ij}$.
In general, $\dx$ could be different from $\dy$, in which case the GP stencil
becomes a stretched-cross consisting of the five
volume-averaged data around $\bx_{ij}$.

\begin{figure}[h!]
    \centering
  \includegraphics[scale=0.3]{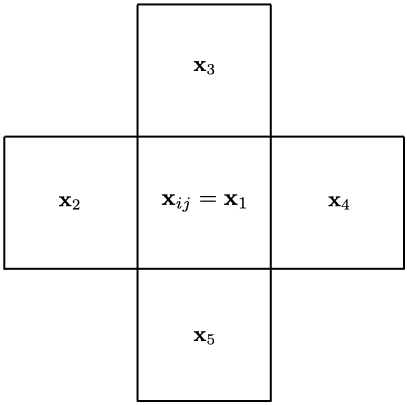}
    \caption{The five-point GP stencil of radius $R=1$ for the 3rd-order GP 
    reconstruction method. The ordered labeling
    illustrates how the local volume-averaged conservative variables at $t=t^n$
    are rearranged into a one-dimensional five-point array, $\overline{\bq}_{ij}$.}
    \label{fig:gp_stencil_GP-R1}
\end{figure}

\cref{fig:gp_stencil_GP-R1} illustrates how the GP-R1 stencil is
configured around the central cell $\bx_{ij}$. For exposition
purposes, we also show ordered labeling, which is used
for reshaping the five volume-averaged quantities 
at those cells in \cref{fig:gp_stencil_GP-R1}
into a one-dimensional
array, denoted by $\overline{\bq}_{ij}$. 
At each timestep $t=t^n$, the five local volume-averaged conservative variables,
$\bar{q}_m=\bar{q}(\bx_m, t^n)$, $m=1, \dots, 5$,
are cast into $\overline{\bq}_{ij}$ in the
order as shown in \cref{fig:gp_stencil_GP-R1}, 
starting from the data $\bar{q}_1$ at the central cell $\bx_{ij}$,
\beq
\overline{\bq}_{ij} = 
\left(
\bar{q}_1, \bar{q}_2, \bar{q}_3, \bar{q}_4, \bar{q}_5
\right)^T.
\eeq

\begin{figure}[h!]
    \centering
    \begin{subfigure}{80mm} 
    \centering
    \includegraphics[scale=0.25]{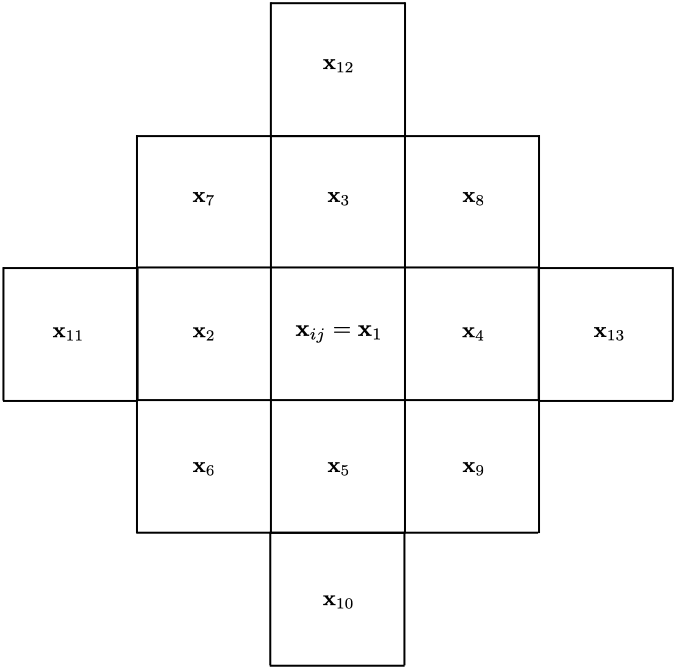} \caption{}\label{subfig:GP-R2}
    \end{subfigure}
    \begin{subfigure}{80mm}
    \centering
    \includegraphics[scale=0.25]{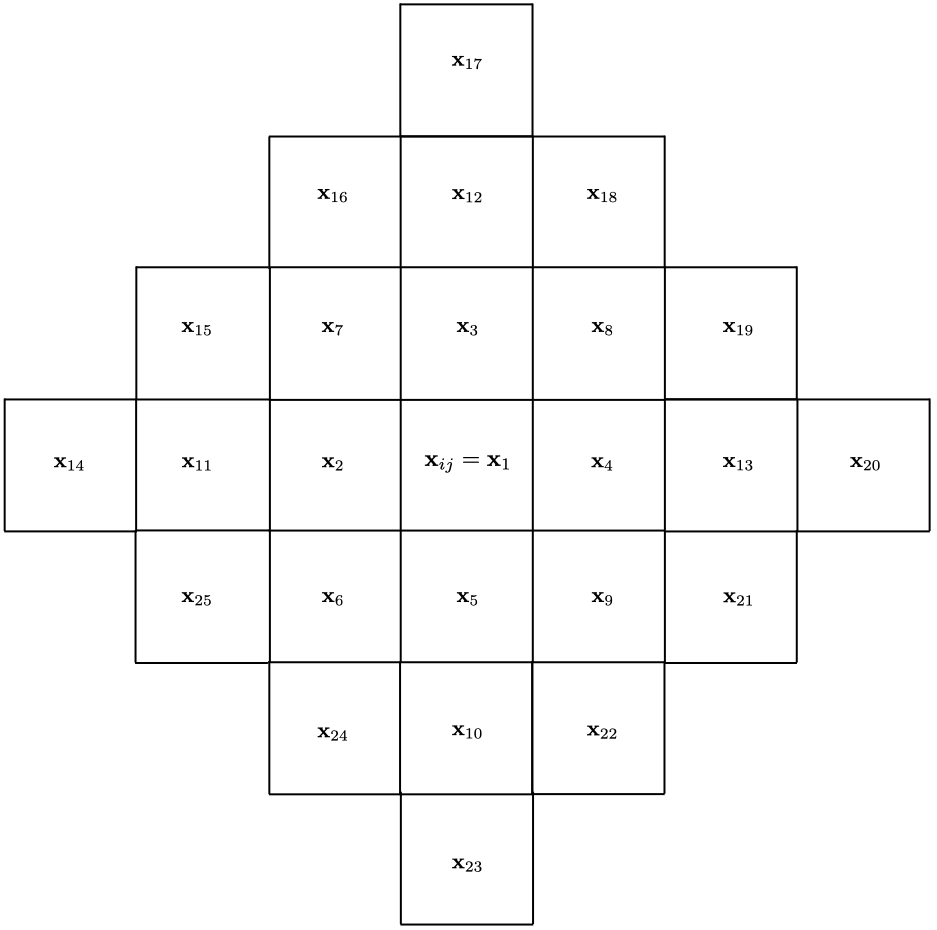}  \caption{}\label{subfig:GP-R3}
    \end{subfigure}
    \caption{(a) The 13-point GP stencil of radius $R=2$ for the 5th-order GP 
    reconstruction method. 
    (b) The 25-point GP stencil of radius $R=3$ for the 7th-order GP 
    reconstruction method. 
    In both, the ordered labeling illustrates how the local volume-averaged conservative variables at $t=t^n$
    are rearranged into a one-dimensional five-point array, $\overline{\bq}_{ij}$.}
    \label{fig:gp_stencil_GP-R2-R3}
\end{figure}

In a similar fashion, we form the blocky-diamond GP stencils that group
the neighboring 13 cell data for the GP-R2 stencil in \cref{subfig:GP-R2}
and the 25 cell data for the GP-R3 stencil in \cref{subfig:GP-R3}.
The one-dimensional array $\overline{\bq}_{ij}$ will hold 13 cell data for GP-R2
and 25 for GP-R3, following the orderly fashion indicated in \cref{fig:gp_stencil_GP-R2-R3}.
This finalizes our discussion on the GP stencils for the 3rd-, 5th-, and 7th-order methods.
We remark that the size of our GP stencils is smaller than the local stencil size
of the existing polynomial-based MOOD methods at the same or comparable accuracy in 2D
\cite{diot2012improved}, e.g., 
5 cells for the 2nd-order polynomial MOOD methods,
8 for 3rd-order, 
20 for 5th-order, and
26 for 6th-order.

\begin{figure}[h!]
    \centering
  \includegraphics[scale=0.25]{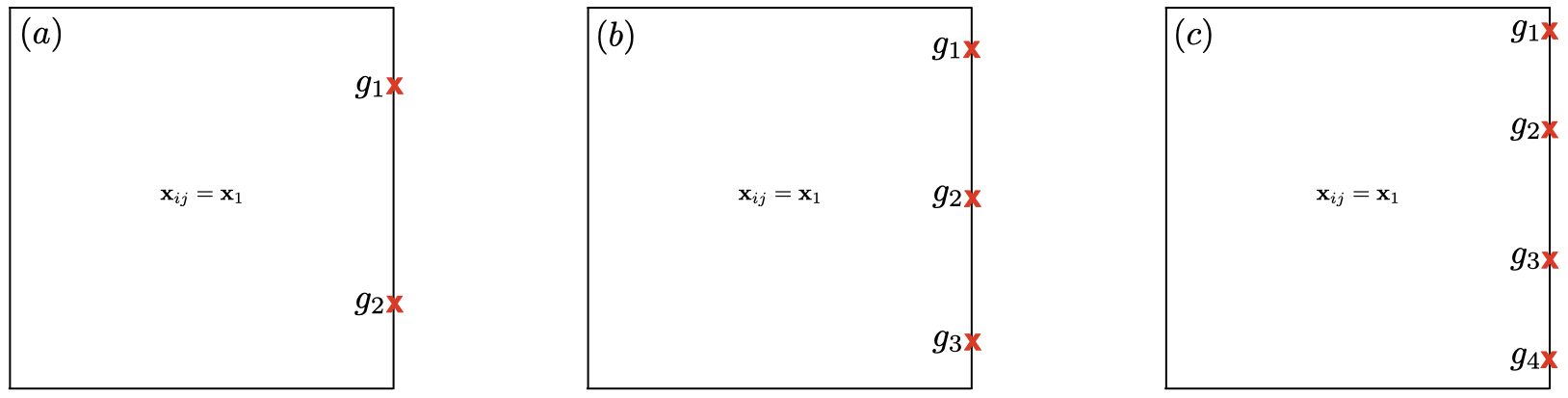}
    \caption{Multi-point Gaussian quadrature points, $g_m$, at the $x$-normal cell-face
    of the central cell $\bx_{ij}=\bx_1$. 
    (a) The 4th-order two-point quadrature rule with two quadrature points is
    used with the 3rd-order GP-R1 method.
    (b) The 6th-order three-point quadrature rule with three quadrature points is
    used with the 5th-order GP-R2 method.
    (c) The 8th-order four-point quadrature rule with four quadrature points is
    used with the 7th-order GP-R3 method.
    The values of the quadrature points $g_m$ and the corresponding weights $\omega_m$
    are given in \cref{tab:QR-pt-weight}.
    }
    \label{fig:multipt-QR}
\end{figure}

The orderly-grouped GP stencil data in $\overline{\bq}_{ij}$ are used as inputs
to the GP reconstructor. We extend our earlier work in \cite{reyes2018new} to designing two-dimensional
GP reconstruction schemes that convert the local volume-averaged 
data $\overline{\bq}_{ij}$ to a $(2R+1)$th-order accurate pointwise conservative quantity
at each Gaussian quadrature point, $g_m$, as displayed in \cref{fig:multipt-QR}.

%
%

\begin{table}[ht!]
    \footnotesize
    \centering
    \caption{Multi-point Gaussian quadrature rules (QR) are
    used in combination with GP.
    The quadrature points, $g_m$, are tabulated over the reference interval $[-0.5, 0.5]$ that
    maps to the unit length of each cell-face, e.g., $[y_{j-1/2}, y_{j+1/2}]$ at the $x$-normal cell-face.
    See also \cref{fig:multipt-QR}.
    }\label{tab:QR-pt-weight}
    \begin{tabular}{@{}cccccccccccc@{}}
        \toprule
        {QR}             &    $g_1$        & $\omega_1$     &  
                            &     $g_2$       & $\omega_2$     & 
                            &     $g_3$       & $\omega_3$     & 
                            &     $g_4$       & $\omega_4$    \\                                                                  
        \midrule
        2-point QR  & $ \frac{1}{2\sqrt{3}}$                   & $\frac{1}{2}$   &  
                            & $-\frac{1}{2\sqrt{3}}$                  & $\frac{1}{2}$    & 
                            &  ---                                               & ---                  &
                            &  ---                                               & --- \\
        3-point QR & $\frac{1}{2}\sqrt{\frac{3}{5}}$     & $\frac{5}{18}$    &   
                           & 0.0                                              & $\frac{8}{18}$    & 
                           & $-\frac{1}{2}\sqrt{\frac{3}{5}}$   & $\frac{5}{18}$  &  
                           &  ---                                              & --- \\
        4-point QR & $ \frac{1}{2}\sqrt{\frac{3}{7} + \frac{2}{7}\sqrt{\frac{6}{ 5} }}$   & $\frac{18-\sqrt{30}}{72}$ &
                           & $ \frac{1}{2}\sqrt{\frac{3}{7} - \frac{2}{7}\sqrt{\frac{6}{ 5} }}$   & $\frac{18+\sqrt{30}}{72}$ &  
                           & $-\frac{1}{2}\sqrt{\frac{3}{7} - \frac{2}{7}\sqrt{\frac{6}{ 5} }}$   & $\frac{18+\sqrt{30}}{72}$ & 
                           & $-\frac{1}{2}\sqrt{\frac{3}{7} + \frac{2}{7}\sqrt{\frac{6}{ 5} }}$   & $\frac{18-\sqrt{30}}{72}$  \\
    \end{tabular}
\end{table}

\subsection{Basic theory on GP regression}\label{sec:gp_regression}
Formally speaking, a GP is a collection of random variables (agnostic functions in our case),
any finite number of which have a joint Gaussian distribution~\cite{rasmussen2005}.
We follow \cite{reyes2018new} to derive the GP reconstructor from 
the GP regression models~\cite{rasmussen2005}.
The perspective of the GP regression defines GP as a distribution 
over functions, by which GP can make a probabilistic inference on those
sample functions. 
In this way, a GP constructs a random function space
from which the GP can draw samples (i.e., random variables or agnostic functions) 
probabilistically based on the built-in stochastic property in the function space.
A GP is fully specified by the two functions:
\bit
\item a mean function $m(\mathbf{x})$ defined as an expectation of $f(\bx)$, i.e., 
$m(\bx) = \mathbb{E}\Big[f(\bfx)\Big]$ 
over $\mathbb{R}^N$, and
\item a covariance function that is a symmetric, positive-definite kernel
 $K(\mathbf{x},\mathbf{y}) = 
 \mathbb{E}
 \Big[\Big(f(\mathbf{x})-m(\mathbf{x})\Big)
        \Big(f(\mathbf{y})-m(\mathbf{y})\Big)
  \Big]$ 
 over $\mathbb{R}^{N}\times\mathbb{R}^{N}$.
\eit

The two functions give rise to the GP \textit{prior}, leading us to write 
$f \sim \mathcal{GP}\Big(m(\bx),K(\bx,\by)\Big)$
\footnote[$\dagger$]
{It reads as ``a function $f$ is in $\mathcal{GP}.$''}, meaning that
any randomly drawn functions $f$ from this distribution of functions (or GP's function space)
are sampled with mean function $m(\bx)$ and covariance kernel
function $K(\bx, \by)$ probabilistically.

Following \cite{reyes2018new,reyes2019variable,reeves2020application}, 
we choose a constant zero mean function, $m(\bx) = 0$. 
As will be shown, this choice simplifies the GP reconstruction procedure without impacting the solution accuracy.
Also, we set 
the squared exponential (SE) kernel function to be our default choice of $K(\bx,\by)$, defined by
\beq\label{eq:SE}
K(\mathbf{x}, \mathbf{y})=K_{\mathrm{SE}}(\mathbf{x}, \mathbf{y})= 
\exp \left[-\frac{(\mathbf{x}-\mathbf{y})^{2}}{2 \ell^{2}}\right],
\eeq
where the length hyperparameter, $\ell$, controls the characteristic length scale 
of the functions in the GP function space
distributed with $m(\bx)$ and $K(\bx, \by)$.
Operationally, we have found the choice $\ell = 6\Delta$ or $\ell = 12\Delta$, where $\Delta = \min\{\dx, \dy\}$,
produces satisfactory results
~\cite{reyes2018new,reyes2019variable,reeves2020application}, 
which will continue to be our default choice.

The GP approach to regression (or interpolation) begins with making a set of
observed data (or training points). 
To this end, let us first begin with a simple example in the context of the present study, where
the \textit{pointwise}  training points $\bq_{ij}$ are \textit{pointwise} function values $f(\bx_m)$
evaluated at each GP stencil point $\bx_m$. 
Note that we will discuss the \textit{volume-averaged}
training points $\overline{\bq}_{ij}$ in \cref{sec:gp_reconstruction}, which 
is more relevant to the FV reconstruction.
These function values are stored in a one-dimensional array,
$\bq_{ij}$, in the orderly fashion described in 
\cref{fig:gp_stencil_GP-R1,fig:gp_stencil_GP-R2-R3}.
The observed data, $\bq_{ij}$,
are assumed to be \textit{probabilistically known}, or
${\bq}_{ij} \sim \mathcal{GP}\Big(m(\bx),K(\bx,\by)\Big)$,
in terms of the \textit{prior} GP distribution.
Trained on the observed data, the GP regression yields
a pointwise \textit{posterior} distribution on the function values $f(\bx_*)$ at any new point $\bx_*$
by applying the conditioning property of Bayes' Theorem 
to the joint Gaussian distribution on the observed data ${\bq}_{ij}$,
that is, the GP makes inference on  $f(\bx_*)$ given ${\bq}_{ij}$.
Here, $\bx_*$ refers to a new point where no observation has been made,
or in other words, the function $f$ 
has no observed information at $\bx_*$. 

Assuming a zero mean, the conditioning property furnishes 
a new pointwise \textit{posterior mean function} $\tilde{m}(\bx_*)$
\footnote[$\S$]{The resulting GP \textit{posterior} is 
fully defined by the posterior mean function $\tilde{m}(\bx_*)$
and the posterior covariance function $\tilde{\Sigma}$. 
The explicit form of $\tilde{\Sigma}$ is not discussed here as
we do not utilize it for the current study. For details,
see~\cite{rasmussen2005,bishop2007pattern}.}
given by
\beq\label{eq:posterior_mean}
\tilde{m}(\bx_*) 
= \bk_*^T \bK^{-1} {\bq}_{ij}
= \bz_*^T {\bq}_{ij},
\eeq
where $[\bK]_{mn}\equiv K(\bx_m,\bx_n)$ and $[\bk_*]_m = K(\bx_*,\bx_m)$.
The data-independent vector $\bz_*^T= \bk_*^T \bK^{-1}$ is 
called the prediction vector, following the same convention 
in~\cite{reyes2018new,reyes2019variable,reeves2020application}.
For us, $\bx_m$ and $\bx_n$ correspond to the GP stencil coordinates in 
\cref{fig:gp_stencil_GP-R1,fig:gp_stencil_GP-R2-R3}, while
$\bx_*$ is the Gaussian quadrature point locations on each cell-face in
\cref{fig:multipt-QR,tab:QR-pt-weight}. Again, for FV GP-MOOD, we will need to
consider the the volume-averaged quantities $\overline{\bq}_{ij}$ as inputs
instead of the pointwise inputs ${\bq}_{ij}$ in \cref{eq:posterior_mean},
together with the kernel modifications
(see \cref{sec:gp_reconstruction}).
The size of the resulting linear system in \cref{eq:posterior_mean} 
is characterized to be
%
$\bK \in \mathbb{R}^5\times\mathbb{R}^5$ for GP-R1,
$\mathbb{R}^{13}\times\mathbb{R}^{13}$ for GP-R2, and
$\mathbb{R}^{25}\times\mathbb{R}^{25}$ for GP-R3, while
$\bk_*$ and 
${\bq}_{ij} \in \mathbb{R}^5\times\mathbb{R}^1$ for GP-R1,
$\mathbb{R}^{13}\times\mathbb{R}^{1}$ for GP-R2, and
$\mathbb{R}^{25}\times\mathbb{R}^{1}$ for GP-R3.

We emphasize here that the new posterior mean function 
in \cref{eq:posterior_mean} is to be viewed
as a new high-order GP \textit{interpolator} (but not a reconstructor yet)
\cite{reyes2018new,reyes2019variable},
which makes a probabilistic statement about the unobserved 
pointwise function value at $\bx_*$. 
As a CFD interpolation algorithm, 
\cref{eq:posterior_mean} demonstrates a $(2R+1)$th convergence
rate on a smooth input data set, $\bq_{ij}$.
Of importance is the same data type in the input and the output data,
e.g., pointwise function values in both.
For this reason, \cref{eq:posterior_mean} cannot be directly used as
a finite volume reconstruction algorithm where the input and output
data types are different. The required modification for FV will be
discussed in \cref{sec:gp_reconstruction}.

So far, we have briefly outlined the underlying GP-based Bayesian
prior and posterior distributions. As will be shown in the following sections, 
our discussion on the basic Bayesian theory 
is sufficient for us to construct the proposed GP reconstruction methods.
Interested readers are encouraged
to refer to our former studies~\cite{reyes2018new,reyes2019variable,reeves2020application}
for more detailed discussions on the relevant mathematical derivations. 
For a more general discussion on GP theory, see ~\cite{rasmussen2005,bishop2007pattern}.

\subsection{The $(2R+1)$th-order GP reconstruction}\label{sec:gp_reconstruction}
One of the attractive properties in GP regression is that
the underlying GP prediction models, such as \cref{eq:posterior_mean}, are
preserved under a linear operation.
Let us denote a linear operator by $\mathcal{L}$ over the function space
GP has generated.
It can be shown that 
$\mathcal{L}(f) \sim \mathcal{GP}
\bigg(
m \Big(   \mathcal{L}(  \bx ) 
    \Big),
\bK \Big( \mathcal{L}(\bx),  
              \mathcal{L}( \by)
       \Big)
\bigg)$ 
if
$f \sim \mathcal{GP}\Big(m(\bx),\bK(\bx,\by)\Big)$, 
i.e., GP is invariant under linear operations.

This nice invariant property enables us to modify the GP regression model in \cref{eq:posterior_mean} 
to a GP finite volume reconstruction scheme that constitutes two operations in parallel into
one compact linear system, similar to  \cref{eq:posterior_mean}. The two operations
include (i) a high-order data \textit{discretization} at $\bx_*$ using input data, and at the same time
(ii) a high-order data type \textit{conversion} from the volume-averaged input data $\overline{\bq}_{ij}$
to a pointwise output at $\bx_*$.

We apply the same linear operations
taking volume averages of a given function to $m(\bx)$.  
This task can be done by conducting the following three steps \cite{reyes2018new}:
\bit
\item [(a)] Take integral averages of the pointwise function $q$ on each $\bx_m = \bx_{ij}=(x_{i_m},y_{j_m})$ to get
\beq
\overline{q}_{m} 
= \frac{1}{\Delta_{m}}\int_{\mathcal{V}_{m}} q(\bx) d\mathcal{V}(\bx),
\eeq
which simplifies to $\frac{1}{\dx\dy}\int_{I_{j_m}} \int_{I_{i_m}} q(x,y) dx dy$ in a uniform 2D Cartesian geometry
with $\mathcal{V}_{m} = I_{i_m} \times I_{j_m} = [x_{i_m-1/2}, x_{i_m+1/2}] \times [y_{j_m-1/2}, y_{j_m+1/2}]$.
Nothing needs to be done if the volume-averages are given as initialized values; otherwise, 
each pointwise value needs to be converted explicitly to 
the corresponding volume-averaged values; see \cref{sec:isentropic_vortex}.
\item [(b)] Take integral averages of the GP covariance kernel function $\bK$
over both $\bx$ on each $\mathcal{V}_m$ and $\by$ on each $\mathcal{V}_n$ 
to get
a new volume-averaged covariance kernel $\bfC$. Writing in component-wise form,
\beq
[\bfC]_{mn} = \frac{1}{\Delta_m \Delta_n} 
\int_{\mathcal{V}_{n}} \int_{\mathcal{V}_{m}} K(\bx, \by) d\mathcal{V}(\bx) d\mathcal{V}(\by).
\eeq
\item [(c)] Take integral averages of $\bk_* = \bk(\bx_*,\bx)$ over $\bx$ on each $\mathcal{V}_m$, 
with $\bx_*$ being fixed, to obtain a new volume-averaged covariance kernel vector $\btt_{*}$.
In componentwise form,
\beq
[\btt_{*}]_m
= \frac{1}{\Delta_{m}}\int_{\mathcal{V}_{m}} q(\bx) d\mathcal{V}(\bx).
\eeq
\eit

We take a computational benefit of using the SE kernel in \cref{eq:SE}:
their exact integrations are readily available, in which the procedures become simplified
since the SE kernel can be split in a product per dimension.
The result is to obtain two integral-averaged covariance kernels, $\bfC$ and $\btt_*$.
The explicit forms in 1D are obtained in~\cite{reyes2018new}, 
which are straightforwardly extended to a general 3D case in this study,
%
\begin{equation}
  \label{eq:SE-pred}
  [\btt_{*}]_m =\prod_{d=x,y,z} \sqrt{\frac{\pi}{2}}\frac{\ell}{\Delta_d} \left \{
        \erf{\frac{\Delta_{m*,d}+1/2}{\sqrt{2}\ell/\Delta_d}}
     - \erf{\frac{\Delta_{m*,d}-1/2}{\sqrt{2}\ell/\Delta_d}}
     \right \}, 
\end{equation}
and
\begin{align}
  \label{eq:SE-cov}
  [\mathbf{C}]_{mn} = 
  \prod_{d=x,y,z} \sqrt{\pi}\left ( \ \frac{\ell}{\Delta_d} \right)^2 & 
  \left \{ \left (
     \frac{\Delta_{mn,d}+1}{\sqrt{2}\ell/\Delta_d}\erf{\frac{\Delta_{mn,d}+1}{\sqrt{2}\ell/\Delta_d}}
  + \frac{\Delta_{mn,d}-1}{\sqrt{2}\ell/\Delta_d}\erf{\frac{\Delta_{mn,d}-1}{\sqrt{2}\ell/\Delta}}
   \right ) \right .
   \nonumber \\
  +  \frac{1}{\sqrt{\pi}} & \left . 
      \left (
              \expo{-\frac{(\Delta_{mn,d}+1)^2}{2(\ell/\Delta_d)^2}} + 
              \expo{-\frac{(\Delta_{mn,d}-1)^2}{2(\ell/\Delta_d)^2}} \right
      ) \right . \nonumber \\ 
      -2 & \left . 
      \left (
               \frac{\Delta_{mn,d}}{\sqrt{2}\ell/\Delta_d}\erf{\frac{\Delta_{mn,d}}{\sqrt{2}\ell/\Delta_d}}
            + \frac{1}{\sqrt{\pi}}\expo{-\frac{\Delta_{mn,d}^2}{2(\ell/\Delta_d)^2}}\right ) 
       \right \},
\end{align}
where
\begin{eqnarray}
    \Delta_{mn,d} = \frac{\be_d \cdot (\bx_{n}-\bx_{m})}{\Delta_d}, \;\;
    {\bx_m} =  (x_{i_m}, y_{j_m}, z_{k_m})^T, 
    {\bx_n} =  (x_{i_n},y_{j_n},z_{k_n} )^T, \;\;
\label{eq:delta_mnd}
\end{eqnarray}
with $\be_d$ and $\Delta_d$ being the unit vector 
and the grid delta (e.g., $\dx, \dy, \dz$) in each $d$-direction,
respectively.

Finally, we obtain a new GP \textit{reconstructor} for FV by
taking the volume-averaged quantities as input and computes a $(2R+1)$th accurate
pointwise function value $\tilde{m}$ at $\bx_*$ 
by way of designing the integral version of the posterior mean function, 
which reads as
\begin{equation}\label{eq:gp_reconstructor}
\tilde{m}_* =
\tilde{m}({\bx_*})= 
\mathbf{t}_*^{T}\mathbf{C}^{-1}\overline{\bq}_{ij}
=\mathbf{z}_*^{T}\overline{\bq}_{ij},
\end{equation}
where $\mathbf{z}_*^{T} = \mathbf{t}_*^{T}\mathbf{C}^{-1}$ is called the prediction vector.
To ensure the interpolation of constants function by GP, i.e., $\mathbf{z}_*^{T}\bm{1}=1$,
we further normalize the prediction vector by its $L_1$ norm,
\begin{equation}\label{eq:normalized_z}
    \mathbf{z}_*^{T} \longrightarrow {\mathbf{z}_*^{T}}/{||\mathbf{z}_*^{T}||_1}.
\end{equation}
As noted in~\cite{reyes2018new,reyes2019variable,reeves2020application}, 
the prediction vector $\bz_*$ is data-independent, only depending on the grid configuration. 
Therefore, in practice, $\mathbf{z}_*$ can (and should) be pre-computed before each simulation 
as soon as the grid geometry is configured
\footnote[$\ddagger$]{When an adaptive mesh refinement (AMR) is considered for a simulation,
$\mathbf{z}_*$ can be computed for each different AMR level and they can be saved for reuse
during the simulation.
}.

We note the close resemblance between 
\cref{eq:posterior_mean} and \cref{eq:gp_reconstructor}.
They  both are defined as a linear system whose size is characterized by
the sizes of the $N\times N$ square covariance kernel matrix 
and two one-dimensional vectors of
the input data and the covariance kernel vector that are $N \times 1$.
At first glance, the total number of operations seems 
to be $\mathcal{O}(N^3)$, including
a vector-matrix multiplication of  $\mathcal{O}(N^2)$ 
to get the prediction vector $\mathbf{z}_*^T$,
followed by the $\mathcal{O}(N)$ dot-product calculation
between  $\mathbf{z}_*^T$ and the input vector.
In practice, however, the overall operation count in \cref{eq:posterior_mean} 
and \cref{eq:gp_reconstructor} is much lower 
since the calculation of the prediction vector, $\bz_*^T$, is  pre-computed only
once and for all for the entire domain. 
$\bz_*^T$ is then saved and reused throughout the simulation,
requiring only the $\mathcal{O}(N)$ dot-product calculation
between  $\mathbf{z}_*^T$ and the input vector during the
simulation.

Besides, we use Cholesky decomposition 
to invert the $N\times N$ covariance
kernel matrices, $\bK$ and $\bfC$, as they both are symmetric positive definite,
which further reduces the
computational load by half compared to the other general non-symmetric
solvers (e.g., Gaussian elimination, LU decomposition).
The consequence is the GP interpolator in \cref{eq:posterior_mean} 
or the GP reconstructor in \cref{eq:gp_reconstructor},
whose computational expense is simply governed by the
$\mathcal{O}(N)$ dot-product calculation per cell 
between the constant prediction vector and 
the time-space-varying input vector.

In the case of a 2D regular Cartesian mesh using a $q$-point Gaussian quadrature rule, 
we need $4q$ prediction vectors ($q$ per cell-face). 
An example of values of $\Delta_{mn,d}$ for each $d=x,y$ can be found in \ref{apdx:delta_kh} 
for the GP radius $R=1$ stencil and the two-point Gaussian quadrature rule.

%

\subsection{Singularity of the covariance kernel matrix}
\label{sec:gp_singularity}
The covariance kernels $\bfK$ and $\bfC$ in~\cref{eq:posterior_mean,eq:gp_reconstructor} become
nearly singular when the covariance kernel flattens out in the limit of
$(\bx - \by)^2/2\ell^2 \to 0.$ This can happen in the following cases: 
(i) the hyperparameter $\ell \to \infty$, or
(ii) the computational grid is progressively refined, approaching $\Delta_d \to 0$.
As reported in~\cite{reyes2018new}, this issue is persistent at high resolutions and negatively impacts 
the quality of GP's reconstruction. The outcome is manifested by non-convergent solution behaviors
in all forms of grid convergence studies when the covariance kernel's condition number 
reaches a sufficiently large value,
e.g., $\kappa \gtrsim 10^8.$ 

There are two workarounds to resolve this computational issue.
The first approach is found in the radial basis function (RBF) community,
where practitioners have suggested 
alternative computational methods that help
effectively avoid the ill-conditioning issue
~\cite{fornberg2004stable,fornberg2011stable,
fornberg2008stable,fasshauer2012stable,fornberg2013stable}.
One of the well-known methods is the Contour-Pad\'{e} type algorithm
~\cite{wright2003radial,fornberg2004stable}, originally referred as RBF-CP,
which has been proposed to improve RBF approximations,
where the same ill-conditioning issue arises in the limit of flat basis functions.
The noble idea in RBF-CP is to interpret
the target RBF interpolant at a finite number of evaluation points as
a complex vector-valued function of the length scale parameter,
called the shape-parameter, $\epsilon$
\footnote[$\sharp$]{The shape parameter $\epsilon$ 
of RBF's Gaussian radial kernel -- the equivalence of our SE kernel --
is inversely related to the hyperparameter $\ell$ for SE, e.g.,
$\epsilon = 1/\ell.$
}.
Putting in the context of GP, RBF-CP is
equivalent to considering the prediction vector $\bfz_*$ a complex function of $\ell$
in the complex $\ell$-plane by way of considering a contour path,
from which a vector-valued Pad\'{e} rational approximation
of $\bfz_*$ is derived and is used as a proxy for computing 
the function evaluation at $\bx_*$ stably in the limit of
Im$(\ell) = 0$ and Re$(\ell)\to \infty$.
A recent study~\cite{wright2017stable} 
extends the original RBF-CP~\cite{wright2003radial,fornberg2004stable} 
to a new RBF-RA (RA for rational approximation) that has higher accuracy
for the same computational cost, is simpler in code implementation,
and is more robust for computing the poles of the rational approximation.
%
%
%

The second approach is the use of higher precision floating-point values, 
e.g., quadruple precision. 
Practically speaking, this approach is much more straightforward than the first approach with
the Contour-Pad\'{e} algorithm, but it comes at the price of expensive precision handling
since the approach involves compiling computer codes with quadruple precision.
Fortunately in GP, there are only a couple of routines in need of 
quadruple precision, and there is no negative impact as a result (see below); 
hence the second approach is our choice in the current study.
They are the routines
that correspond to the calculation of the prediction vector, $\bfz^T_*$, 
in~\cref{eq:gp_reconstructor}, including the calculations 
of $\bfC$, the Cholesky decomposition to compute its inverse $\bfC^{-1}$, 
$\btt^{T}_*$, and finally the
multiplications of $\btt^{T}_*$ and  $\bfC^{-1}$.
These routines are compiled with quadruple-precision 
\textit{once-and-for-all} and \textit{saved} at the beginning of each simulation 
since they are independent of local fluid data.
They only depend on the grid configuration
of the GP stencil for each $R$
and the grid distance to $\bx_*$ from each grid location,
but nothing else.
%
%
The computed results with quadruple-precision are stored and saved in
double-precision arrays (therefore truncated to double-precision)
once-and-for-all. For the rest of the simulation, they are re-used 
for the multiplication
with the volume-averaged data, $\overline{\bq}_{ij}$, in~\cref{eq:gp_reconstructor}.
For the scope of our study, the added computational cost with quadruple-precision
is found to be very minimal~\cite{reyes2018new,reyes2019variable} since
the GP kernels, $\btt^{T}_*$ and  $\bfC^{-1}$, are configured only once
initially and multiplied to get $\bz_*^T$. What remains during the simulation 
is the $\mathcal{O}(2R+1)$ dot-product operations between the saved 
prediction vector $\bz_*^T$ and the local flow data vector $\overline{\bq}_{ij}^n$
in \cref{eq:gp_reconstructor}. 
The extra precision handling barely impacts the overall performance. 
%
%
As such, we continue following our previous work~\cite{reyes2018new,reyes2019variable}
and take the second approach in the present study.  
We state that, in another ongoing study, we have successfully implemented and
tested the Contour-Pad\'{e} algorithm for GP, the work of which will be
reported in a forthcoming study.

\subsection{The full set of spatial reconstruction methods}\label{sec:full_set}
To this point, we have designed a set of three high-order GP reconstruction methods,
the 3rd-order GP-R1, the 5th-order GP-R2, and the 7th-order GP-R3.
They are all genuinely multidimensional, unlimited, and their order of accuracy varies
by choosing one of the three GP stencil configurations described in~\cref{sec:gp_stencil}
depending on the GP radius, $R=1, 2, 3$.

By construction, our GP methods deliver meaningfully 
high-order solutions (i.e., third-order or higher) with $R \ge 1$; 
otherwise, the GP regression
with $R=0$ becomes a constant function, $\bK = K(\bx,\bx) = 1$, which is equivalent to
the classic first-order Godunov (FOG) method~\cite{godunov1959difference}.
In this case, we simply switch to the FOG reconstruction that
directly copies 
the cell-centered conservative variables to the Riemann states at cell-faces
without any further consideration of high-order spatial approximations such as our GP
reconstructions.
As known, the FOG method is very robust and is a positivity-preserving method,
always guaranteeing the physically admissible conditions without unphysical oscillations
at shocks and discontinuities.

This completes our discussion on different reconstruction methods, consisting
of FOG, GP-R1, GP-R2, and GP-R3, where the order of accuracy ranges from
as low as the 1st-order to as high as the 7th-order in the odd integer sequence.
In preparation for designing our GP-MOOD algorithm,
we group these methods into three different classes: 
(i) GP-MOOD3 (FOG and GP-R1), 
(ii) GP-MOOD5 (FOG, GP-R1, and GP-R2),
and (iii) GP-MOOD7 (FOG, GP-R1, and GP-R3).

For the sake of comparison studies of GP, 
we introduce the 4th group, labeled as (iv) POL-MOOD3 (FOG and Poly3),  
which uses two polynomial-based
methods only, including FOG and the unlimited 
3rd-order polynomial reconstruction, Poly3, on the
five-point stencil in~\cref{fig:gp_stencil_GP-R1} defined by,
\beq\label{eq:poly3}
p_3(x,y) = a_0 + a_1 x + a_2 x^2 + a_3 y + a_4 y^2 + a_5 xy.
\eeq
The coefficients are determined by matching the volume averages, $\overline{\bq}_{i_m j_m}$,
on each of the five cells in the five-point stencil, giving rise to the following $5 \times 5$ system,
\beq\label{eq:p3_coeffs}
\frac{1}{\dx \dy} \int_{I_{j_m}} \int_{I_{i_m}} p_3(x,y) dx dy = \overline{\bq}_{i_m j_m}, \;\;\; 1 \le i_mj_m \le 5.
\eeq
At first glance, there seem to be more unknowns, $a_0, \dots, a_5$, than knowns,
$\overline{\bq}_{i_m j_m}$, $1 \le i_mj_m \le 5$.
However, the last term with $a_5$ always cancels out on Cartesian grids,
leaving the remaining five coefficients to be uniquely determined using the five knowns.
%
%
The solution to the system \cref{eq:p3_coeffs} provides the five coefficients
(see \ref{apdx:p3_coeffs}),
from which
the Riemann states are obtained directly by 
evaluating $p_3(x,y)$ at the 4th-order Gaussian quadrature points, 
e.g., $g_1$ and $g_2$ from the left panel of \cref{fig:multipt-QR}. 
Using the indexing from \cref{fig:gp_stencil_GP-R2-R3}, we have,
\begin{eqnarray} 
p_3(x_{g_1}, y_{g_1})=  \frac{5}{6}\overline{\bq}_{1}   -\frac{1}{6}\overline{\bq}_{2} +  \frac{1}{4\sqrt{3}}\overline{\bq}_{3} +  \frac{1}{3}\overline{\bq}_{4}   -\frac{1}{4\sqrt{3}}\overline{\bq}_{5},\\
p_3(x_{g_2}, y_{g_2})=  \frac{5}{6}\overline{\bq}_{1}   -\frac{1}{6}\overline{\bq}_{2}   -\frac{1}{4\sqrt{3}}\overline{\bq}_{3} +  \frac{1}{3}\overline{\bq}_{4} +  \frac{1}{4\sqrt{3}}\overline{\bq}_{5}.
\end{eqnarray}
The Riemann states on the other cell faces (top, bottom, and left) are computed 
in a similar way by rotating the above expressions correspondingly.

\section{Integrating GP into the MOOD framework}\label{sec:gp_into_mood}
The GP reconstruction methods are now ready to be integrated into the MOOD framework. 
As briefly introduced in~\cref{sec:intro}, the main idea in the MOOD method is the \textit {a posteriori}
limiting strategy~\cite{clain2011high,diot2012improved,diot2013multidimensional,diot2012methode},
which updates each cell with the highest accurate solver available first, followed by 
the cell-by-cell inspection to see if a set of MOOD admissibility
conditions are met locally. For example,
suppose the updated solution at  $\bx_{ij}$ after the first pass 
with the highest accurate solver fails to meet the admissibility constraints. 
In that case, the process is repeated until the constraints are met with the next highest accurate solver.
In the worst case, a local solution could end up with 
the most diffusive -- but most stable -- solver, e.g., FOG, 
in the regions where shocks
and discontinuities are present. 
Reportedly, and also will be seen in our results in~\cref{sec:results},
the regions of such troubled cells are 
only about a few percent (e.g., less than 10\% in practice) of the entire domain
~\cite{diot2013multidimensional,diot2012methode,bourriaud2020priori}.

\begin{figure}[ht!]
    \centering
	\centering
        \includegraphics[width=11cm]{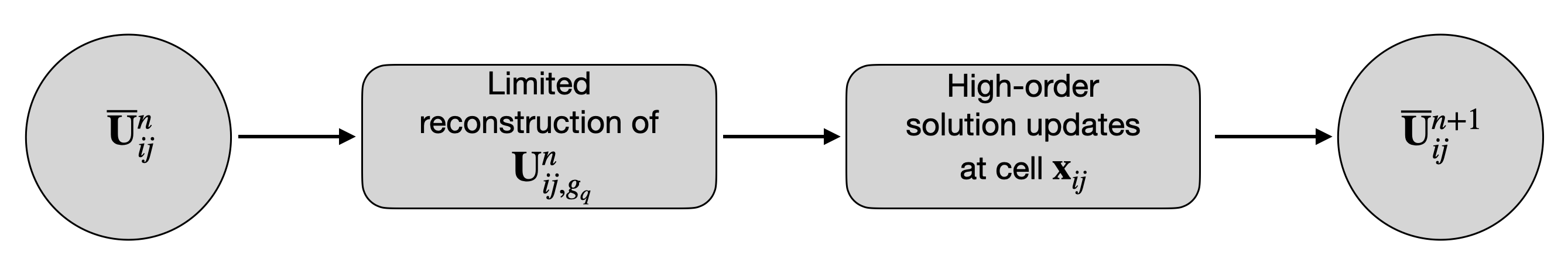}
    \caption{The logical flow line of the solution updating procedure.
    Shown is the principle flow line in conventional \textit{a priori} high-order methods 
    where limited spatial data reconstructions are applied to \textit{all} available cells
    as a fundamental building-block, regardless of local flow smoothness.}
    \label{fig:a_priori_flowChart}
\end{figure}

The MOOD method, by design, is endowed with the positivity-preserving property of FOG
near sharp flow gradients while utilizing high-order solutions away from
the cells that experience gradients build-ups. This paradigm of MOOD's \textit {a posteriori}
limiting is conceptually different from the conventional \textit {a priori} limiting strategies
that depend on computationally expensive 
nonlinear shock/discontinuity controls required on every single cell in the
simulation. \cref{fig:a_priori_flowChart} 
displays the logical pipeline
in conventional \textit {a priori} shock-capturing FV methods, where
nonlinear controls play an essential role in a stable evolution of discrete solutions.
The nonlinear nature of their operations inevitably increases computational intensity.

For this reason, the MOOD method is considered to be an alternative paradigm that has 
computational benefits over the classic \textit {a priori} compressible flow algorithms.
This paper explores this attractive MOOD paradigm by integrating the high-order GP
reconstruction methods in the MOOD framework. 
We emphasize that the resulting GP-MOOD method seamlessly provides the pairwise 
integration between the GP reconstruction methods and
MOOD's order cascading mechanism.
The primary advantage of our GP-MOOD method is the performance gain
and the paradigm simplification
by removing the need for least-squares solves in the existing polynomial-based MOOD methods.
Also will be shown are the smaller GP stencils compared 
to the corresponding polynomial MOOD
methods at the same solution accuracy,
further decreasing the computational workload and 
the associated memory footprint in GP-MOOD
relative to the existing polynomial-based MOOD methods.

\subsection{Three building blocks in GP-MOOD}\label{sec:three_bldg-blocks}

The fundamental principle of the MOOD paradigm centers around three important
building blocks, the detection criteria, the safe scheme, and the scheme cascade
\cite{semplice2018adaptive}.
We adopt these three basic components
to devise the following relaxed version for our GP-MOOD methods:

\begin{enumerate}
    \item[(i)] \textit{The detection criteria:}
    The first component is a sequence of prescribed properties 
    that the discrete numerical solution has to fulfill to be considered 
    acceptable. These conditions are of two types: ``Physical Admissibility Detection'' (PAD) and 
    ``Numerical Admissible Detection'' (NAD). In most fluid dynamics simulations, 
    PAD ensures that the numerical solution represents an admissible flow state
    (e.g., positivity in pressure and density). More generally, they ensure that the numerical solution 
    makes sense from the physical model's point of view. They only depend on the set of PDE solved. 
    
    On the other hand, NAD ensures that the solution produced by the solver is essentially non-oscillatory. 
    It is based on a relaxed discrete maximum principle (DMP). 
    It also includes the detection of non-numeric values such as 
    \texttt{NAN}'s and \texttt{Inf}'s, that is, the admissibility of the state from a computer science point of view 
    (Computer science Admissibility Detection or CAD). CAD is identical for all sets of PDEs. 
    If a candidate solution does not satisfy either of the PAD and NAD criteria in some cells, such cells are
    recorded as \textit{troubled cells} and their discrete updates are repeated with a lower order method (see (iii) below).
    
    Newly added to the conventional detection criteria is the so-called 
    Compressibility-Shock Detection (CSD) check in our study. 
    Operationally, this new condition is run after PAD and CAD (but not DMP) to measure
    the strengths of fluid compressibility and local shocks
    by measuring the local gradient of pressure ($\nabla p$) and the divergence of the velocity field
    ($\nabla \cdot \bV$).
    We proceed to the next step of DMP
    only when cells undergo strong compressibility with rapid build-ups of pressure gradients.
    The CSD check is motivated by the early work on 1D GP reported in~\cite{lee2017new}, where
    GP's unlimited reconstructed solutions are used away from shocks.
    In contrast to the current GP-MOOD methods,
    the second-order limited piecewise linear method (PLM) was used 
    where shocks are present in~\cite{lee2017new}.
    In this early work, local shocks were tracked by 
    a generic shock detector~\cite{balsara1999staggered} to selectively choose 
    GP versus PLM depending on the local shock strengths.
    
    To summarize, we have the following three categories:
    \bit
    \item[a)] PAD: positivity preservation on density and pressure variables, 
    \item[b)] NAD: numerical validity that monitors CAD (\texttt{NAN \& Inf}) and DMP (see more \cref{sec:DMP}),
    \item[c)] CSD: compressibility and shock strengths, i.e., $\nabla p > \sigma_p$ and $\nabla \cdot \bV < -\sigma_v$,
    where $\sigma_p>0$ and $\sigma_v>0$ are (heuristically) tunable threshold parameters.
    \eit
    \item[(ii)] \textit{The safe scheme:}
    The second component is the choice of a numerical method used as the last resort 
    when all the other high order schemes have failed to produce an acceptable solution according to the detection criteria in (i). 
    Therefore, the selection choice has to focus on scheme's robustness and stability that guarantee to 
    produce an admissible solution state. To this end, the first-order Godunov (FOG) scheme is most popular, 
    while the second-order MUSCL method could be used as well to improve the results 
    on contact discontinuity (e.g., see \cite{padioleau:tel-03130146}).
    In this study, we use FOG for the choice of the safe scheme.
    \item[(iii)] \textit{The scheme cascade:}
    A family of reconstruction schemes is the third component that provides a sequential pipeline of  
    different reconstruction methods, starting from the most accurate available method to the safe scheme.
    The conventional MOOD method uses a set of unlimited \textit{polynomial} reconstruction methods
    in different orders up to the 6th-order accuracy~\cite{diot2012improved,diot2013multidimensional}.
    Alternatively, for the present study, we use the three unlimited linear GP reconstruction methods 
    of 3rd-, 5th-, and 7th-order, namely, GP-R1, GP-R2, and GP-R3.
    As briefly introduced in \cref{sec:full_set}, 
    these GP methods and the simple 3rd-order polynomial MOOD method 
    are used in this study.
    
    To summarize, we has the following four families of methods, each of which
    have the specific ordered cascading scheme depicted by the arrows on troubled cells:
    \bit
    \item[a)] the 7th-order GP-MOOD7: the cascading follows as GP-R3 $\to$ GP-R1 $\to$ FOG,
    \item[b)] the 5th-order GP-MOOD5: the cascading follows as GP-R2 $\to$ GP-R1 $\to$ FOG,
    \item[c)] the 3rd-order GP-MOOD3: the cascading follows as GP-R1 $\to$ FOG,
    \item[d)] the 3rd-order POL-MOOD3: the cascading follows as Poly3 $\to$ FOG.
    \eit    
    The decrement pattern of our order-cascading within each group
    is consistent with the finding in~\cite{diot2013multidimensional}, 
    which simplifies the original, long-listed one-by-one polynomial degree decrement procedure
    (e.g., Poly 5 $\to$ Poly4 $\to$ Poly3 $\to$ Poly2 $\to$ FOG)
    to a shorter list of two or three methods.
 \end{enumerate}

\begin{figure}[ht!]
    \centering
    \centering
    \begin{subfigure}{14cm}
        \centering
        \includegraphics[width=\textwidth]{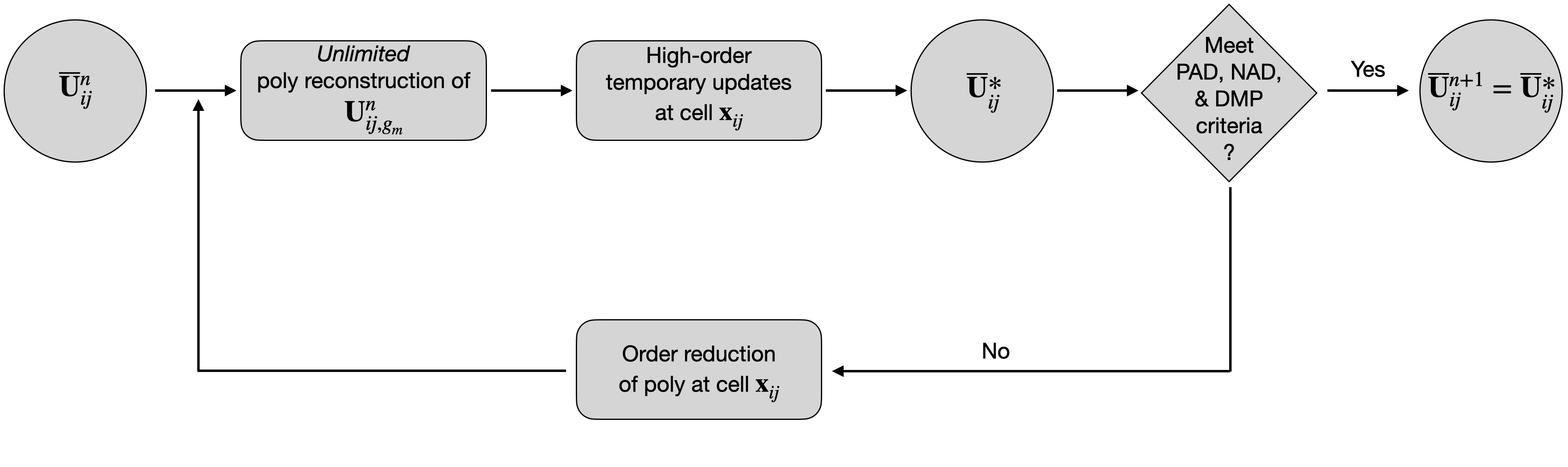}
        \caption{}\label{fig:poly_mood_flowChart}
    \end{subfigure}
    \begin{subfigure}{16cm}
        \vspace{0.5cm}
	\centering
        \includegraphics[width=\textwidth]{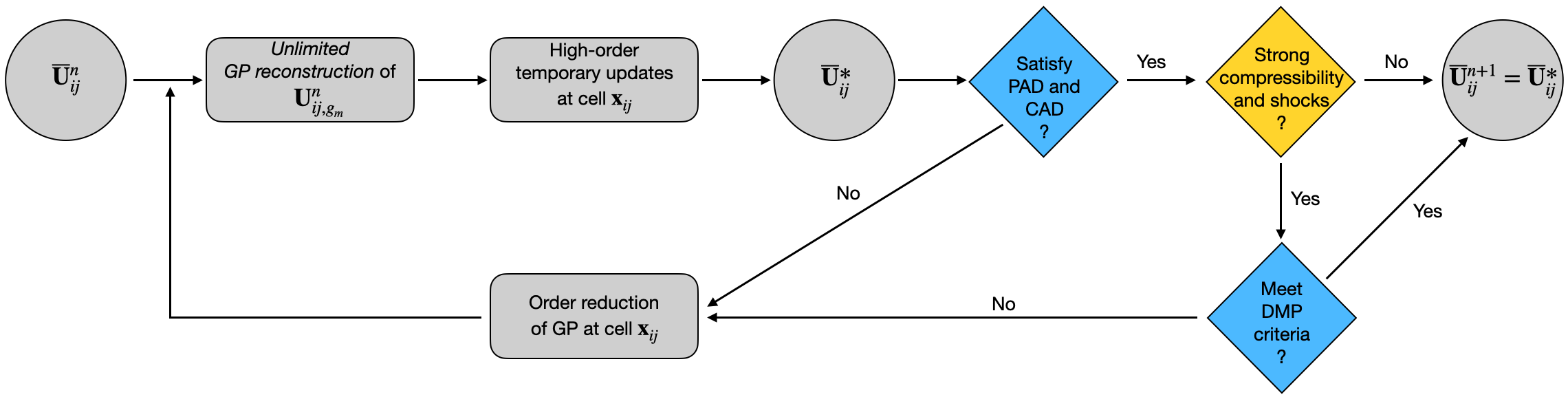}
        \caption{} \label{fig:gp_mood_flowChart}
    \end{subfigure}    
    \caption{The logical flow of the solution updating procedure in the MOOD loop.
    (a) The flow chart in the existing \textit{a posteriori} polynomial MOOD method.
    The solution accuracy of a group of unlimited polynomial data reconstruction methods
    cascades down from high to low until all the MOOD criteria
    are met in each troubled cell.
     (b) The flow chart of the GP-MOOD method. 
     An extra CSD condition that checks the strength of
     flow compressibility (the yellow diamond)
    is added between the positivity and \texttt{NAN \& Inf} check 
    (i.e., PAD and CAD in the top sky blue diamond) 
    and the rest MOOD criteria on DMP (the bottom sky blue diamond).
    The conditions in the grey diamond in \cref{fig:poly_mood_flowChart}
    are split into the two sky blue diamonds 
    in \cref{fig:gp_mood_flowChart} for GP-MOOD.
    The overall order-reduction MOOD concept is relaxed 
    with the compressibility-shock admissibility (CDS) check
    to maximize the use of GP's unlimited high-order reconstructed solutions
    before they are reduced to lower order ones.
    In both, $\bU_{ij,g_m}^n$ denotes the pointwise Riemann states
    at each Gaussian quadrature point, $g_m$, $1\le m \le q$, in \cref{fig:multipt-QR};
    $\overline{\bU}^{\ast}_{ij}$ denotes a volume-averaged, pre-validated
    candidate solution using
    the highest-order method available in each MOOD cascading loop.
    }
    \label{fig:a_posteriori_flowChart}
\end{figure}

The logical MOOD loop pipelines of the standard polynomial and the GP-MOOD
methods are summarized and compared in~\cref{fig:a_posteriori_flowChart}.
The primary difference between the polynomial-MOOD and GP-MOOD methods
is the split of the detection criteria 
(the grey diamond in \cref{fig:poly_mood_flowChart})
of the polynomial-MOOD method
into two groups, placing PAD and CAD in the first group and 
the rest DMP in the second group
(the two sky blue diamonds \cref{fig:gp_mood_flowChart}) 
in the GP-MOOD loop.
The newly added compressibility-shock detection check (CSD) is monitored
between the two separated detection criteria groups, allowing
a relaxed MOOD framework that promotes 
the use of GP's high-order solution updates
as much as possible.
 
We detail the GP-MOOD strategies described above in the next following subsections.
For this paper to be self-contained, we provide concise descriptions 
of those existing concepts while we give sufficient details on new strategies for GP-MOOD.
Interested readers refer to 
\cite{clain2011high,diot2012improved,diot2013multidimensional,diot2012methode,
semplice2018adaptive,bourriaud2020priori} for reviewing 
the existing MOOD methodologies.

\subsection{PAD -- Physical Admissibility Detection}\label{sec:PAD}
The PAD criterion is drawn from considering the most constitutional condition inferred from the underlying physics.
As we solve the Euler equations, it is indispensable to ensure positivity in both density and pressure
variables to guarantee stable discrete evolutions of numerical solutions regardless of how
extreme the flow in simulations may turn into. 
The candidate solution $\overline{\bq}^*_{ij}$ on each $\bx_{ij}$ satisfies the PAD criterion if:
\begin{equation}\label{eq:PAD}
    \overline{\rho}^*_{ij} > 0 \quad \text{and} \quad 
    \overline{p}_{ij}^* > 0.
\end{equation}
For those cells that fail to meet this positivity condition, the solution order of accuracy
gets reduced according to the decrement patterns, a) through d), described in 
(iii) \textit{The scheme cascade} in the previous
section. The logical flow repeats the MOOD loop for the next iteration.

It is important to re-emphasize that obtaining the derived primitive variable,
$\overline{p}_{ij}^*$, from the
updated conservative quantities, $\overline{\bq}^*_{ij}$, involves
nonlinear conversion processes such as an EoS call,
which can have an
impact on the solution accuracy \textit{if} $\overline{p}_{ij}^*$ is
subsequently re-used to construct any relevant conservative variables
without respecting the difference between the pointwise and volume-averaged
quantities.
As clearly noted in \cref{sec:governing_eqns}, such a conversion 
is strictly prohibited. However, there is no negative impact
on the solution accuracy as long as the derived pressure
is used for PAD and is no longer used afterward, 
which is the case here.

\subsection{NAD -- Numerical Admissibility Detection}\label{sec:NAD}
Below, we discuss two types of criteria that are monitored in NAD, 
including the CAD and the relaxed DMP criteria.

\subsubsection{CAD -- Computer science Admissibility Detection}\label{sec:CAD}
The CAD criterion
\footnote[$\P$]{Our numerical experiments show that it is sufficient to check
the CAD criterion on density and pressure only, although one can check CAD on all 
updated variables. Our selective choice is consistent with the single variable choice
in \cite{semplice2018adaptive}, which uses the fluid density to avoid
the unnecessary detection processes.
We use the \texttt{ISNAN} command for CAD 
in our Fortran implementation.}
ensures that the updated candidate solutions
$\overline{\bq}^*_{ij}$ do not represent either \texttt{NAN}s or \texttt{Inf}s, which are
direct outcomes of invalid floating-point operations, such as
division by zero, 
square root of a negative number, 
arithmetics with $\pm \infty$, etc.
In most hydrodynamics codes that use \textit{a priori} limitings, 
these invalid states are considered \textit{physically violated} 
in a strong sense 
because the process of \textit{a priori} limitings should explicitly
prohibit such states. As such, these codes are often
compiled with ``crash-if-\texttt{NAN/Inf}'' flag options.

The situation is different in \textit{a posteriori} MOOD paradigms, where
there is no such an explicit monitoring process until a candidate solution is
available, only after which invalid solutions are rejected 
via MOOD's post processes.
Another view to see this is that the Riemann states 
computed by the unlimited MOOD reconstruction
are not necessarily physically valid. Indeed, those invalid
states are well justified (and should be allowed) 
as \textit{proper outputs} within the MOOD paradigm.
Thus, ``crash-if-\texttt{NAN/Inf}'' compiler flags are not an option
for the \textit{a posteriori} MOOD-type codes. Instead,
these invalid states are separately controlled via the CAD process, i.e.,
we say the candidate solution $\overline{\bq}^*_{ij}$ satisfies the CAD criterion if:
\begin{equation}\label{eq:CAD}
    \overline{\rho}^*_{ij} \ne (\texttt{NAN} \mbox{ or } \texttt{Inf}) \quad \text{and} \quad 
    \overline{p}_{ij}^*      \ne (\texttt{NAN} \mbox{ or } \texttt{Inf}).
\end{equation}

We also make an important note on CAD in connection to 
preserving the assumed symmetry in symmetry-preserving simulations.
Our experiments have shown that a broken symmetry can be
induced if any one of the Riemann state input pair,
$(U^n_L, U^n_R)$ at each Gaussian quadrature point
happens to be an invalid state. For example, assume that
the left state vector $U^n_L$ contains several state variables
that are $\texttt{NAN}$ but a valid normal velocity.
If the normal velocity is negative, 
and all the state variables
in $U^n_R$ are physically valid with a negative normal velocity as well, 
the Riemann problem will
still be able to calculate the interface flux based on the valid right
states based on the upwinding nature of the Riemann problem.
Since this upwind flux has been computed
with such a physically non-admissible input pair, 
the resulting flux has no guarantee to produce a physically
admissible candidate solution, albeit it qualifies all of the MOOD criteria.
If this type of incident continues in simulations, 
this inconsistency can trigger the onset of unphysical asymmetries,
which could lead to the development of erroneous flow instabilities.

To prevent this issue from happening, we check the CAD condition
also on the reconstructed densities (and nonlinearly derived pressures) in
all Riemann states from the high-order GP reconstruction,
in addition to the updated candidate solutions at cell centers.
To our knowledge, the previous studies on the MOOD methods
\cite{clain2011high,diot2012improved,diot2013multidimensional,diot2012methode}
are mostly intact from this asymmetry issue, primarily because their
applications heavily focus on unstructured grid geometries where
symmetry is not rigorously expected in the first place.

\subsubsection{Relaxed DMP criteria: plateau + DMP + \textit{u2}}
\label{sec:DMP}
The discrete maximum principle (DMP) is the crux of the MOOD method, which
was single-handed to detect the solution admissibility in the original MOOD method
\cite{clain2011high}. For the exposition purposes, 
let us assume that we use the density variable to run the DMP check.
The original DMP monitors if there is any excessive numerical oscillation produced
in the resulting candidate solution $\overline{\rho}^{*}_{ij}$ at $\bx_{ij}$ in comparison to
the adjacent neighboring input states, $\overline{\rho}^{n}_{m}$, 
\begin{equation}\label{eq:DMP}
    \underset{m \in \mathcal{N}(ij)}{\min}\Bigl(\overline{\rho}_m^{n},\overline{\rho}_{ij}^{n}\Bigr) 
    \leq \overline{\rho}_{ij}^* \leq  
    \underset{m \in \mathcal{N}(ij)} {\max}\Bigl(\overline{\rho}_m^{n},\overline{\rho}_{ij}^{n}\Bigr),
\end{equation}
where $\mathcal{N}(ij)$ is a set of all indices, including the immediate neighbors
that share a common cell-face with the cell $\bx_{ij}$.
For example, in our structured 2D Cartesian grid configuration
for GP-MOOD3
depicted in~\cref{fig:gp_stencil_GP-R1},
the set $\Bigl\{\overline{\rho}^n_m : m\in\mathcal{N}(ij)\Bigr\}$ is given as
$\Bigl\{\overline{\rho}^n_m : m=2, 3, 4, 5\Bigr\}$.
This concept is extended to the following DMP check
for a $p$-stage RK method,
\begin{equation}\label{eq:DMP_RK}
    \underset{m \in \mathcal{N}(ij)}{\min}\Bigl(\overline{\rho}_m^{n, rk(s)},\overline{\rho}_{ij}^{n,rk(s)}\Bigr) 
    \leq \overline{\rho}_{ij}^{*,rk(s)} \leq  
    \underset{m \in \mathcal{N}(ij)} {\max}\Bigl(\overline{\rho}_m^{n,rk(s)},\overline{\rho}_{ij}^{n,rk(s)}\Bigr),
    \quad \text{for each } 1\le s \le p,
\end{equation}
where $\overline{\rho}_m^{n, rk(s)}$ and $\overline{\rho}_{ij}^{*,rk(s)}$ are respectively
the adjacent initial states and the resulting candidate solution at each stage $s$
over the course of the $p$-stage RK updates.

In the follow-up studies
\cite{diot2012improved,diot2013multidimensional,diot2012methode,semplice2018adaptive},
this single-handed DMP criterion has been further revised and relaxed by two additional
detections, namely the Plateau detection and the  \textit{u2} detection.
Hereunder, we describe them briefly in the order they are used during the NAD procedure.
Following~\cite{semplice2018adaptive}, the Plateau detection is executed \textit{before}
the original DMP in \cref{eq:DMP_RK}
in order to avoid the loss of precision on constant flat plateau states
as a consequence of the criterion in \cref{eq:DMP_RK} 
that detects not only large unphysical oscillations but also
micro-oscillations. 
To improve such situations, the first relaxation takes place
in the form of skipping the DMP criterion in \cref{eq:DMP_RK}  if
micro-oscillations (hence the name Plateau detection) are present:
\begin{equation}\label{eq:plateau}
\underset{m \in \mathcal{N}(ij)}{\max}\Bigl(\overline{\rho}_m^{n,rk(s)},\overline{\rho}_{ij}^{n,rk(s)}\Bigr) -
\underset{m \in \mathcal{N}(ij)}{\min}\Bigl(\overline{\rho}_m^{n, rk(s)},\overline{\rho}_{ij}^{n,rk(s)}\Bigr) 
< \min(\dx^3,\dy^3).
\end{equation}

The second relaxation called the \textit{u2} detection
aligns with overcoming the second-order-accuracy 
bottleneck if the MOOD loop strictly satisfies the original DMP in \cref{eq:DMP_RK}
\cite{diot2013multidimensional}.
A resolution is made available by allowing the violation of \cref{eq:DMP_RK} on smooth extrema
via the \textit{u2} detection.
Executed \textit{after} \cref{eq:DMP_RK}, the  \textit{u2} detection 
reinstates the admissibility of the rejected candidate solutions from \cref{eq:DMP_RK}.
The  \textit{u2} detection checks if the rejection is the consequence of the situation where 
a new extremum state $\overline{\rho}_{ij}^{*,rk(s)}$ happens to be
present in the region where local flows experience smooth variations.
If this is the case, the rejection is revoked, and its admissibility is reinstated.
To meet this, \textit{u2} calculates local second derivative quantities 
to measure  local curvatures,
%
%
%
\begin{equation}\label{eq:u2_curvature}
    \mathcal{C}_{d,ij}^{\min} = 
    \underset{m \in \mathcal{N}(ij) \cup \{ij\}}  
    {\min}\left(\frac{\partial^2 \overline{\rho}_{m}^{n,rk(s)}}{\partial{\eta^2}}\right)
    \;\; \text{, } \quad
    \mathcal{C}_{d,ij}^{\max} = 
    \underset{m \in \mathcal{N}(ij) \cup \{ij\}}
    {\max}\left(\frac{\partial^2 \overline{\rho}_{m}^{n,rk(s)}}{\partial{\eta^2}}\right),
     \quad \text{for each } 1 \le s \le p, 
\end{equation}
where $\eta = \be_d\cdot \bx$ denotes each of $x, y,$ respectively for each $d=x,y$.
With these curvatures, a candidate solution
$\overline{\bq}_{ij}^{*,rk(s)}$ 
rejected by \cref{eq:DMP_RK}
after the $s$th
RK substage update, 
nonetheless regains its admissibility if the initial states 
$\overline{\rho}_{m}^{n,rk(s)}$ at each $s$th RK stage 
fall on a smooth extrema region, checked by either of the followings:
\beq\label{eq:u2_a}
\mathcal{C}_{d,ij}^{\min}\mathcal{C}_{d,ij}^{\max} > -\delta 
\;\; \mbox{ and }\;\;
\max
\Biggl(
\Big|\mathcal{C}_{d,ij}^{\max}\Big|,\Big|\mathcal{C}_{dir,ij}^{\min}\Big|
\Biggr)<\delta,
\eeq
or
\beq\label{eq:u2_b}
\mathcal{C}_{d,ij}^{\min}\mathcal{C}_{d, ij}^{\max} > -\delta 
\;\; \mbox{ and }\;\;
\frac{\Big|\mathcal{C}_{d,ij}^{\min}\Big|}
       {\Big|\mathcal{C}_{d,ij}^{\max}\Big|} \geq \frac{1}{2}, 
\eeq
for all $d = x,y$ with $\delta = \min\{\dx, \dy\}$.
One can use the standard centered differencing for 
the second partial derivatives in \cref{eq:u2_curvature}.
To provide another application of GP,
we instead use a new second-derivative GP formula that can be derived 
by taking second derivatives 
of the GP prediction vector, $\btt_*$, in~\cref{eq:gp_reconstructor}.
We show its derivation in \ref{apdx:gp_2ndDer}.

\subsection{CSD -- Compressibility-Shock Detection}\label{sec:CDS}
As shown in~\cref{fig:gp_mood_flowChart}, 
an extra layer called CSD is added to the MOOD loop 
between the PAD/CAD checks and the relaxed DMP check
in order to early-accept a candidate solution as the final high-order solution
without DMP if the candidate solution does not reside inside a strong shock.
A multidimensional shock detection switch is employed for this purpose.
Following ~\cite{mignone2011pluto}, a candidate solution $\overline{\bq}^{*,rk(s)}_{ij}$
is recorded as valid and exits the MOOD loop 
if the local flow is weakly compressible,
\beq\label{eq:divV}
\nabla \cdot \overline{\bV}_{ij}^{n,rk(s)}
\equiv
\frac{\overline{u}_{i+1,j    }^{n,rk(s)} - \overline{u}_{i-1,j  }^{n,rk(s)}}{2\dx}+
\frac{\overline{v}_{i,    j+1}^{n,rk(s)} - \overline{v}_{i,   j-1}^{n,rk(s)}}{2\dy}
\ge - \sigma_v
\eeq
and
the local (normalized) pressure gradient is weak, 
\beq\label{eq:gradP}
\widetilde{\nabla}\overline{p}_{ij}^{n,rk(s)}
\equiv
\frac{\Big|   \overline{p}_{i+1,j    }^{n,rk(s)} -\overline{p}_{i-1,j  }^{n,rk(s)}\Big|}
{2\dx \min\{\overline{p}_{i+1,j    }^{n,rk(s)},  \overline{p}_{i-1,j  }^{n,rk(s)} \}}+
\frac{\Big|  \overline{p}_{i,    j+1}^{n,rk(s)} - \overline{p}_{i,   j-1}^{n,rk(s)}\Big|}
{2\dy \min\{\overline{p}_{i,   j+1}^{n,rk(s)},   \overline{p}_{i,   j-1}^{n,rk(s)} \}}
\le \sigma_p,
\eeq
where $\sigma_v$ and $\sigma_p$ are both positive tunable parameters. 
Similar to \cite{mignone2011pluto}, our default choice is heuristically set to 5 for both, 
which sufficiently and satisfactorily
reroutes candidate solutions to early acceptance
at each $s$th RK sub-stage according to
the two conditions in~\cref{eq:divV,eq:gradP},
promoting the use of GP's high-order solutions without further order decrement.

We have observed that this new CSD criterion is necessary to control
numerical dissipation on the right scale. Otherwise, the final validated
solutions channeled directly to the next DMP check without early acceptance 
become too diffusive to trigger subsequent nonlinear flow patterns 
that have been well-understood as unique signatures of such benchmark problems.
Examples are discussed in~\cref{sec:shu_osher,sec:Implosion}.

\subsection{The safe scheme}\label{sec:safe_scheme}
The first order Godunov (FOG) scheme is our choice for ``the safe scheme'' to guarantee
the inherent solution robustness with strong positivity preservation if all tried high-order solutions
turn out to be inadmissible during the MOOD loop. As will be seen in \cref{sec:results},
the FOG solutions are employed near shocks and discontinuities in most shock
dominant simulations, taking less than 10\% of the entire domain in practice.
Another possible choice could be to use a \textit{limited} 
2nd-order TVD linear method as a safe scheme
to improve the solution quality around contact discontinuities 
with extra care for positivity preservation
(e.g., see \cite{padioleau:tel-03130146}).
However, compared to the safe scheme with FOG,
this choice could potentially 
compromise the strong positivity preservation as well as the rest desirable
properties monitored in the MOOD loop since those properties are
no longer to be checked once the solution cascades down 
to the safe scheme but to admit it.

\subsection{The scheme cascade with GP}\label{sec:cascade}

In GP-MOOD, 
the original strategy in polynomial-based MOOD methods
with polynomial reconstruction schemes, their grouping, and
the order decrement  in the MOOD loop 
are replaced by the GP alternatives of high-order 
unlimited linear GP reconstruction with different choices of $R=1,2,3$.
Three GP groups for the MOOD cascade, introduced in \cref{sec:three_bldg-blocks},
are our GP alternatives to the polynomial-based MOOD methods. They are
the 7th-order GP-MOOD7 (GP-R3 $\to$ GP-R1 $\to$ FOG),
the 5th-order GP-MOOD5 (GP-R2 $\to$ GP-R1 $\to$ FOG), and
the 3rd-order GP-MOOD3 (GP-R1 $\to$ FOG) for GP-MOOD.
The additional 3rd-order POL-MOOD3 (Poly3 $\to$ FOG) will be used for comparison.

In~\cite{clain2011high}, three systematic strategies called \textbf{EdgePD} have been introduced.
The main idea is to provide a consistent way to assign 
a specific choice of reconstruction order of accuracy
at each cell face when the two neighboring cells
that share the cell face in common 
are at different orders or accuracy 
as a result of different paths in the MOOD order decrement.
\begin{figure}[hbt!]
    \centering
    \includegraphics[scale = 0.22]{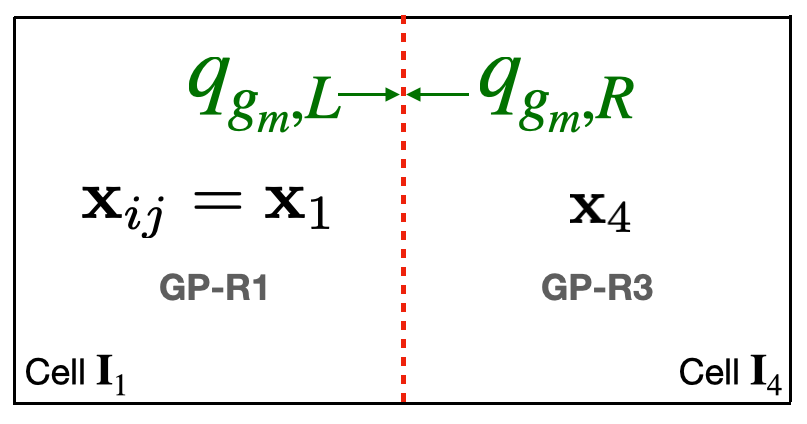}
    \caption{A GP-MOOD7 example for determining 
    reconstruction orders for two Riemann states,
    $q_{g_m,L}$ and $q_{g_m,R}$, which respectively 
    belong to the cell $\bI_1$ on the left and  $\bI_4$ on the right.
    The Riemann states are displayed at one of the 
    Gaussian quadrature points, $g_m$, $m=1, \dots, q$,
    on the cell face $\bx_{i+1/2,j} = (x_{i+1/2},y_j)$ (red dotted line)
    for the $q$-point Gaussian quadrature rule. 
    In this example with GP-MOOD7, we specifically 
    use the 4-point quadrature rule (i.e., $q=4$) that
    provides the 8th-order quadrature
    rule, sufficient for the 7th-order GP-MOOD7.}
    \label{fig:gp_EDP_1}
\end{figure}
For example, 
assume that the computation is being done with
the 7th-order GP-MOOD7, which is integrated with 
the 4-point Gaussian quadrature rule (i.e., $q=4$, see \cref{fig:multipt-QR}).
Consider a situation in~\cref{fig:gp_EDP_1}, where 
the cell face at $(x_{i+1/2}, y_j)$ is depicted in a red dotted line.
Two cells that share this face, namely $\bI_1$ and $\bI_4$, are
centered at $\bx_1 = (x_i, y_j)$ and $\bx_4 = (x_{i+1}, y_j)$,
respectively. Readers can see that 
these two cells are extracted from \cref{subfig:GP-R3}. 
What can happen is that, at any $s$th RK sub-stage,
the reconstruction order at
$\bI_4$ remains at the highest 7th-order with GP-R3
while it has been reduced to the 3rd-order with GP-R1
on $\bI_1$. There are three different ways to assign a 
reconstruction order at each of the left and right Riemann state pair 
$(q_{g_m,L}, q_{g_m,R})$
at each Gaussian quadrature point, $g_m$. See also~\cref{fig:multipt-QR}.
The first option called {$\mbox{EDP}_0$} is to assign 
the Riemann state that belongs to each of the cells
the same order each cell is at. That is, 
GP-R1 is assigned to $q_{g_m,L}$ while GP-R3 to $q_{g_m,R}$.
The second option, {$\mbox{EDP}_1$}, takes the minimum of the two, GP-R1 and GP-R3,
and assigns the minimum to both, i.e., GP-R1 is assigned to both $q_{g_m,L}$ and $q_{g_m,R}$.
The last option, {$\mbox{EDP}_2$}, considers other three neighboring cells around $\bx_1$ as well
(e.g., they are $\bI_2, \bI_3$, and $\bI_5$ in  \cref{subfig:GP-R3}) and
extend the minimum search to these cells to assign the resulting minimum 
order to all 16 Riemann states (i.e., four faces times four Riemann states per face)
that belong to $\bI_1$.
The study in~\cite{clain2011high} points out that {$\mbox{EDP}_0$} should not be used
to properly exit the MOOD loop within a finite number of cascades; hence it is not our choice.
On the other hand, {$\mbox{EDP}_2$} is more aggressive than {$\mbox{EDP}_1$} 
in reducing the reconstruction order at the cell under consideration. 
For this reason, our choice for GP-MOOD is {$\mbox{EDP}_1$} in this study.

%
%
%

\subsection{Quick summary of the GP-MOOD procedures}\label{sec:gp-mood_summary}
\begin{figure}[h!]
    \centering
    \includegraphics[scale = 0.31]{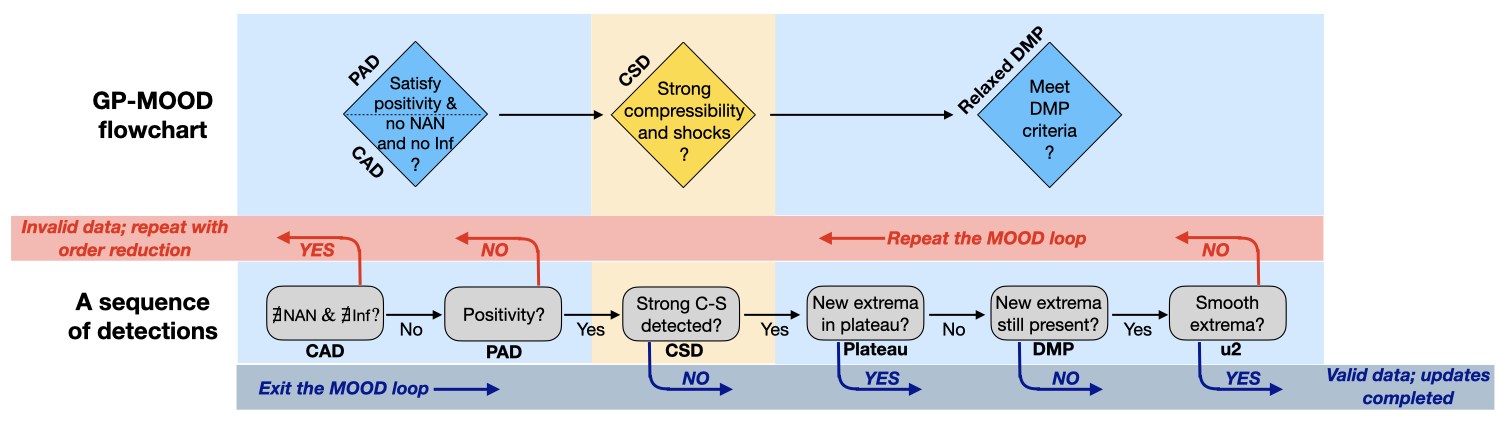}
    \caption{Flowchart of the detection at each sub-step used in the GP-MOOD loop.
    The sub-steps in the bottom portion with ``A sequence of detections'' belong to the 
    corresponding steps in the top portion with ``GP-MOOD flowchart.'' 
    Accordingly, the candidate solutions are rerouted from each of those sub-steps
    to either the next MOOD loop (invalid data on the red horizontal stripe in the middle)  or
    the exit of the MOOD loop (valid data on the grey-blue horizontal stripe in the bottom).
    All criteria except for PAD and CSD belong to NAD.}
    \label{fig:mood_loop}
\end{figure}

Putting all things together, 
we summarize the GP-MOOD procedures in~\cref{fig:mood_loop}, focusing on the
individual detection criterion at each stage.
\cref{fig:gp_mood_flowChart} (the top portion) 
is revisited in~\cref{fig:mood_loop} to map to
the corresponding sub-steps (the bottom portion).

\subsection{ Time integration}\label{sec:rk}
To achieve the anticipated target solution accuracy 
at 3rd-, 5th-, and 7th-order in GP-MOOD,
it is necessary that the consideration of temporal accuracy 
also needs to be accounted for correspondingly. 
Any simulations that are integrated with a temporal solver, 
whose accuracy is lower than
the spatial solver, will be penalized by the lower temporal accuracy.
The situation is easily understood when 
a $p$th order spatial method is integrated
with a $q$th order temporal solver. 
The error in discrete updates will be dominated 
by the leading error determined by either of the two cases, 
$\mathcal{O}(\Delta s^p) > \mathcal{O}(\Delta t^q)$ or 
$\mathcal{O}(\Delta s^p) < \mathcal{O}(\Delta t^q)$
with $\Delta s = \min\{\Delta x, \Delta y, \Delta z\}$ in 3D generally.
Two recent studies by Lee \textit{et al.}
~\cite{lee2021single,LEE2021100098}
have shown that,
in the case with the lower order temporal method, $q<p$, 
there exist a critical grid delta $\Delta s_\text{crit}$
beyond which (i.e., the grid delta $\Delta s$ becomes smaller than $\Delta s_\text{crit} $) 
the numerical error is severely penalized and follow
the convergence rate dictated by the lower temporal accuracy.
In numerical simulations, this temporal error dominance 
becomes more influential at high resolutions, creating
a computational dilemma of not making the expected
solution improvements by refining spatial grids.

As such, our strategy is to employ SSP-RK solvers
and match temporal accuracy with GP's spatial accuracy 
in all our convergence studies on smooth flows. 
In other shock dominant simulations, 
where a rigorous convergence rate can be lifted, 
we choose 
a relevant SSP-RK solver whose accuracy 
may be lower than the spatial GP solvers
to reduce the computational cost.

Two choices of our SSP-RK solvers include 
the optimal three-stage, 3rd-order SSP-RK3 solver
~\cite{shu1988efficient,gottlieb1998total} given as 
(we drop the cell index $ij$ for simplicity),
\begin{eqnarray}\label{eq:rk3}
\overline{\bU}^{n+1}=\frac{\overline{\bU}^{n}+2
\overline{\bU}^{(3)}}{3}, \;\; \text{ where } \;\;\; 
\begin{cases}
     \overline{\bU}^{(1)} = \overline{\bU}^{n}\  - \Delta t \, \mathbb{F}_{\nabla}(\overline{\bU}^{n}),\\
     \overline{\bU}^{(2)} =\overline{\bU}^{(1)} - \Delta t \,\mathbb{F}_{\nabla}(\overline{\bU}^{(1)}),\\
      \overline{\bU}^{(3)} = \hat{{\bU}}^{(2)} - \Delta t \,\mathbb{F}_{\nabla}(\hat{{\bU}}^{(2)} )\;
      \mbox{ with } \;\hat{{\bU}}^{(2)}= \frac{1}{4}(3\overline{\bU}^{n}+\overline{\bU}^{(2)}),\\
      \end{cases}
\end{eqnarray}
and
the optimal five-stage 4th-order SSP-RK4 solver~\cite{spiteri2002new},
\begin{eqnarray}\label{eq:rk4}
\begin{aligned}
\overline{\bU}^{n+1} & =  0.517231671970585\, \overline{\bU}^{(2)} \\
  & +  0.096059710526147\,\overline{\bU}^{(3)}  + 0.063692468666290\,\Delta t \,\mathbb{F}_{\nabla}(\overline{\bU}^{(3)} )\\
  &  + 0.386708617503269\,\overline{\bU}^{(4)} - 0.226007483236906\,\Delta t \,\mathbb{F}_{\nabla}(\overline{\bU}^{(4)}),
\end{aligned}
\end{eqnarray}
where the intermediate solutions at each sub-step are given as
\begin{eqnarray}\label{eq:rk4_substeps}
\begin{cases}
      \overline{\bU}^{(1)} = \overline{\bU}^{n} - 0.391752226571890\,\Delta t \,\mathbb{F}_{\nabla}(\overline{\bU}^{n}),\\
      \overline{\bU}^{(2)}  = 0.444370493651235\,\overline{\bU}^{n}
                                       +0.555629506348765\, \overline{\bU}^{(1)}
                           	      -0.368410593050371\,\Delta t \, \mathbb{F}_{\nabla}( \overline{\bU}^{(1)}),\\
     \overline{\bU}^{(3)} = 0.620101851488403\,\overline{\bU}^{n}
                    	            +0.379898148511597\, \overline{\bU}^{(2)}
                           	    -0.251891774271694\,\Delta t \, \mathbb{F}_{\nabla}( \overline{\bU}^{(2)} ),\\
     \overline{\bU}^{(4)} = 0.178079954393132\,\overline{\bU}^{n}
      		        		   +0.821920045606868\,\overline{\bU}^{(3)}
				    -0.544974750228521\,\Delta t \, \mathbb{F}_{\nabla}(\overline{\bU}^{(3)}).
\end{cases}
\end{eqnarray}
The SSP property ensures that 
if each sub-step solution is admissible, so is
the final updated solution $\overline{\bU}^{n+1}_{ij}$.

\subsection{Numerical stability and performance comparison}\label{sec:stability}
For the simulations reported in this paper, we find empirically that
the GP-MOOD methods are readily stable up to a CFL number of 0.8.
This stability limit is larger than 
the theoretical limits of the above SSP-RK methods
set by the so-called CFL coefficient $\mathcal{C}$ (e.g., \cite{spiteri2002new})
(or also known as the SSP coefficient, e.g., \cite{gottlieb2011strong}),
with which the stable time step is defined as $\dt = \mathcal{C} \dt_{\mbox{\tiny{FE}}}$.
The stability theory prescribes that 
$\mathcal{C} = 1/D$ for SSP-RK3 in Eq.~\eqref{eq:rk3} 
and  $\mathcal{C} = 1.508/D$ for SSP-RK4 in Eq.~\eqref{eq:rk4},
where $D$ is the spatial dimensionality of each problem.
As pointed out in~\cite{gottlieb2011strong}, 
$\mathcal{C}$ depends solely on time discretization,
and the first-order accurate 
Forward-Euler time step $\dt_{\mbox{\tiny{FE}}}$ depends solely on spatial
discretization. 
In our MOOD approach, the dependency of $\dt_{\mbox{\tiny{FE}}}$
is related to all participating high-order GP methods as well as
FOG. The \textit{a posteriori} nature of the MOOD order 
cascading algorithm makes it hard to perform a formal theoretical 
stability analysis of our methods.
As such, we instead conduct a numerical test to determine 
the numerical stability bound of our methods empirically.

To see this, we set up a 2D Sedov problem (see \cref{sec:sedov}) 
on a $400\times 400$ grid resolution and monitor the stability of a
baseline model chosen for the stability test purpose.
The choice of our baseline model consists of the 3rd-order setting,
namely GP-MOOD3 and SSP-RK3, solved with the HLLC Riemann solver.
The 4th-order accurate two-point quadrature rule is used by default for 
the 3rd-order setting. Besides, we tested our baseline model with the
2nd-order one-point quadrature rule to illustrate the stability sensitivity to
different quadrature rules. Lastly, we compare our baseline model
to the FOG-only model with one-point quadrature to help provide 
the stability gain afforded by the 3rd-order GP solution in GP-MOOD3.
Note that, for FOG, the use of the two-point (and other multi-point) 
quadrature rule is identical
to the use of the one-point quadrature rule since each volume-averaged
quantity is constant on each cell.

Let us first realize that there are a couple of typical signatures when 
a numerical method becomes unstable.
They include the presence of checkerboard patterns (or odd-even decoupling) 
that appear due to the lack of numerical stability used in a simulation;
unphysical oscillations near shocks and discontinuities, which grow rapidly in
time. The former issue can be cured by imposing added numerical dissipation,
which can be alleviated by switching to a more diffusive reconstruction method
or a Riemann solver, adding multidimensional wave information,
or reducing a CFL number. Similarly,
the latter situation can also be improved by adopting a more limited
reconstruction (irrelevant for GP-MOOD), 
a more diffusive Riemann solver, or reducing a CFL number.
In our GP-MOOD methods, we expect solutions to evolve with unlimited
high-order GP methods as long as numerical stability is sufficiently provided, 
except at shocks and discontinuities where FOG is 
to be chosen as the last resort after all unlimited high-order GP solutions
(should) fail there.

\begin{figure}[h!]
    \centering
  \includegraphics[scale=0.8]{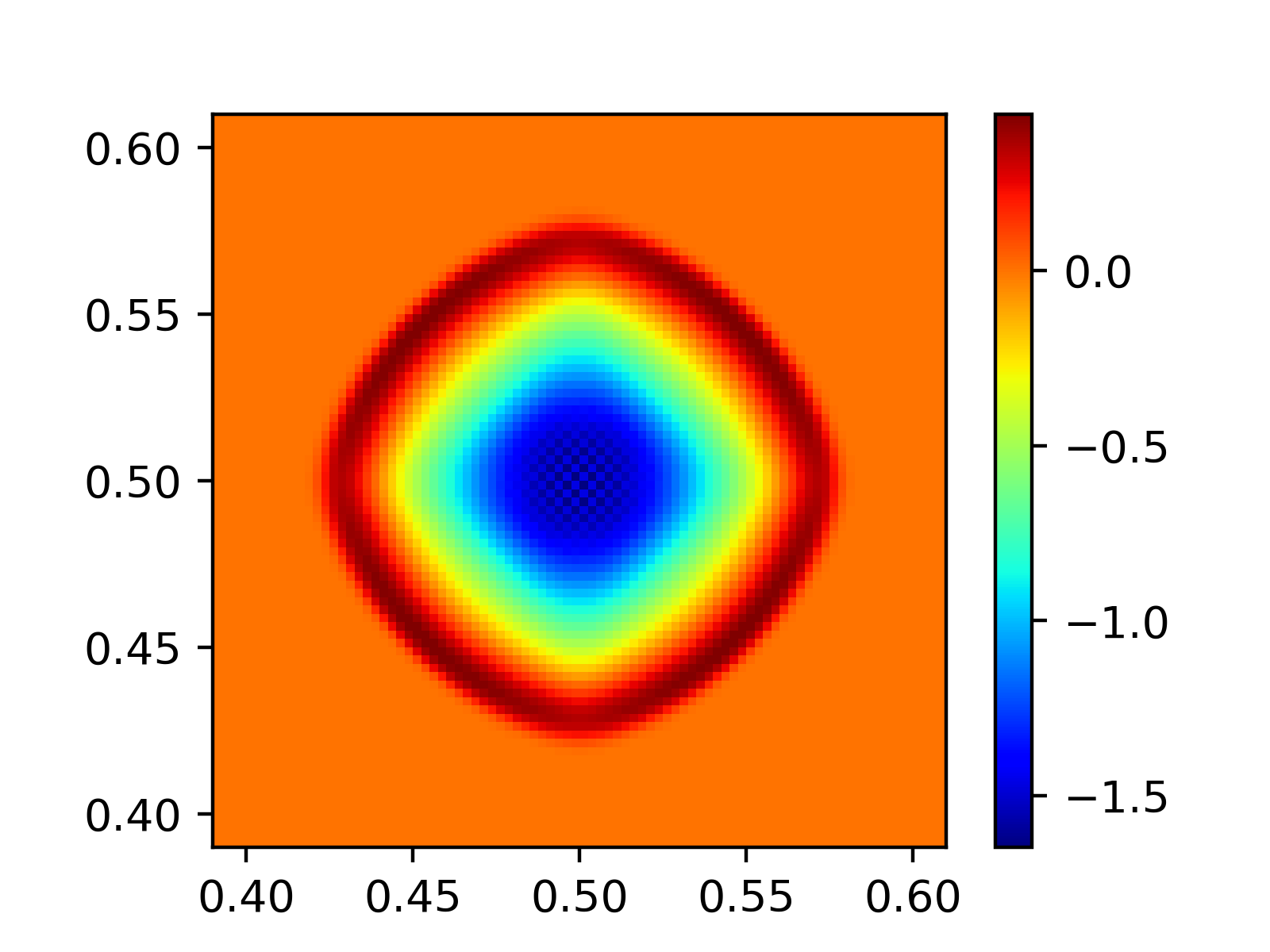}
  \includegraphics[scale=0.8]{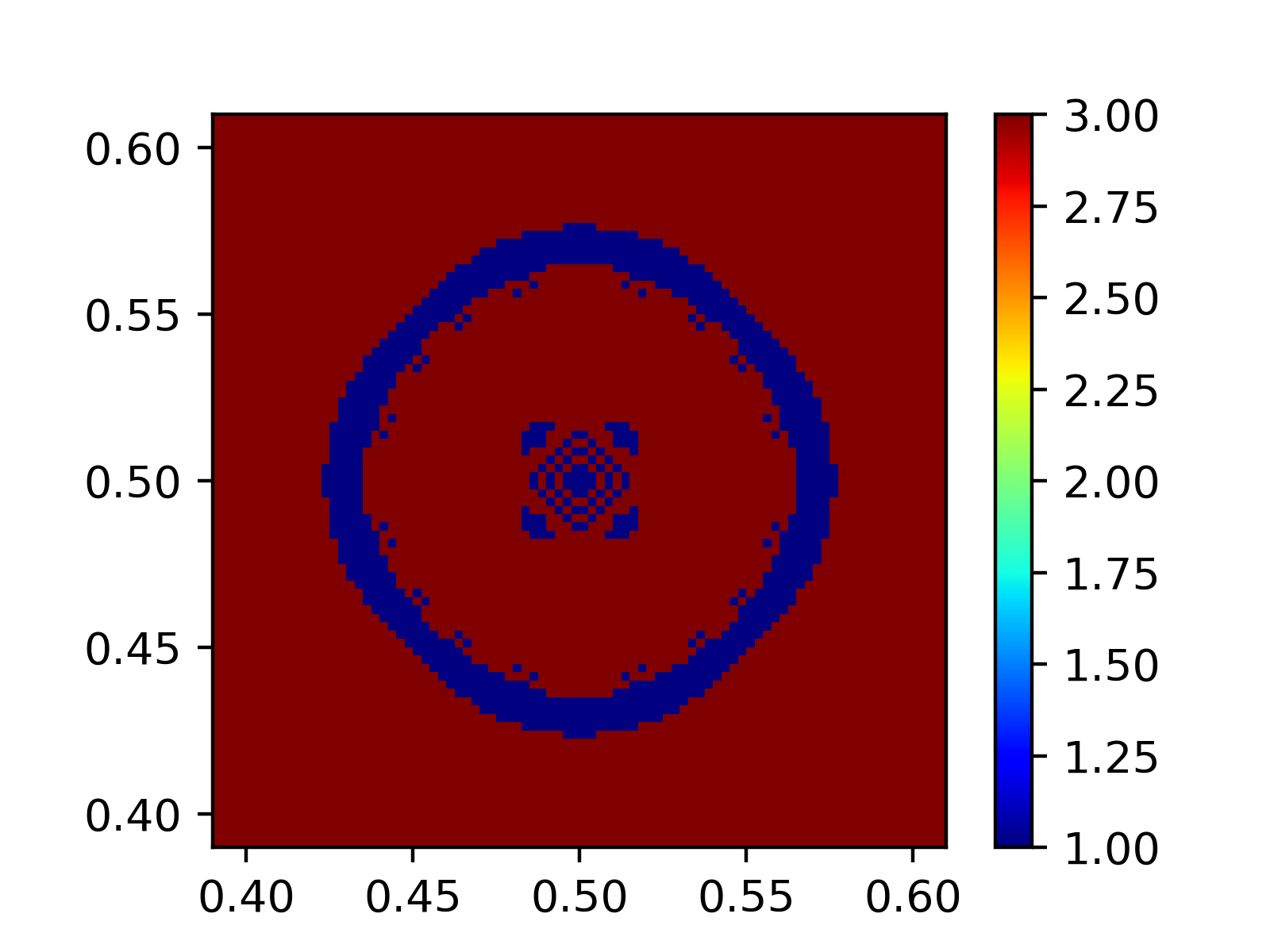}
    \caption{Examples of unstable runs of the Sedov explosion problem.
    (\textbf{left}) A checkerboard pattern appears in the central low-density region
    due to insufficient numerical stability in the simulation,
    configured with the FOG-only model using CFL=0.8. The plot displays
    the zoomed-in density at $t=5.5225 \times 10^{-3}$.
    The checkerboard pattern grows in time and eventually leads the run to
    crash.
    (\textbf{right}) Another unstable run using the baseline model with 
    1-point QR and CFL=0.85. The solution order distribution 
    (i.e., either ``1'' or ``3'' to display the local regions
    of the flow corresponding to the solutions with FOG and GP-R1, respectively)
    is plotted at $t=3.9044 \times 10^{-3}$.
    The FOG solutions (dark blue) are present 
    not only at shocks but also in the low-dense, smooth central region,
    the latter of which is related to a checkerboard-patterned 
    density formation (not shown) similar to the left panel.
    }
    \label{fig:sedov_unstable}
\end{figure}

Based on these considerations, we grant a method
is \textit{numerically stable} if the method satisfies the following criteria:
(i) there is no checkerboard pattern, particularly in the low dense central region
(see \cref{fig:sedov_unstable}),
(ii) there is no unphysical growth of variables in magnitude, 
(iii) there is no crash caused by \texttt{NAN} or \texttt{Inf},
(iv) there is no FOG solution present in the smooth central region,
i.e., the unlimited 3rd-order GP solver successfully produced 
stable solutions in the central region 
(see the right panel in \cref{fig:sedov_unstable}),
and finally 
(v) the Sedov explosion retains self-similar and symmetry.

\begin{figure}[h!]
    \centering
  \includegraphics[scale=0.38]{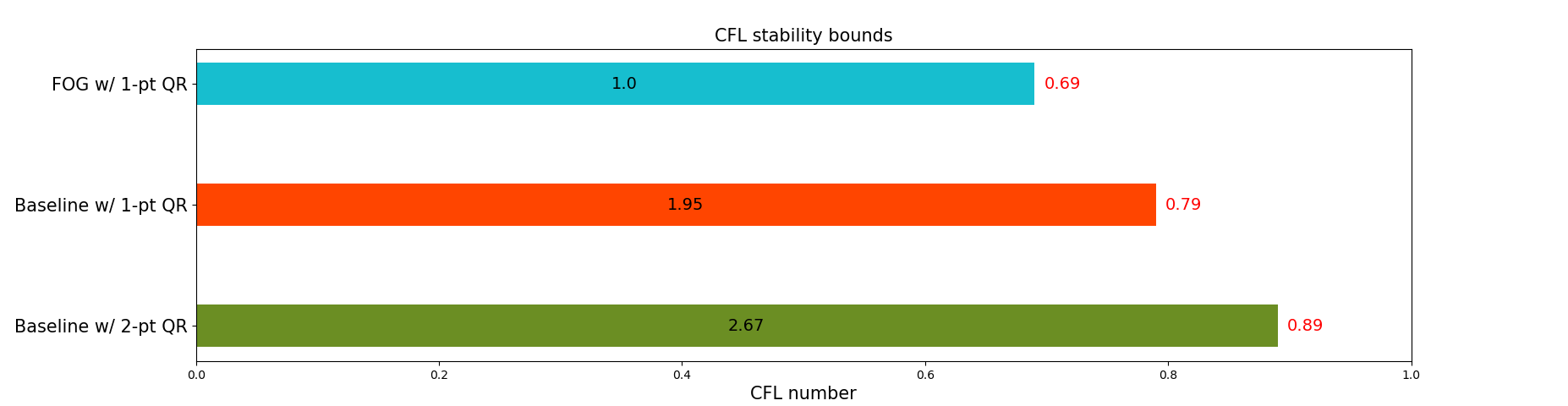}
    \caption{Empirical CFL stability bounds for GP-MOOD3 and FOG. The baseline
    model is the combination of GP-MOOD3 with HLLC and SSP-RK3.
    The FOG-only model (top) is solved using 1-point quadrature rule (QR), while
    the baseline models are solved using 1-point (middle) and 2-point (bottom)
    QRs to test the 2D Sedov problem 
    resolved on a $400\times 400$ grid resolution.
    The numbers in black labeled on each bar graph represent the CPU performance
    for each corresponding CFL bound, normalized by the CPU performance of the FOG-only
    model (top). The numbers in red next to each bar 
    are the empirical CFL bounds that ensure numerical
    stability in each test case. The performance results are the
    average of five independent runs
    using non-parallel, serial calculations on 
    a 2.4 GHz 8-Core Intel Core i9 CPU.
    }
    \label{fig:CFL_stability}
\end{figure}

The results are displayed in \cref{fig:CFL_stability}. The FOG-only model is shown
to be stable up to CFL=0.69, beyond which the central region exhibits severe
checkerboard patterns, leading to crash. This odd-even decoupling in the
central region is significantly improved using the \textit{multidimensional}
3rd-order GP-R1, which gains more stability by adding
wave information from each transverse direction. As a result,
the CFL stability bound grows from CFL=0.69 with FOG-only to CFL=0.79 
with the baseline model with the 1-point quadrature rule (QR).
Doubling the number of quadrature points to the 2-point QR further
enhances the CFL bounds to CFL=0.89 with the baseline model.
We also measured the CPU performance of the time-to-solution in each model,
respectively, with CFL=0.69, 0.79, and 0.89.
Normalized by the FOG-only model, the baseline run with the
1-point QR is 1.95 times more expensive than FOG-only, while 
the ratio becomes 2.67 for the baseline run with the 2-point QR.
In addition, the baseline run with the 2-point QR and CFL=0.89 
(the bottom dark green bar in  \cref{fig:CFL_stability})
 is 1.37 times more 
expensive than the baseline run with the 1-point QR and CFL=0.79
(the middle red bar in \cref{fig:CFL_stability}).
On the other hand, the CPU performance ratio at CFL=0.69
of the baseline model with the 1-point QR (not shown) 
to the FOG-only model (the top cyan bar in \cref{fig:CFL_stability}) is
observed to be 2.07, while the ratio at CFL=0.79 of the baseline
model with 2-point QR (not shown) to 
the baseline model with the 1-point QR 
(the middle red bar in \cref{fig:CFL_stability})
is 1.60.

Note that the baseline model with the 2-point QR 
(the bottom dark green bar in  \cref{fig:CFL_stability})
is the least
stable option among the test cases in \cref{sec:results}.
From the results herein, we conclude that the GP-MOOD methods studied
in this paper are sufficiently stable with CFL = 0.8; hence is
the default choice in \cref{sec:results}.
In all stable GP-MOOD3 runs tested in this section,
the percentage of those cells solved by FOG
as the result of MOOD's order reduction
never exceeds 2\% of the entire cells,
leaving about 98\% of the domain solved by the unlimited 
3rd-order GP method.
This fact can also be viewed as an algorithmic efficacy of GP-MOOD's
\textit{a posteriori}  approach over the classical \textit{a priori} paradigm,
where, in the latter, the use of computationally intensive nonlinear limiters 
is required to be always fully activated and calculated on all cells for stability.
This experiment
indicates that the nonlinear limiters are not required on 98\% of them.
Such an unnecessary computation is efficiently circumvented in GP-MOOD.

\section{Stepwise implementation of the GP-MOOD method}\label{sec:stepwise}
In this section, we summarize the GP-MOOD method proposed in this study in a stepwise fashion.
We intend to provide a bird-eye view of GP-MOOD for ease of practical implementation.

\begin{enumerate}
    \item[] \textbf{Step 1:}
    Choose a hyperparameter, $\ell$, for each simulation.
    We recommend a constant value (see~\cref{sec:isentropic_vortex}) for a convergence study while
    our default choice $\ell = 6\Delta$ or $\ell = 12\Delta$ with $\Delta = \min\{\dx, \dy\}$ works well for most 
    of the problems in ~\cref{sec:results}.
    \item[] \textbf{Step 2:}
    Once a grid is configured, calculate the GP covariance kernels in~\cref{eq:SE-pred,eq:SE-cov} 
    and compute the GP prediction vector $\bz_*^T$ according to
    ~\cref{eq:gp_reconstructor,eq:normalized_z} for $R=1$ with the chosen $\ell$.
    Save the computed $\bz_*^T$ for later use.
    Repeat the same for $R=2$ if GP-MOOD5 is considered.
    Instead, repeat the same for $R=3$ if GP-MOOD7 is considered.
    Group these GP methods according to the GP-MOOD3, GP-MOOD5, and GP-MOOD7
    cascading families alongside FOG.
    \item[] \textbf{Step 3:}
    Start a simulation. The 3rd-order GP-R1 is solved with the 2-point quadrature rule;
    the 5th-order GP-R2 with the 3-point rule; the 7th-order GP-R3 with the 4-point rule.
    Choose a proper SSP-RK method. After each $s$th RK sub-stage at $n$th timestep,
    each of the procedures in the MOOD loop in~\cref{fig:gp_mood_flowChart} will be
    tested on the $s$th candidate solution, $\overline{\bU}^{*,rk(s)}_{ij}$.
    \item[] \textbf{Step 4:}    
    Execute the CAD criteria in~\cref{eq:CAD}.
    If there are no \texttt{NAN}'s and no \texttt{Inf}'s, the solution
    moves to \textbf{Step 5}; otherwise, record the solution as inadmissible  and
    conduct the MOOD order reduction. 
    Repeat the MOOD loop with the next accurate
    reconstruction method starting from \textbf{Step 4}.
    \item[] \textbf{Step 5:}    
    Execute the PAD criteria in~\cref{eq:PAD}.
    If the candidate solution passes PAD, go to \textbf{Step 6}; otherwise,
    record the solution as inadmissible and 
    conduct the MOOD order reduction. 
    Repeat the MOOD loop with the next accurate
    reconstruction method starting from \textbf{Step 4}.
    Note that the order of operations 
    between \textbf{Step 4} and \textbf{Step 5} can be
    swapped.
    \item[] \textbf{Step 6:}    
    Run the CSD criteria in~\cref{eq:divV,eq:gradP}.
    If there are no strong compression 
    and no strong pressure
    gradients in the local gas, record the candidate solution as admissible and accept it
    as the final solution. Finish the MOOD loop.
    Otherwise, go to \textbf{Step 7}.
    \item[] \textbf{Step 7:}    
    Check the Plateau criterion in~\cref{eq:plateau}.
    If the solution satisfies \cref{eq:plateau}, record the candidate solution
    as admissible and take it as the final solution. Finish the MOOD loop.
    Otherwise, go to \textbf{Step 8}.
    \item[] \textbf{Step 8:}        
    Check DMP in~\cref{eq:DMP_RK}. If the solution satisfies
    \cref{eq:DMP_RK}, record it as admissible and take it as the final solution.
    Exit the MOOD loop. Otherwise, go to \textbf{Step 9}.
    \item[] \textbf{Step 9:}     
    Run the \textit{u2} checks in~\cref{eq:u2_curvature,eq:u2_a,eq:u2_b} as the last MOOD test.
    If the candidate solution meets either~\cref{eq:u2_a} or ~\cref{eq:u2_b},
    take it as the final admissible solution and exit the MOOD loop.
    Otherwise, conduct the MOOD order reduction by one cascade.
    Repeat the MOOD loop with the next accurate
    reconstruction method starting from \textbf{Step 4}.
   \item[] \textbf{Step 10:}
   Once the MOOD loop reaches the safe method, FOG, 
   take it as the final admissible solution. Exit the MOOD loop.
   No further MOOD check is needed at this point.

\end{enumerate}
The entire GP-MOOD pipeline is pictured as a flowchart in \cref{fig:gp_mood_flowChart}, while 
a zoomed-in view that focuses more on the MOOD criteria is schematically
provided in \cref{fig:mood_loop}.

\section{Numerical results}\label{sec:results}
This section displays three types of test problems, including
\begin{itemize}
\item[(i)] a grid convergence test problem in 2D that assesses the accuracy of 
the proposed GP-MOOD methods (see \cref{sec:isentropic_vortex}),
\item[(ii)]  standard shock-dominant benchmark problems in 
1D and 2D that test the validity of our GP-MOOD solvers
on a set of well-known CFD problems (see \cref{sec:shock_problems}), and finally,
\item[(iii)]  highly compressible, strong shock problems that are known as 
stringent for numerical testing purposes.
These problems require a strong positivity-handling and expensive
nonlinear shock-limiters in the \textit{a priori} FV literature
(see \cref{sec:Mach_jets}).
\end{itemize}

The choice of a Riemann solver is default to be the 
HLLC Riemann solver~\cite{toro1994restoration} in all test problems,
except that the HLL Riemann solver~\cite{harten1983upstream} is used
for the last two test problems in \cref{sec:Mach_jets}
to suppress the grid-aligned carbuncle instabilities~\cite{quirk1997contribution}.
For GP, we follow ~\cite{reyes2018new,reyes2019variable}
to set 
$\ell = 1$ for the 2D convergence test problem in \cref{sec:isentropic_vortex} and
$\ell = 6\Delta$ for the 1D Shu-Osher shock tube problem in \cref{sec:shu_osher},
while we set $\ell = 12\Delta$ in all other problems in this section.
As discussed in~\cref{sec:stability}, we take a CFL number of 0.8 by default.
The ratio of specific heats is set to be $\gamma = 1.4$ in all test problems.


\subsection{Grid convergence -- The isentropic vortex test}\label{sec:isentropic_vortex}
The accuracy of the GP-MOOD schemes is considered on a well-known 2D test problem
called the nonlinear isentropic vortex advection presented by Shu \cite{shu1998essentially}.
As in~\cite{reyes2019variable} we double the original domain size to avoid self-interactions
of the vortex across the periodic domain.
It sets up a circular region centered at $(x_c, y_c) = (10, 10)$ on a periodic square domain,  
$[0,20]\times [0,20]$, 
where a Gaussian-shaped vortex with rotating velocity fields are initialized.
As the problem consists of the smooth advection of the vortex along the diagonal direction, 
any departure from the initial condition (or the exact solution of the problem)
will be considered numerical errors of the numerical method under consideration.

We follow the standard initial condition in \cite{shu1998essentially} to initialize 
\textit{pointwise} values of the primitive variables as
%
%
\begin{eqnarray}
\rho(x,y) &=& \left[ 1 - \left(\gamma -1\right)\frac{\beta^2}{8\gamma\pi^2}e^{1-r^2}\right]^{\frac{1}{\gamma-1}},
\label{eq:vortex_IC_dens}\\
u(x,y) &=& 1 - \frac{\beta}{2\pi}e^{\frac{1}{2}\left(1-r^2\right)}(y-y_c),\label{eq:vortex_IC_velx}\\
v(x,y) &=& 1 + \frac{\beta}{2\pi}e^{\frac{1}{2}\left(1-r^2\right)}(x-x_c),\label{eq:vortex_IC_vely}\\
p(x,y) &=& \rho(x,y)^\gamma,\label{eq:vortex_IC_pres}
\end{eqnarray}
with $r = r(x,y) = \sqrt{(x-x_c)^2 + (y-y_c)^2}$ and the vortex strength $\beta = 5.$
With the mean diagonal flow velocity fields, $(u,v)=(1,1)$, the vortex makes one period
of diagonal advection and returns to the initial position at $t_{\tiny{\max}} = 20,$
at which point we measure numerical errors.

For a successful error analysis, 
it is crucial to convert these pointwise values to the 
\textit{volume-averaged conservative quantities} -- 
the data type evolved in FV --
to ensure that there is no accuracy loss in FV evolutions.
If one initializes conservative quantities with the pointwise values
without converting to volume-averaged quantities,
the simulation will introduce an
error of order two (i.e., $\mathcal{O}(\dx^2 + \dy^2)$) 
right at the initial step even before any solution evolution.
This initial error will deter us from assessing the correct convergence rates
of our GP-MOOD methods.
This concern was clearly pointed out in our earlier discussion, 
(i), (ii), and (iii) in \cref{sec:governing_eqns}.
To this end, we perform the followings:

\begin{enumerate}
    \item[a)] We convert the pointwise primitive variables 
    (i.e., the initial conditions in Eq.~\eqref{eq:vortex_IC_dens} -- Eq.~\eqref{eq:vortex_IC_pres})
     to the corresponding pointwise conservative variables and use 
     a sufficiently highly accurate quadrature rule to convert them to 
     the corresponding volume-averaged conservative quantities.
     We use the 10th-order accurate 5-point Gauss-Legendre quadrature for the conversion 
     to guarantee that the initial setup accuracy is much higher than all of the rest discrete operations in
     the simulation.
     \item[b)] To match the temporal accuracy of SSP-RK3 and SSP-RK4 with those of GP's spatial accuracy, 
     we further reduce time steps $\dt$ so that the temporal error is scaled down to (or smaller than or equal to) 
     the spatial error. For example, the time step size of SSP-RK4 is scaled as constant according to
     $\Delta t = \Delta^{\frac{5}{4}}$, $\Delta = \min\{\dx,\dy\}$,
     when used with the 5th-order GP-MOOD5.
\end{enumerate}

\begin{figure}[h!]
    \centering
    \begin{subfigure}{8cm}\centering
    \includegraphics[width=\textwidth]{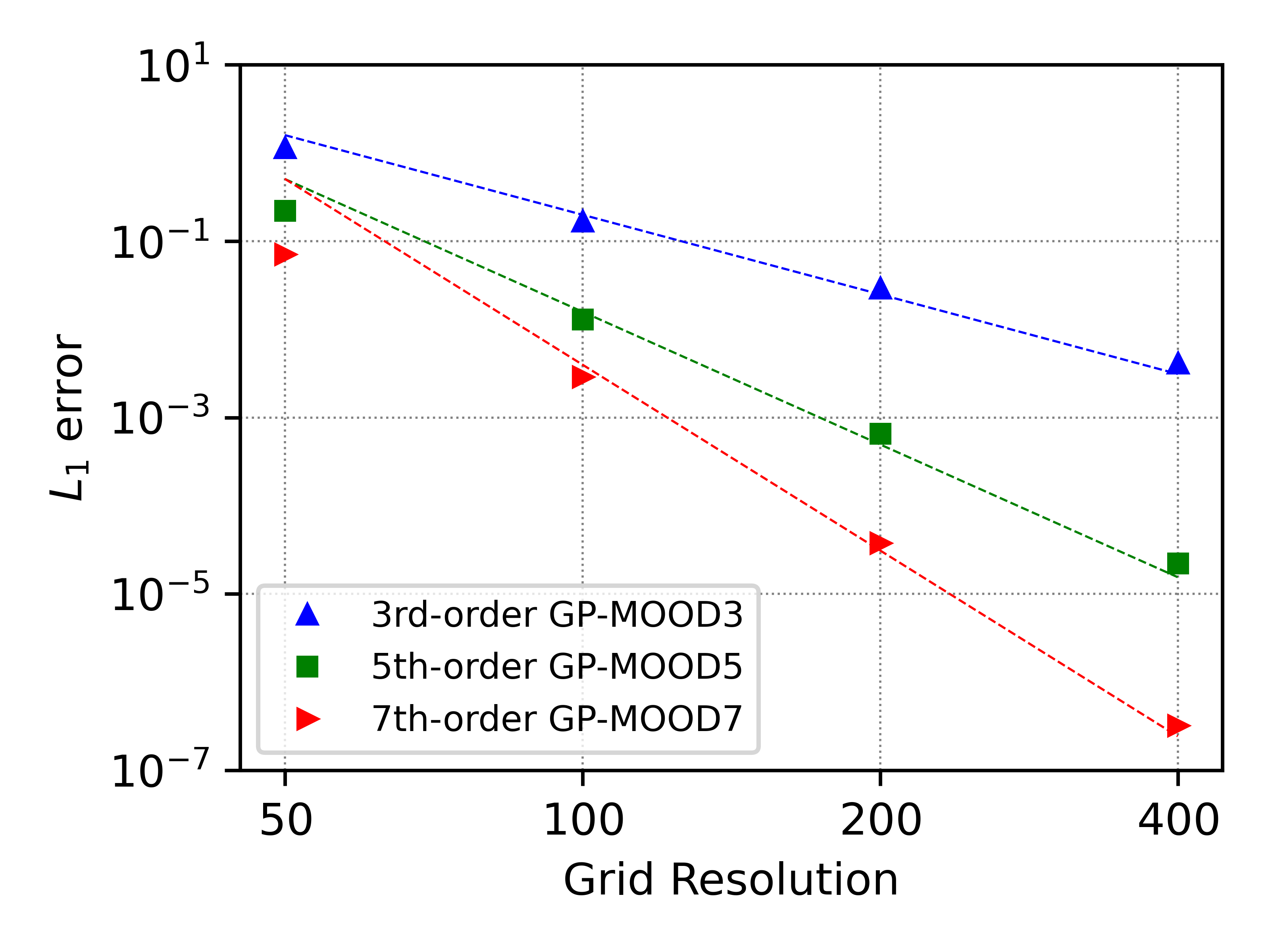}
    \caption{}\label{fig:gp_conv_sphere}
    \end{subfigure}
    \begin{subfigure}{8cm}\centering
    \includegraphics[width=\textwidth]{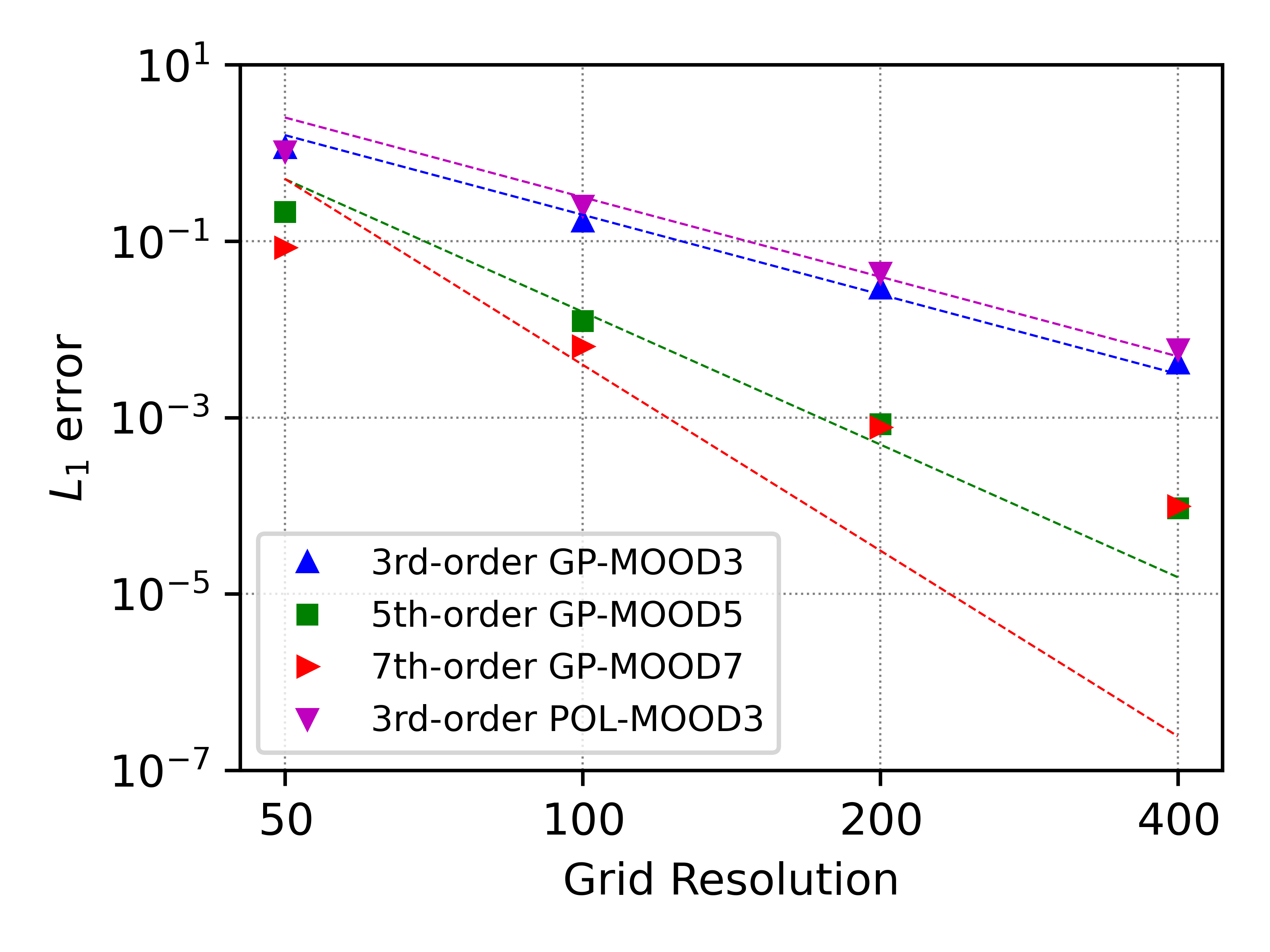}
    \caption{}\label{fig:gp_conv_cross}
    \end{subfigure}
    \caption{Convergence study on the 2D isentropic vortex advection
    for three different radii, stencils, and reconstructions. 
    (\textbf{left}) The convergence rates on the blocky-diamond
    GP stencils.
    (\textbf{right}) The convergence rates on the cross-shaped stencils.
    CFL=0.8 is used for all cases with proper time step size reductions.}
    \label{fig:gp_convergence}
\end{figure}

For GP, we set $\ell = 1$ in all test cases considered here, based on our former studies
on the \textit{a priori} GP methods \cite{reyes2018new,reyes2019variable}.
Tested in \cref{fig:gp_convergence} include 
the 3rd-order GP-MOOD3 (GP-R3 $\to$ FOG),
the 5th-order GP-MOOD5 (GP-R5 $\to$ GP-R3 $\to$ FOG), and
the 7th-order GP-MOOD7 (GP-R7 $\to$ GP-R3 $\to$ FOG).
In~\cref{fig:gp_conv_sphere}, the results are solved on the blocky-diamond GP stencils 
as depicted in~\cref{fig:gp_stencil_GP-R1,fig:gp_stencil_GP-R2-R3}.
In~\cref{fig:gp_conv_cross} though, the same tests are repeated on cross-shaped GP stencils
that are the direct extension of the GP-R1 stencil in \cref{fig:gp_stencil_GP-R1}
to a bigger cross stencil by adding extra cells in each normal direction only 
according to the GP radius $R$.
While the cross GP stencil and the blocky sphere stencil are the same for GP-R1,
the size of the cross increases to the 9-point stencil for GP-R2 and 13-point for GP-R3.
%
%
%

The results of $L_1$ errors are reported in ~\cref{fig:gp_convergence} 
on four different grid resolutions, $N_x=N_y = 50, 100, 200$, and $400$.
In \cref{fig:gp_conv_sphere}, the convergence rates of the three GP-MOOD methods 
on the corresponding diamond GP stencil follow the analytical convergence rates (dotted lines) 
of $(2R+1)$, showing the expected 3rd-, 5th-, and 7th-order rates, respectively.
See also \cref{table:convergence-gp-mood-diamond}. The experimental order
of convergence (EOC) is computed as
\beq
\mbox{EOC}  = \frac{\ln(E_c/E_r)}{\ln(2)},
\eeq
where $E_c$ and $E_r$ are the $L_1$ errors 
on the coarse and the next coarse resolutions 
(e.g., $E_c$ on $50 \times 50$ and $E_r$ on $100 \times 100$),
respectively.
The demonstration of the full convergence rate in each GP-MOOD method 
proves that the GP-MOOD detection algorithms operate successfully without
any erroneous order reduction on this smooth advection problem,
and the solution evolves with the highest accurate, unlimited GP reconstruction method
in each case. 

However, in \cref{fig:gp_conv_cross} and \cref{table:convergence-gp-mood-cross},
the accuracy of GP-MOOD is heavily compromised when the solutions evolve on the smaller
cross-shaped stencils (i.e., 9 vs. 13 stencil points for GP-MOOD5; 
13 vs. 25 for GP-MOOD7.
See \cref{fig:gp_stencil_GP-R2-R3}).
We see that the rate of convergence is disturbed in GP-MOOD5
and GP-MOOD7. In both cases, the orders of accuracy asymptotically converge at
3rd-order, with both errors reduced by about two orders of magnitude compared to
the GP-MOOD3 error.
We also show the convergence rate of the polynomial counterpart, POL-MOOD3, in
\cref{fig:gp_conv_cross} and \cref{table:convergence-gp-mood-cross}.
Similar to that of GP-MOOD3, the 3rd-order rate of convergence is achieved except
that POL-MOOD3's magnitude of the error on each grid resolution is about 1.5 times larger than
the reported error with GP-MOOD3.
%
%

\begin{table}[ht!]
    \footnotesize
    \centering
    \caption{Performance results of three different GP-MOOD algorithms on the blocky-diamond GP stencils.
    }\label{table:convergence-gp-mood-diamond}
    \begin{tabular}{@{}cccccccccccc@{}}
        \toprule
        \multirow{2}{*}{Grid Resolution} & \multicolumn{2}{l}{GP-MOOD3} &  & \multicolumn{2}{l}{GP-MOOD5}  & & \multicolumn{2}{l}{GP-MOOD7} \\ 
                                                              \cmidrule(lr){2-3} \cmidrule(l){5-6} \cmidrule(l){8-9} 
                                                          & $L_1$ errors & EOC &  & $L_1$ errors & EOC & & $L_1$ errors & EOC \\
        \midrule
        \( 50 \times 50 \)     & 1.13746068e+00   & --      &  & 2.22016430e-01    & --      & &  7.11860045e-02   & --    \\ 
        \( 100 \times 100 \) & 1.67459246e-01   & 2.76 &  & 1.29737576e-02   & 4.10  & &  2.90834888e-03  & 4.61 \\ 
        \( 200 \times 200 \) & 2.91650062e-02   & 2.52 &  & 6.60244975e-04   & 4.30  & & 3.76226930e-05   & 6.27  \\ 
        \( 400 \times 400 \) & 4.11626679e-03   & 2.82 &  & 2.23315572e-05    & 4.89 & &  3.21581083e-07  & 6.87  \\ 
    \end{tabular}
\end{table}

\begin{table}[ht!]
    \footnotesize
    \centering
    \caption{Performance results of GP-MOOD5, GP-MOOD7, and POL-MOOD3 on the cross GP stencils.
    The performance of GP-MOOD3 is omitted here since there is no discrepancy 
    between the cross and diamond stencil configurations for GP-MOOD3;
    see \cref{table:convergence-gp-mood-diamond} for its performance.
    }\label{table:convergence-gp-mood-cross}
    \begin{tabular}{@{}cccccccccccc@{}}
        \toprule
        \multirow{2}{*}{Grid Resolution} & \multicolumn{2}{l}{GP-MOOD5}  & & \multicolumn{2}{l}{GP-MOOD7}  & & \multicolumn{2}{l}{POL-MOOD3} \\
                                                              \cmidrule(lr){2-3} \cmidrule(l){5-6} \cmidrule(l){8-9} 
                                                          & $L_1$ errors & EOC &  & $L_1$ errors & EOC & & $L_1$ errors & EOC \\
        \midrule
        \( 50 \times 50 \)     & 2.15245127e-01    & --      & &  8.50162184e-02   & --   &  &  1.04593597e+00 & --  \\
        \( 100 \times 100 \) & 1.24973881e-02   & 4.10  & &  6.39925304e-03  & 3.73  &  &   2.54731899e-01 & 2.04 \\
        \( 200 \times 200 \) & 8.42402789e-04   & 3.89  & & 7.75915293e-04   & 3.04 &  &  4.40709959e-02  & 2.53 \\
        \( 400 \times 400 \) & 9.37519200e-05   & 3.17 & &  9.92175024e-05  & 2.97  &  &  5.96194298e-03 & 2.89 \\
    \end{tabular}
\end{table}

\subsection{Shock-dominant benchmark test problems}\label{sec:shock_problems}
Well-known standard shock-dominant benchmark problems 
in 1D and 2D are tested in this section
to discuss the validity of our GP-MOOD solvers.

\subsubsection{The Shu-Osher shock tube test in 1D}\label{sec:shu_osher}

\begin{figure}[htpb!]
    \centering
    \includegraphics[scale=0.85]{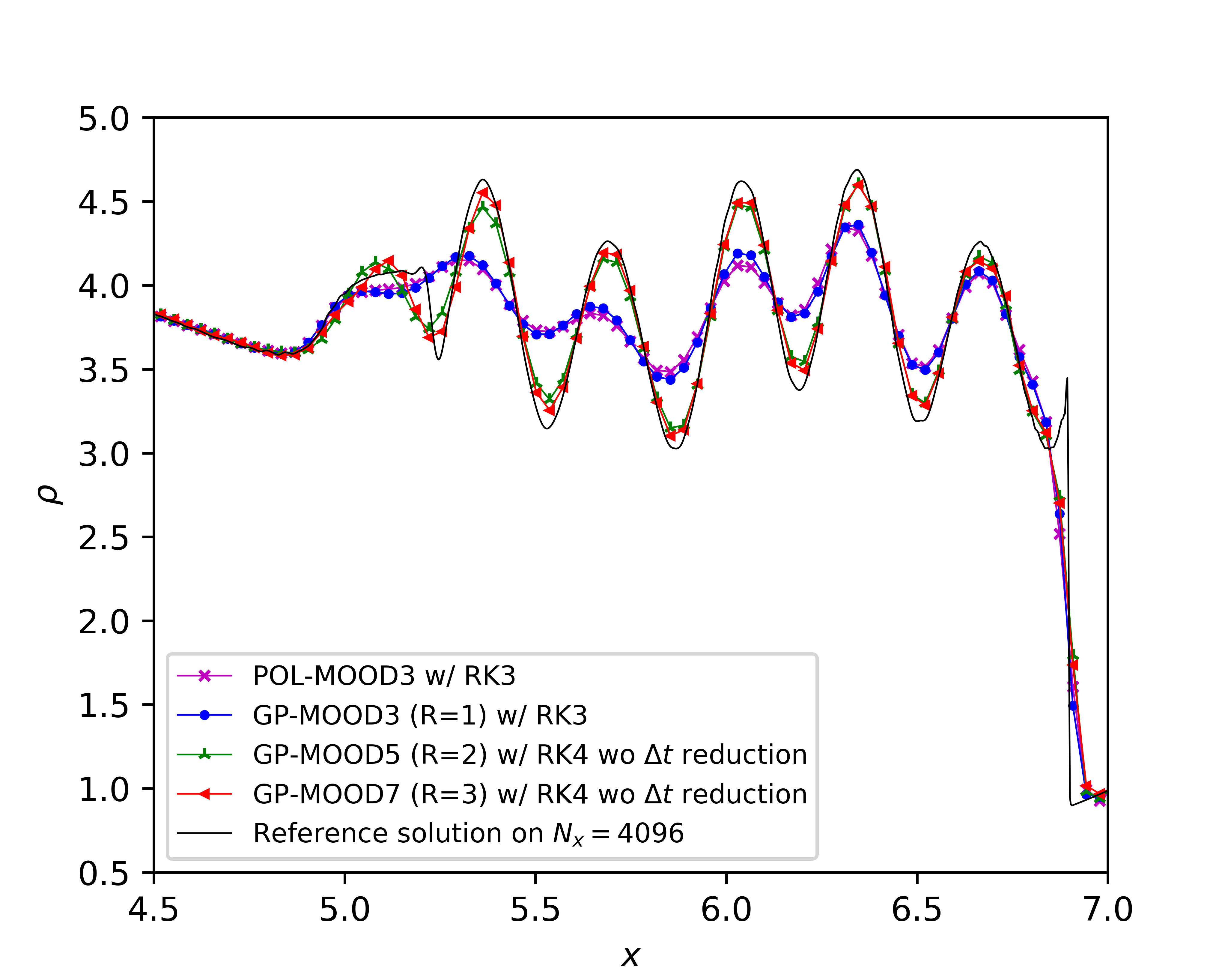}
    \caption{Result of the Shu-Osher shock tube test at $t=1.8$
    using GP-MOOD3, GP-MOOD5, GP-MOOD7 with $\ell = 6 \dx$, 
    CFL =0.8, resolved on 256 grid cells. The GP solutions 
    are plotted along with the POL-MOOD3 solution on the same
    grid resolution. The reference solution is computed using
    POL-MOOD3 with SSP-RK3 without reducing $\dt$ on
    $N_x = 4,096$.
    }
    \label{fig:shuOsher_test}
\end{figure}

The Shu-Osher shock tube problem \cite{shu1989efficient} provides 
a good comparison study of different numerical algorithms to measure relative
accuracy and diffusivity. The point of interest in this problem is to see
how well different methods can produce the high and low-frequency
regions of the density profile, as the Mach 3 shock wave travels 
from left to right into the initial low-frequency sinusoidal profile.
As a result, the perturbed density is compressed at the 
shock-smooth profile interactions, creating a high-frequency 
trail of doubled-frequency
waves to the left of the immediate post-shock region.
Further down to the left in the post-shock region,
attached to the tail of the high-frequency region
is the low-frequency structure -- now substantially sharpened
from the initial sinusoidal profile due to the shock-steeping --
where the perturbation returns to the 
original frequency of the unshocked initial wave.
We follow the standard setup of the initial condition
\cite{shu1989efficient},
\begin{equation}
    (\rho, u, p) =  \begin{cases}
      (3.857143,\ 2.629369,\ 10.33333) &\text{if} \quad x < 0.5,\\
      (1 + 0.2\sin(5(x-4.5)),\ 0,\ 1)  &\text{if} \quad  x >  0.5,
     \end{cases}
\end{equation}
on a one-dimensional domain $[0,9]$ and we evolve the solution until $t=1.8$.
The boundary condition at $x=0$ and $x=9$ is the fixed boundary 
condition that sets the initial values during the simulation.
All simulations are computed on a grid resolution, $N_x = 256$,
except for the reference solution resolved on the 16 times more refined
resolution, $N_x = 4,096$, using POL-MOOD3 and SSP-RK3 without
any time step reduction.
The time steps of SSP-RK4 are reduced for GP-MOOD5 and GP-MOOD7
to match the accuracy of the spatial and temporal solvers.

Simulated results are compared in \cref{fig:shuOsher_test}, 
where we display a zoomed-in view of the
entire profile to focus on the  solution comparison
over the high-frequency region.
Overall, all results produce acceptable density profiles 
capturing the assumed high-frequency amplitudes fairly well,
conforming with the reference solution.
The amplitude closest to the reference profile is achieved by
GP-MOOD7, followed by GP-MOOD5, GP-MOOD3, and POL-MOOD3.
It is quite impressive to see how closely the solutions 
produced by GP-MOOD5 and GP-MOOD7 
follow the reference profile, with only on 
the grid resolution 16 times lower than the reference solution.
To be more quantitative,
we compare the maximum values of the density,
measured at $x\approx 6.34$.
The reference solution reaches $\rho \approx 4.69$,
while both GP-MOOD5 and GP-MOOD7 
give $\rho \approx 4.60$, providing 98\% accuracy in
predicting the largest amplitude with only
utilizing 6.25\% of the reference grid resolution.
For GP-MOOD3, $\rho \approx 4.36$, yielding
93\% accuracy, while $\rho \approx 4.32$ with 
92\% accuracy for POL-MOOD3.

Lastly, we remark that it is crucial to utilize the
CSD criterion described in \cref{sec:CDS} in order to
produce the well-systematic solutions
as in \cref{fig:shuOsher_test}, namely, GP-MOOD7 being the most accurate, 
GP-MOOD5 intermediate, and GP-MOOD3 the least accurate.
If CSD is not employed, numerical diffusivity of the GP solutions 
becomes excessive, and their solution
qualities turn into a disorderly fashion, e.g., 
the GP-MOOD5 solution turns out to be
as diffusive as the GP-MOOD3 solution.

\subsubsection{The Sedov explosion test in 2D}\label{sec:sedov}
Next, we revisit the Sedov blast test \cite{sedov1993similarity}. 
The primary purpose is to test how well the GP methods retain
the symmetry of the self-similar evolution of the spherical shock
propagation outward from a highly pressurized point-source at 
the domain center, $(x,y) = (0,0)$.
We discretize a 2D outflow domain, $[-0.5, 0.5] \times [-0.5, 0.5]$, at a
grid resolution, $256 \times 256$, where the initial condition
is configured by following the description in \cite{fryxell2000flash}.

\begin{figure}[h!]
    \centering
    \includegraphics[scale=0.8]
    {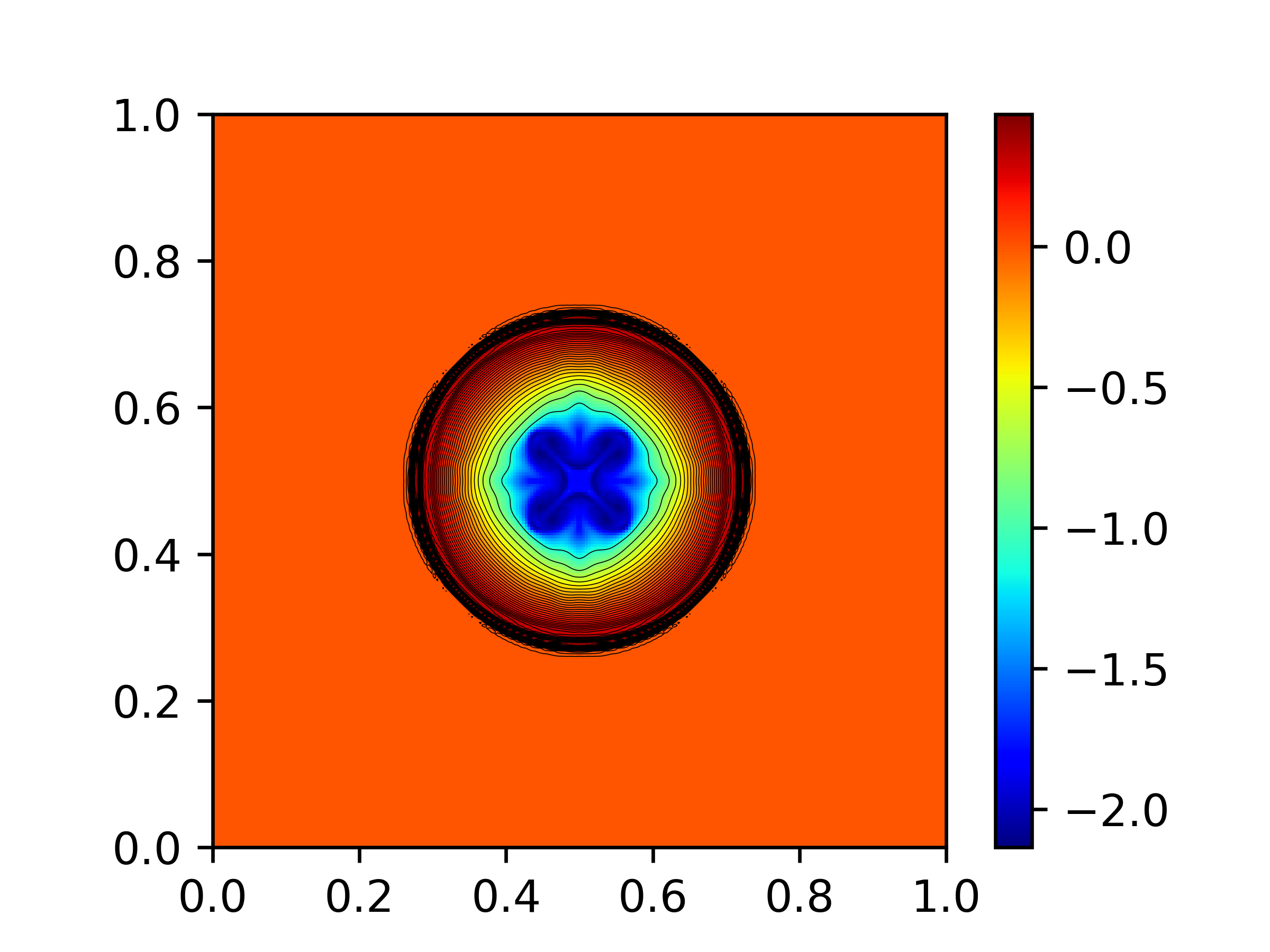}
    \includegraphics[scale=0.75]
    {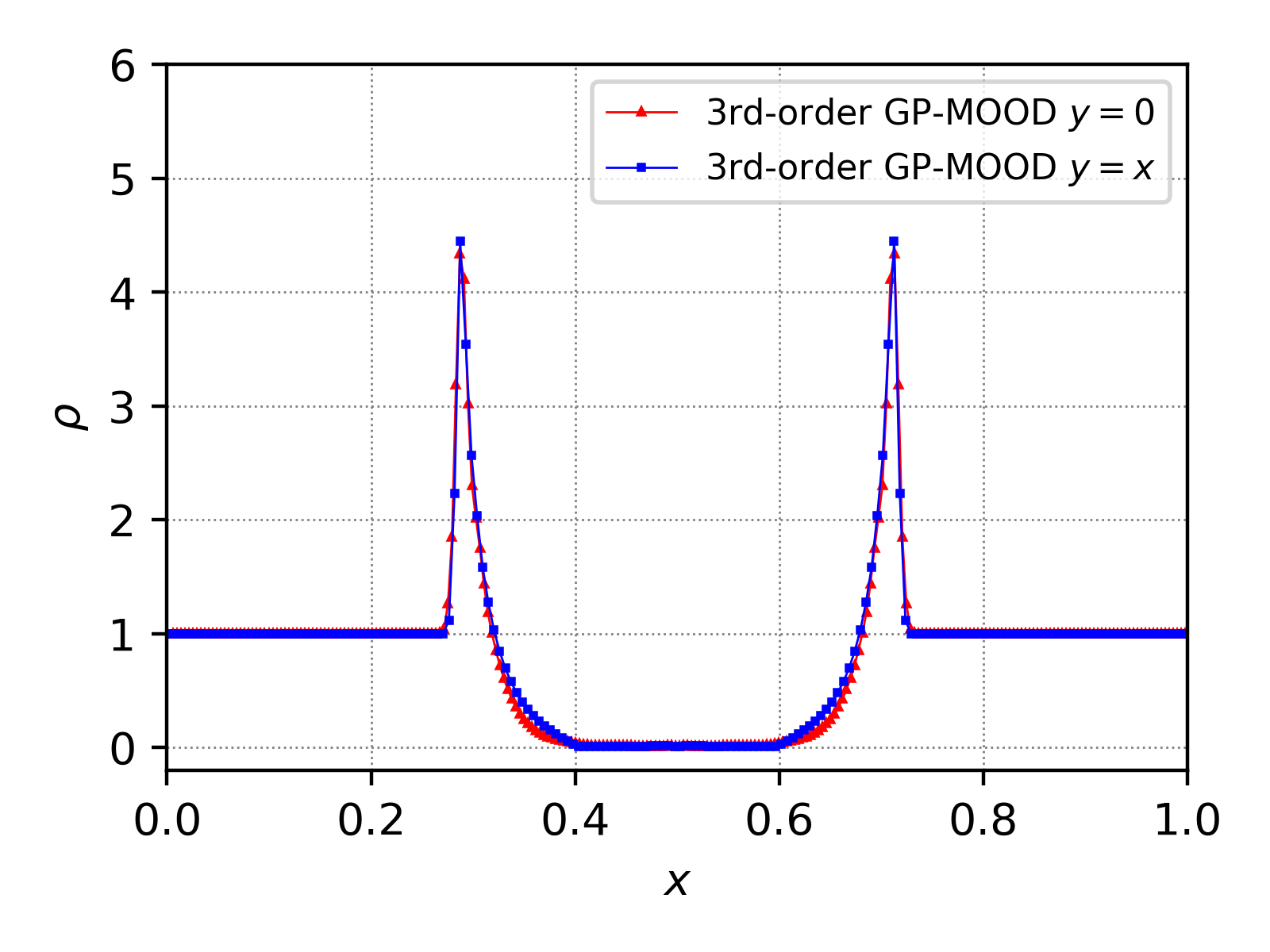}
    \includegraphics[scale=0.8]
    {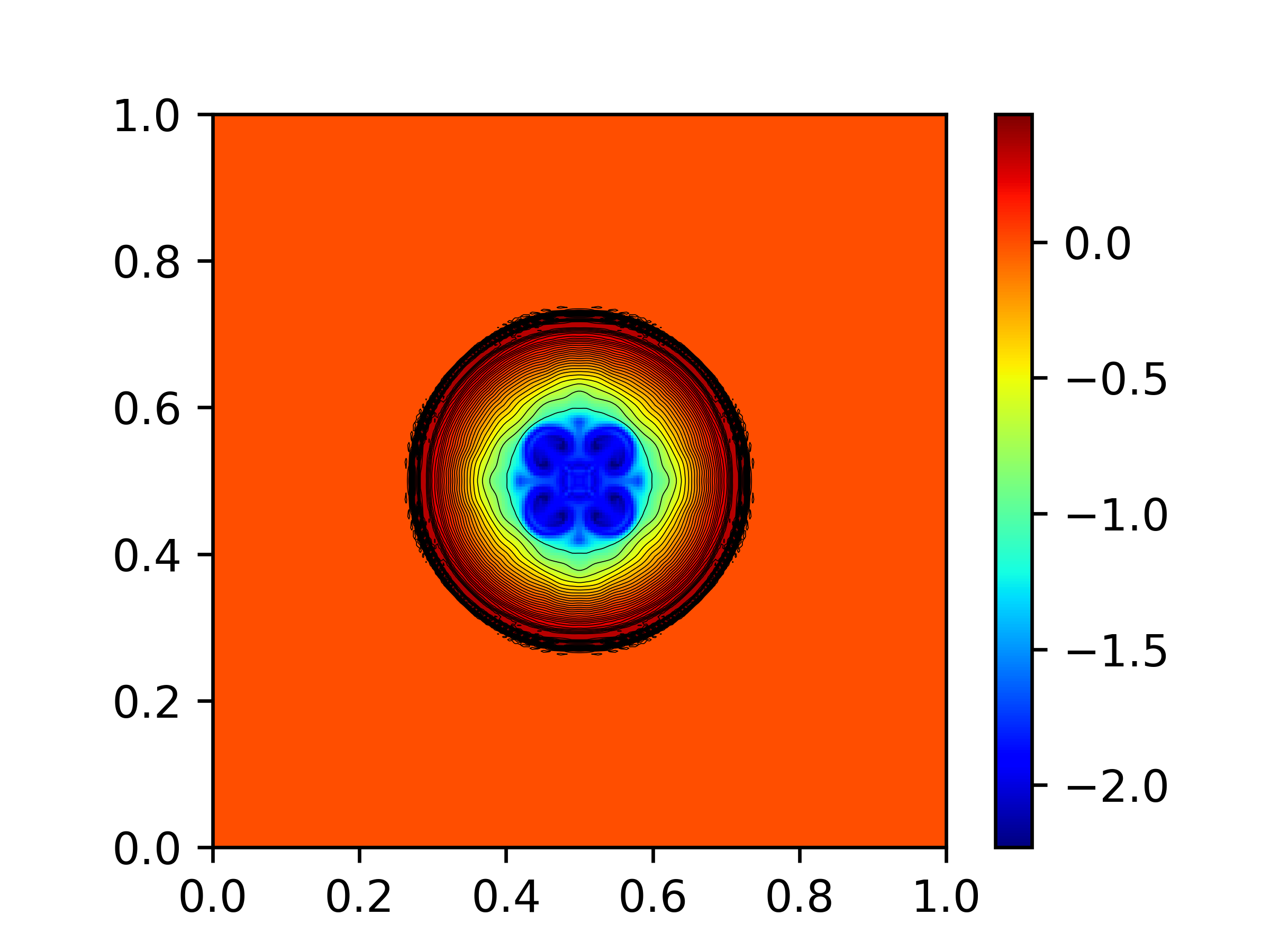}
    \includegraphics[scale=0.75]
    {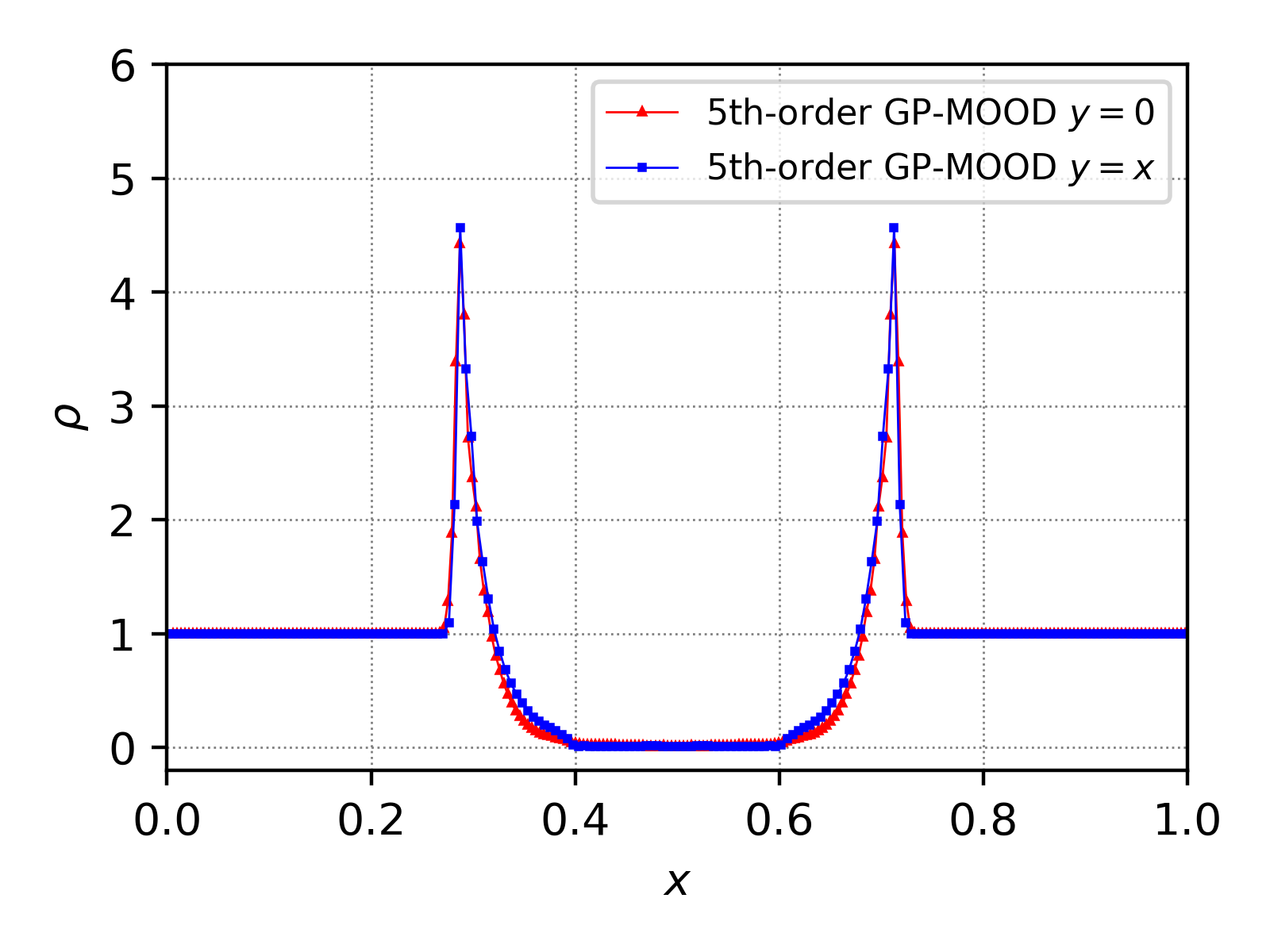}    
    \includegraphics[scale=0.8]
    {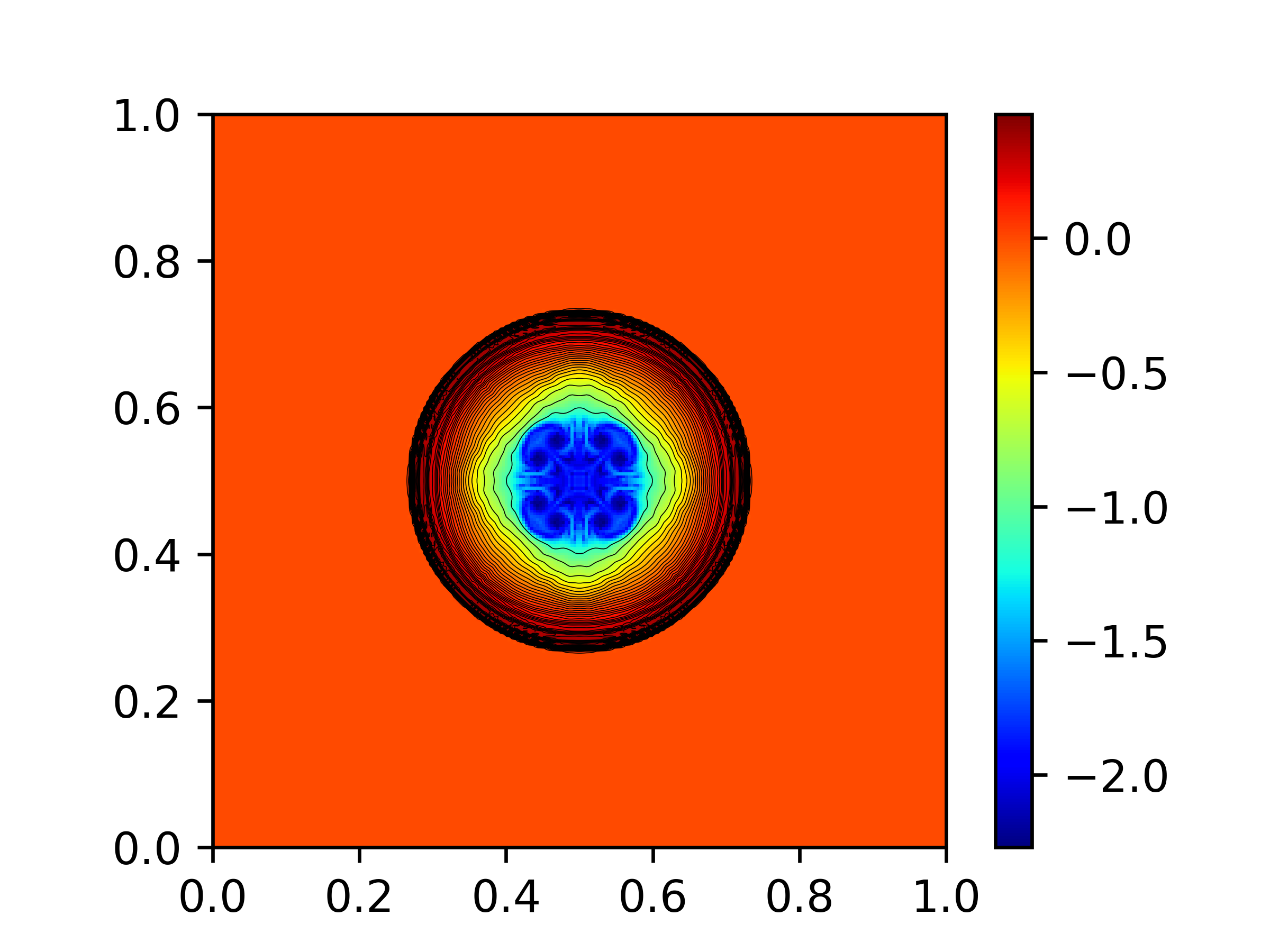}
    \includegraphics[scale=0.75]{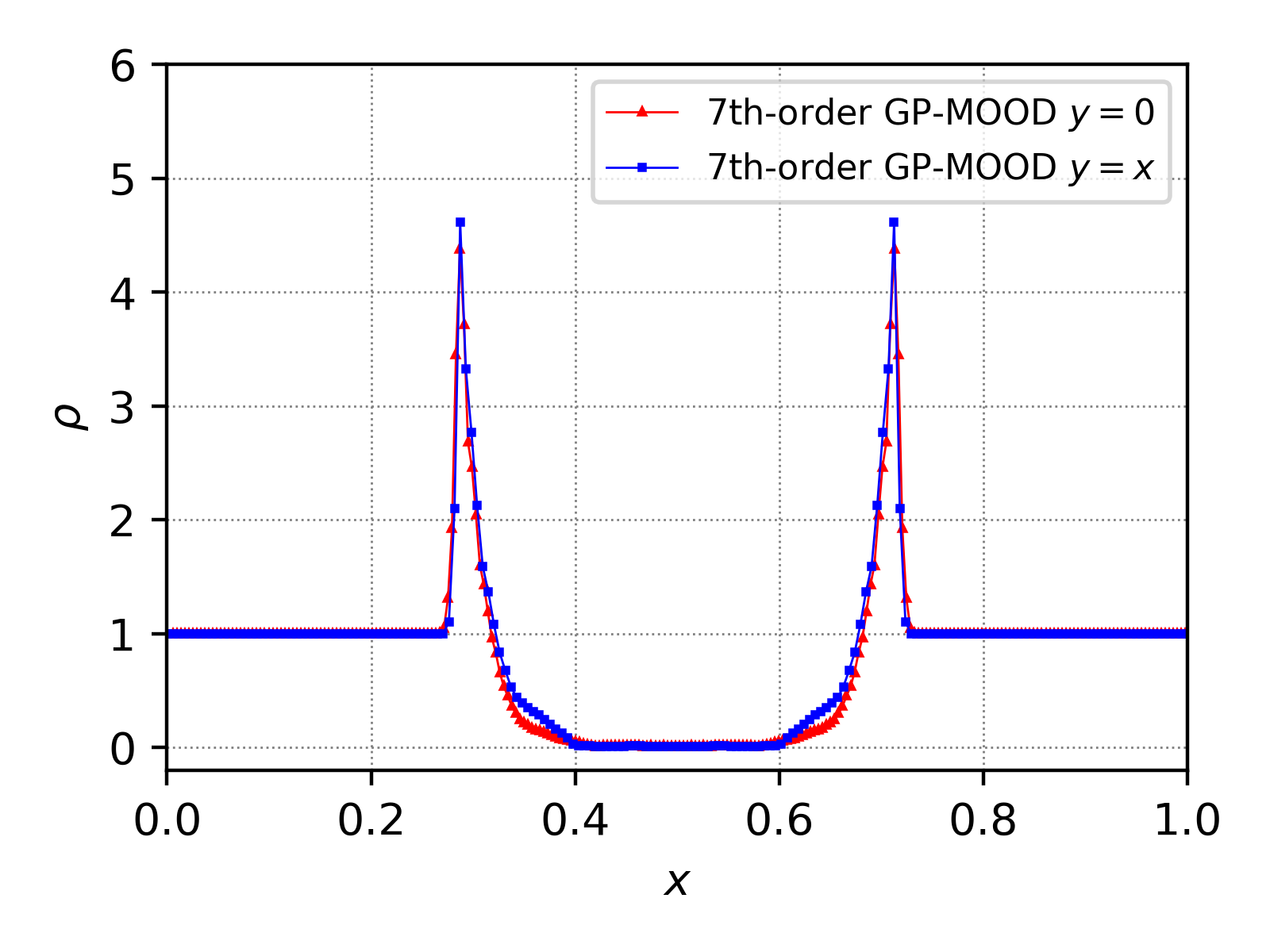}        
 \caption{Density profiles of the Sedov test problem at $t=0.2$ 
 on a $256\times 256$ grid resolution. 
 Three GP-MOOD methods are solved using 
 $\ell = 12 \Delta$ with $\Delta = \dx = \dy$
 and CFL = 0.8.  We over-plot 40 contour lines.
 (\textbf{top}) The 2D density profile (log scale) of GP-MOOD3 (\textbf{left}); 
                      the corresponding sectional profiles (linear scale)
                      along the $x$-axis (red) and the diagonal $y=x$ (blue)
                      (\textbf{right}).
 (\textbf{middle}) The solutions of GP-MOOD5.
 (\textbf{bottom}) The solutions of GP-MOOD7.
}\label{fig:sedov_2d}
\end{figure}

Three test results are shown in \cref{fig:sedov_2d}, including the density
profiles of GP-MOOD3, GP-MOOD5, and GP-MOOD7.
All three results maintain the spherical evolution of the 
exploding shock satisfactorily. Unlike the failed cases in
\cref{fig:sedov_unstable}, the central region
exhibits the details of the low-density flow structures
as the flow is driven radially outward with the expansion.
The expanding flow feels unavoidably the grid-aligned 
dependency, leaving in the central region
the grid-aligned artifacts along 
the $x$- and $y$-axis, as well as 
the round-shaped diagonally-aligned structures.
We observe that those structures become more expressive
as the accuracy order of GP-MOOD increases.
The sectional density profiles show that
the peak amplitudes along the $x$-axis
and the domain diagonal $y=x$ agree
well with each other.

\subsubsection{The double Mach reflection problem in 2D}\label{sec:DMR}

%
%

The next 2D test problem is the double Mach reflection problem by
Woodward and Colella \cite{woodward1984numerical}.
This test problem has been widely chosen to demonstrate
the code capability of handling a strong shock propagation,
reflections, and subsequent developments of vortical roll-ups,
including our former study on the \textit{a priori} GP-WENO finite
difference method \cite{reyes2019variable}.

The initial configuration sets up an inclined planar Mach 10 shock 
initially positioned on
the left side of the rectangular domain.
$[0,4]\times[0,1]$. The shock front makes a
$60^\circ$ incident angle to the 
reflecting wall at the bottom. 
As the shock propagates to the right side of the domain,
the shock is reflected from the bottom wall, which forms
two Mach stems, two contact discontinuities,
and the dense jet along the bottom reflecting wall.
One of the contact discontinuities emanates from
the triple-point, where the contact discontinuity
meets the unreflected Mach 10 shock and
the Mach stem perpendicular to the bottom wall.
The slip surface of this contact discontinuity
undergoes
Kelvin-Helmholtz instabilities, the amount of which
is a sensitive function of the amount of numerical dissipation
in the algorithm under consideration.
Thereby, the problem is often used as an indication
to quantify each method's numerical dissipation
by monitoring the number of Kelvin-Helmholtz
vortical roll-ups at the end of the simulation.
We run the test until $t=0.25$. 
See \cite{woodward1984numerical} for details of 
the initial and boundary setup configurations.

\begin{figure}[ht!]
    \centering
    \includegraphics[width=12cm, trim={0.4cm 2.3cm 0cm 2.7cm},clip]
    {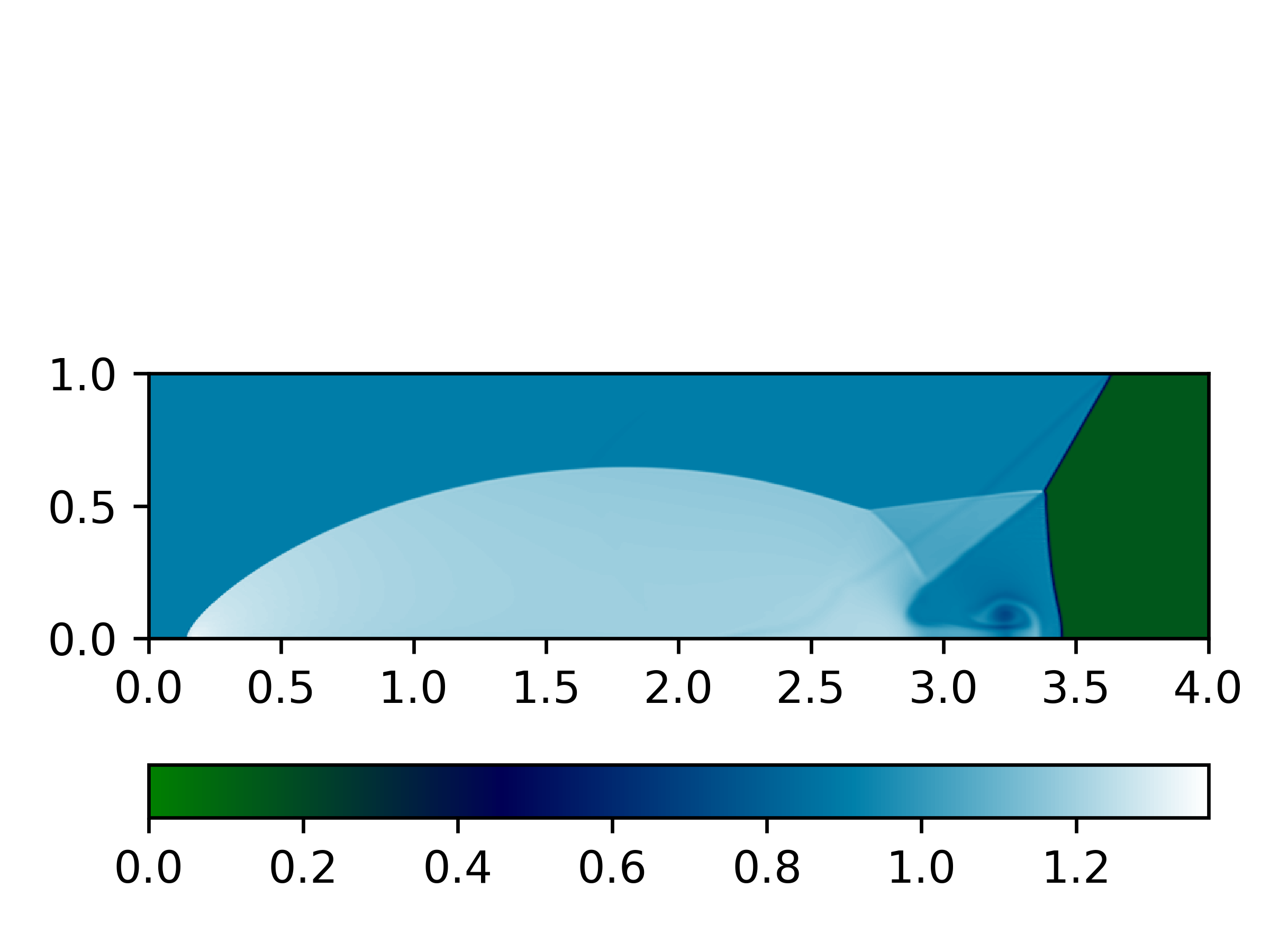}
    \includegraphics[width=12cm, trim={0.4cm 2.3cm 0cm 2.7cm},clip]
    {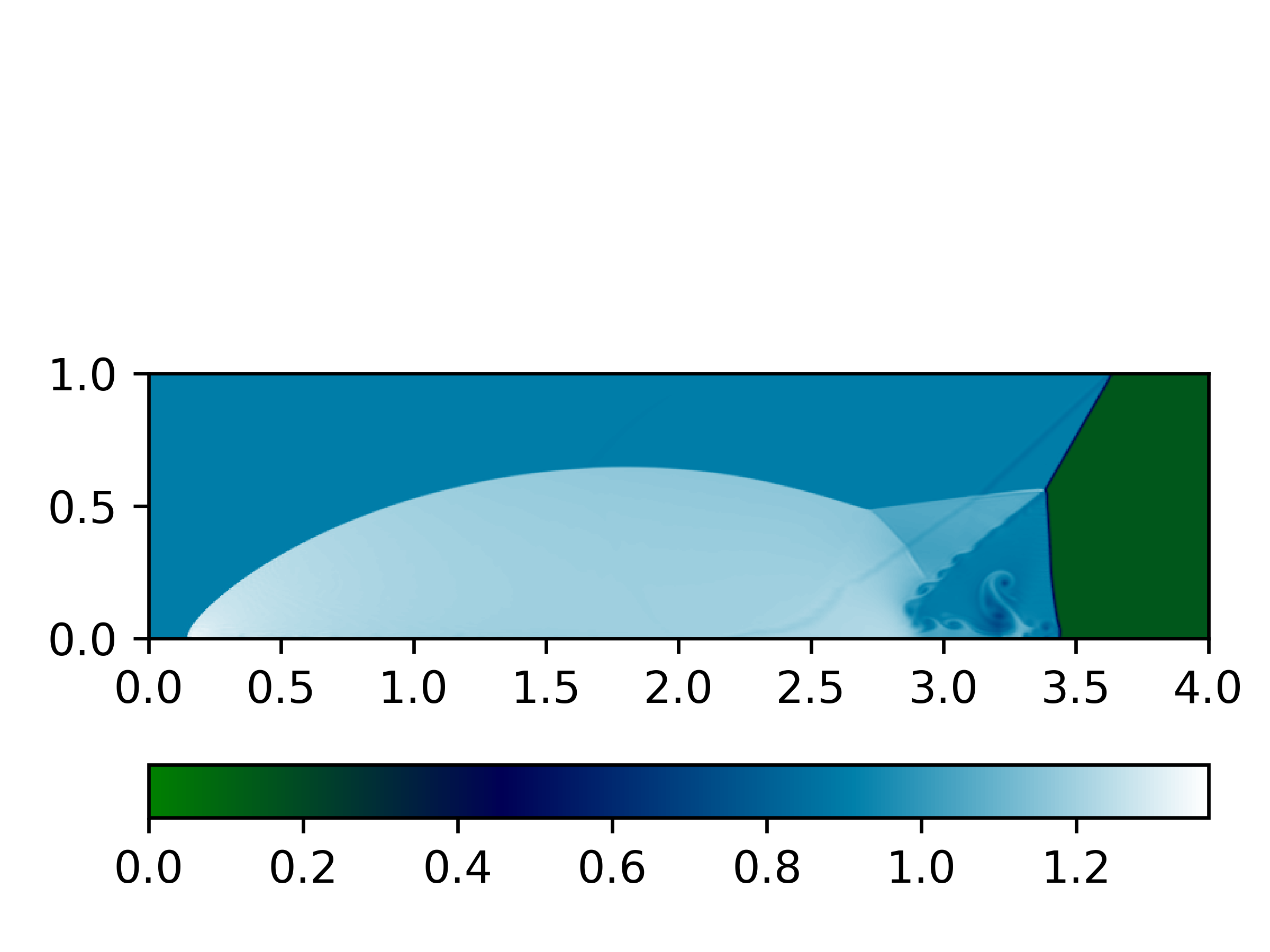} 
    \includegraphics[width=12cm, trim={0.4cm 0cm 0cm 2.7cm},clip]
    {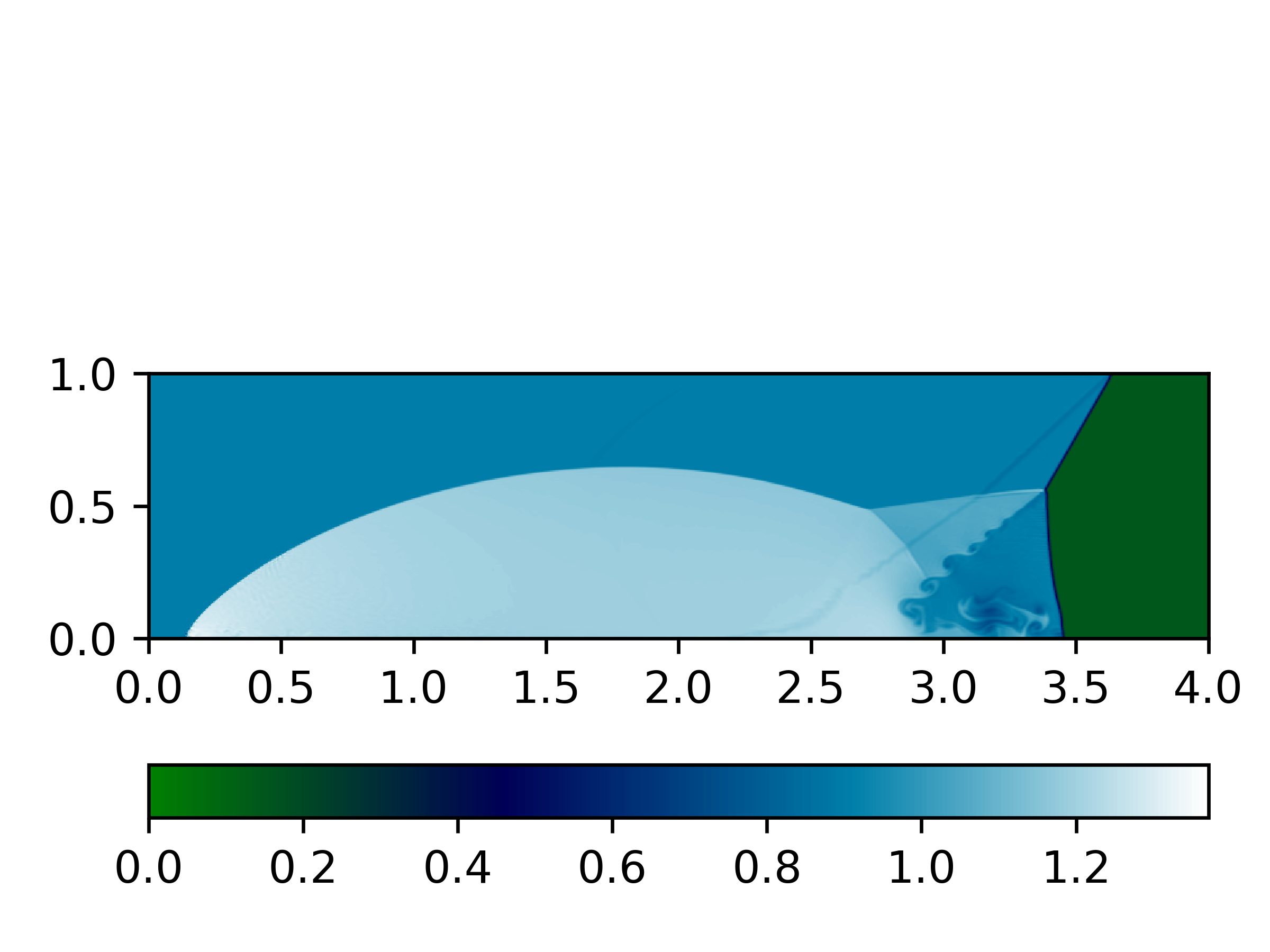}    
 \caption{The 2D double Mach reflection problem on an $800 \times 200$ mesh 
 resolution. The density at $t=0.25$ is plotted using three different GP-MOOD methods,
 all integrated with SSP-RK3, CFL=0.8, and the HLLC fluxes. 
 The choice of $\ell = 12 \Delta$ with $\Delta = \dx = \dy$ is used in all cases.
 (\textbf{top}) The 3rd-order GP-MOOD3. 
 (\textbf{middle}) The 5th-order GP-MOOD5. 
 (\textbf{bottom}) The 7th-order GP-MOOD7.
 }\label{fig:dmr_2d}
\end{figure}

\cref{fig:dmr_2d} displays three density plots of 
GP-MOOD3, GP-MOOD5, and GP-MOOD7,
all integrated using SSP-RK3 without any time step
reduction. They are solved on a grid resolution
of $800 \times 200$ with CFL=0.8. We set
the GP hyperparameter, $\ell = 12 \Delta$,
where $\Delta = \dx = \dy$.
We see that the degree of the 
Kelvin-Helmholtz instability formation
is most pronounced with GP-MOOD7 and
intermediate with GP-MOOD5.
On the other hand, the instability development
is suppressed with GP-MOOD3 at the slip
surface without any sign of vortical roll-up.
The results can also be compared with 
our previous results obtained by
the \textit{a priori} GP-WENO finite
difference method, see Fig. 15
in \cite{reyes2019variable}.
On the same grid resolution, it is seen that
the roll-ups are more evident with GP-MOOD3 
than the result previously predicted by the 5th-order
GP-WENO method therein, namely,
panel (c) of Fig. 15 in \cite{reyes2019variable}.
Although this comparison may not be a
direct point-to-point lateral comparison,
given the two different methods in comparison are
the finite difference method (GP-MOOD5) with CFL=0.8 and
the finite volume method (GP-WENO GP-R2) with CFL=0.4,
this comparison demonstrates that there is smaller
numerical dissipation in GP-MOOD5 than GP-WENO GP-R2.
This result can provide a guideline for our future study 
to conduct a more direct comparison
between the two same consistent discretization methods,
e.g., two finite volume methods of  the 
\textit{a posteriori} GP-MOOD5 scheme and a new
\textit{a priori} GP-WENO GP-R2 scheme.
Such a comparison study will shed light on 
understanding how much numerical dissipation is 
present in the two different shock capturing 
paradigms, \textit{a posteriori} and \textit{a priori} methods.

%
%


\begin{figure}[h!]
    \centering
    \includegraphics[height=5cm, trim={0.5cm 0cm 2.3cm 0cm},clip]
    {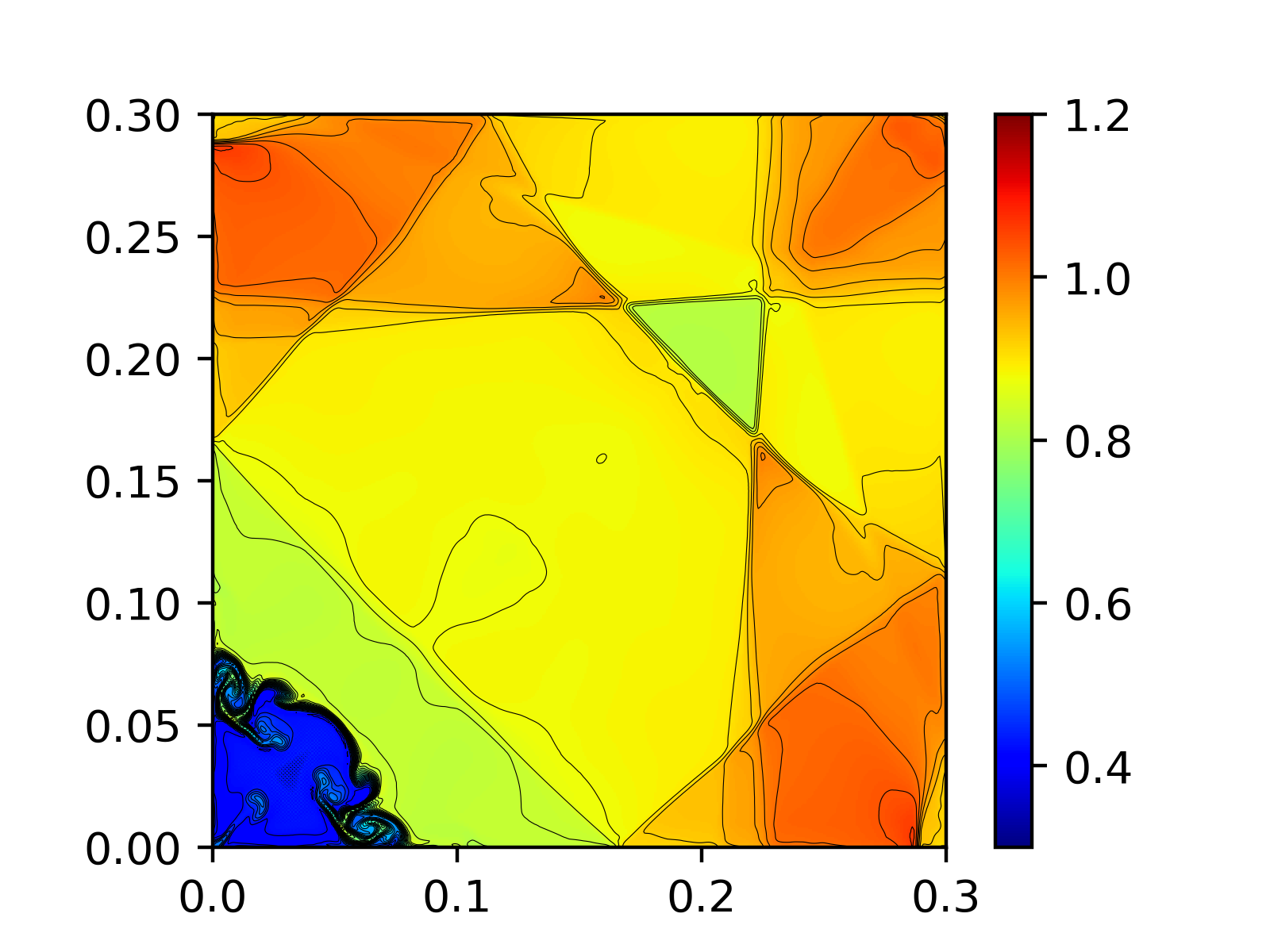}
    \includegraphics[height=5cm, trim={0.5cm 0cm 2.3cm 0cm},clip]
    {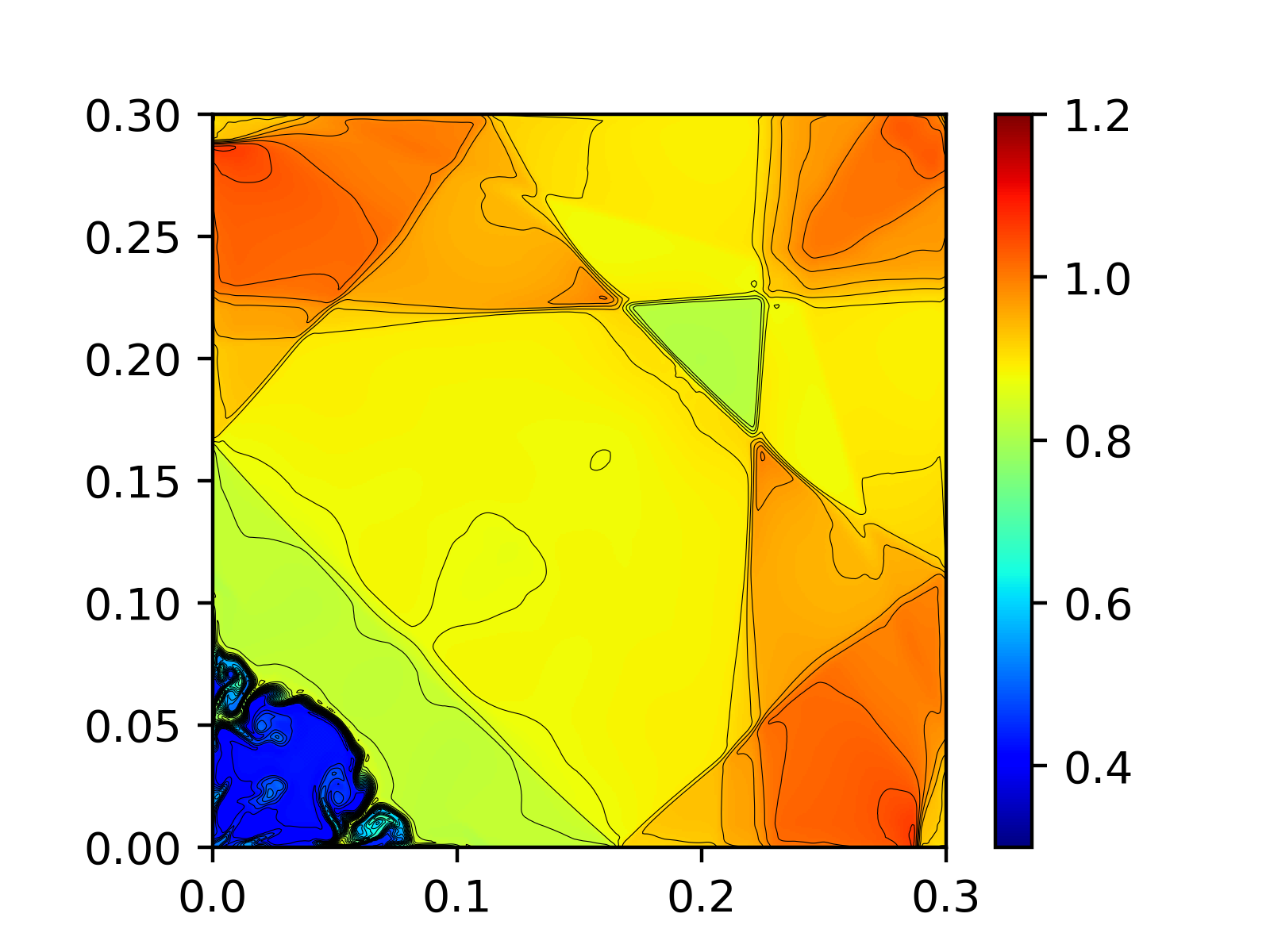}
    \includegraphics[height=5cm, trim={0.5cm 0cm 0cm 0cm},clip]
    {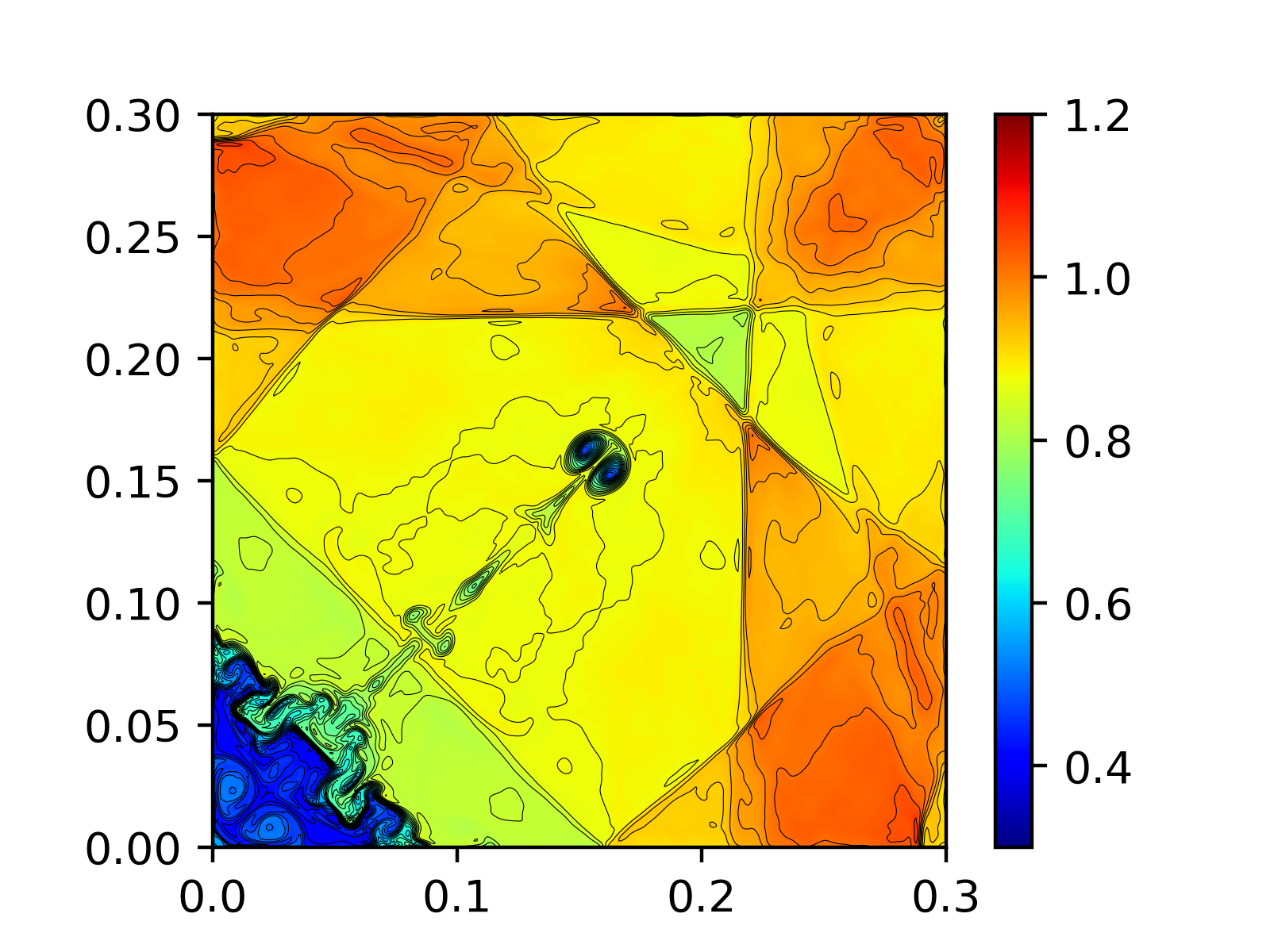}
    \includegraphics[height=5cm, trim={0.5cm 0cm 2.3cm 0cm},clip]
    {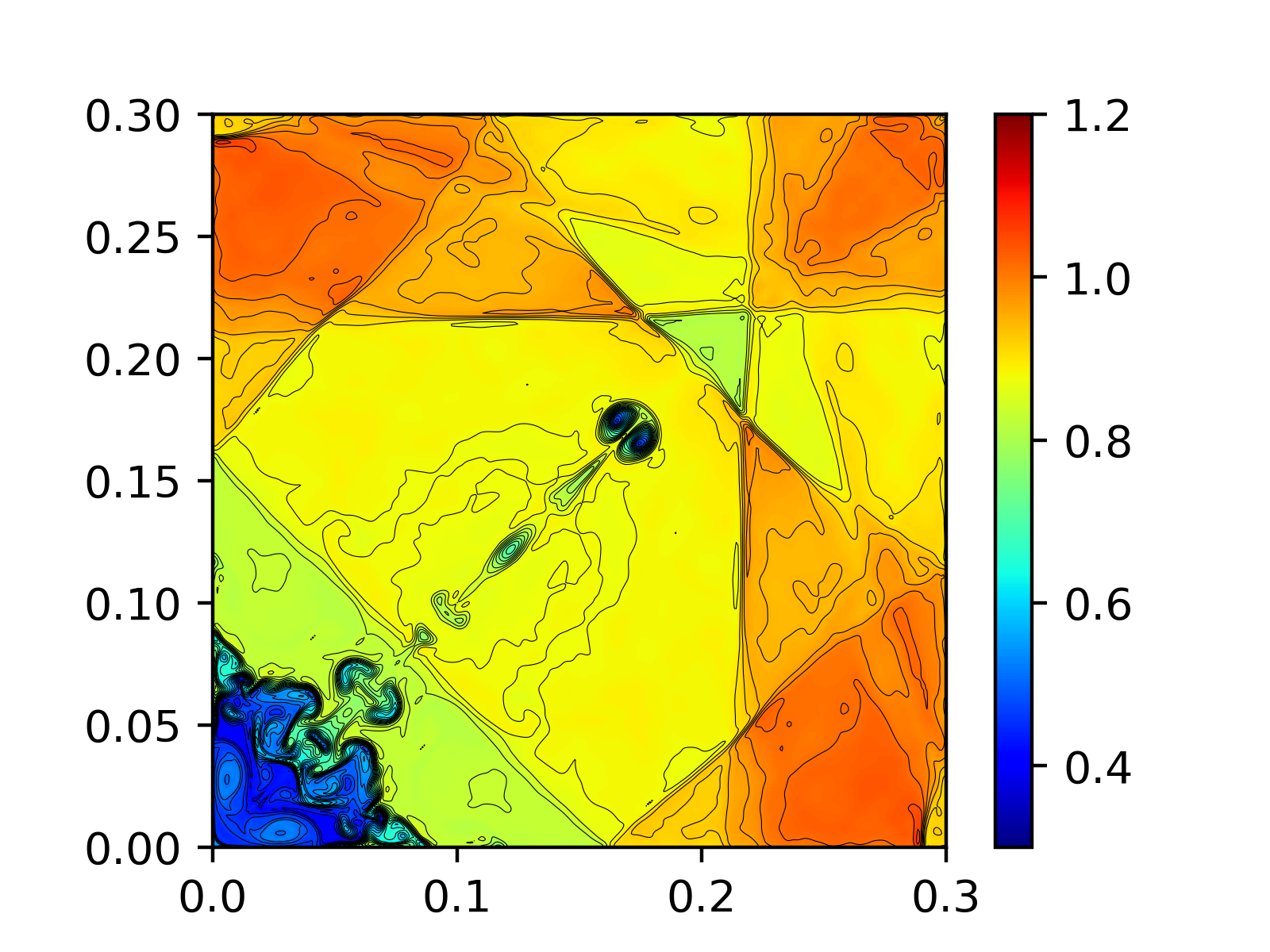}
    \includegraphics[height=5cm, trim={0.5cm 0cm 2.3cm 0cm},clip]
    {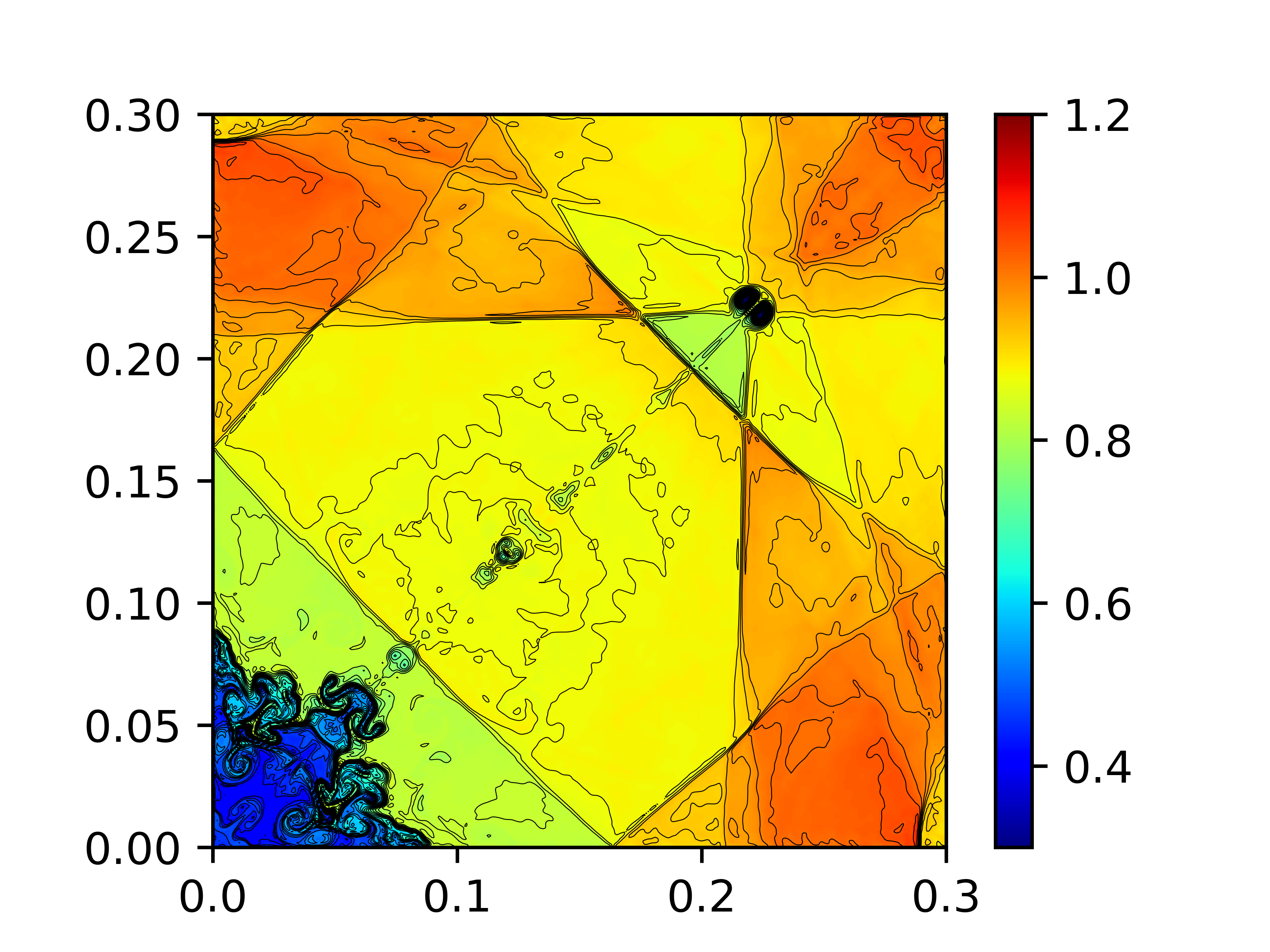} 
    \includegraphics[height=5cm, trim={0.5cm 0cm 0cm 0cm},clip]
    {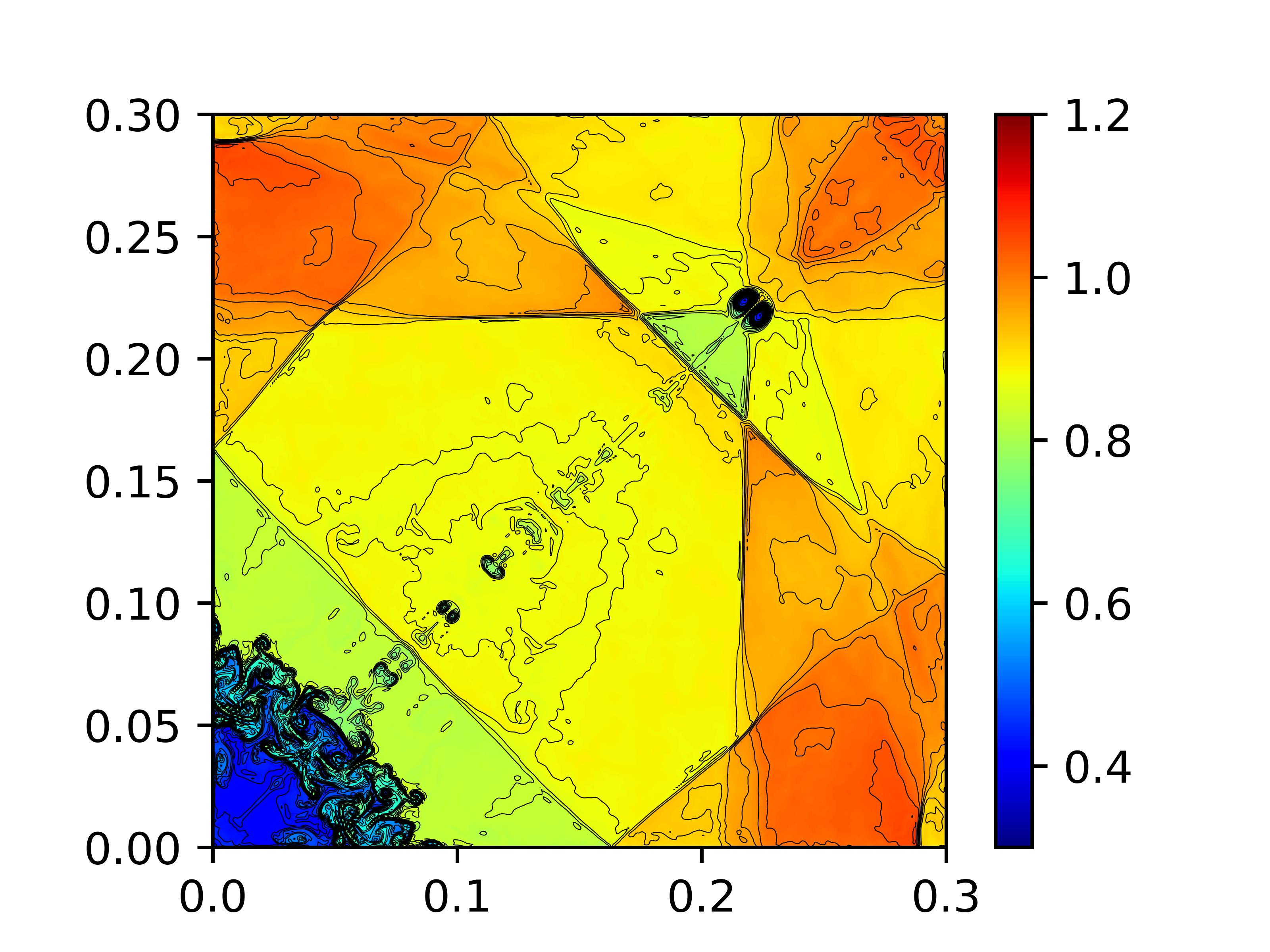} 
 \caption{The density plots at $t=2.5$ of the implosion test in 2D.
 From top left to top right and bottom left to bottom right, the solutions are
 computed using
POL-MOOD3 without CSD,
GP-MOOD3 without CSD,
POL-MOOD3 method with CSD,
GP-MOOD3 with CSD,
GP-MOOD5 with CSD, and
GP-MOOD7 with CSD.
SSP-RK3 is used in all cases without $\dt$ reduction and using CFL=0.8.
 The choice of $\ell = 12\Delta$ with $\Delta = \dx = \dy$ 
 is used in all cases. The contour lines
 are plotted for 40 evenly-spaced density values.
 }
 \label{fig:implosion_2d}
\end{figure}

\subsubsection{The implosion problem in 2D}\label{sec:Implosion}
As for our last benchmark test, we consider the 2D implosion
problem introduced by Hui et al. \cite{hui1999unified} and further simplified
by Liska and Wendroff \cite{liska2003comparison}.
We set up the initial condition by following \cite{liska2003comparison} 
on a square domain, $[0, 0.3] \times [0, 0.3]$. All four boundaries
are reflecting walls, which bounce off a converging
shock wave that is initially launched by the initial planar discontinuity
along $x+y=0.15$. The over-pressurized flow above the line
will send the planar shock wave towards the origin,$(x,y) = (0,0)$.
As it moves, the shock wave interacts with the reflecting walls and 
produces a formation
of a double Mach reflection sliding along the left and the bottom
wall boundaries. This Mach reflection forms
two jets that move into the origin along the boundaries at the exact same time,
ultimately colliding each other simultaneously at $(x,y) = (0,0)$.
This dynamics is 
similar to that of the double Mach reflection test in \cref{sec:DMR}.
The two-jet collision then ejects a new collimated narrow jet from the origin
toward the upward domain diagonal direction.
Due to the reflecting walls, the subsequent shocks
are continuously bounced back into the domain,
creating unceasing nonlinear shock-jet interactions.
These interactions progressively push
the diagonal jet diagonally upward while at the same time
turning its shape into a longer and narrower jet over time.
The mushroom heads, fingers, and the fine elongated filaments
are generated at the diagonal jet surface -- the contact discontinuity
separating two density regions -- as a result of the
Richtmyer-Meshkov instability.
Their morphology is sensitive to numerical dissipation, thereby
can be a good indicator to distinguish differing methods.
Besides, it is crucial to maintain the diagonal symmetry at all times
during the simulation.

In \cref{fig:implosion_2d}, we present six results of the density profiles 
at $t=2.5$ computed on a $400 \times 400$ grid resolution.
We find that the CSD switch introduced in \cref{sec:CDS} 
plays an imperative
role in simulating the diagonal jet dynamics successfully. 
To demonstrate its impact, 
we ran the test with and without the CSD switch.
The test cases in the comparison include, 
in the order from top left to top right and bottom left to bottom right,
the 3rd-order polynomial POL-MOOD3 method without CSD,
the 3rd-order GP-MOOD3 method without CSD,
the 3rd-order polynomial POL-MOOD3 method with CSD,
the 3rd-order GP-MOOD3 method with CSD,
the 5th-order GP-MOOD5 method with CSD, and finally
the 7th-order GP-MOOD7 method with CSD.
All spatial solvers are integrated with SSP-RK3 without any 
time step reduction using CFL=0.8.

As can be seen in the first two panels, 
no diagonal jet is present without CSD
in both of the 3rd-order results
based on polynomial and GP reconstructions. 
On the contrary, the jet propagation
is clearly captured in the same solver calculations with CSD
(the third and fourth panels),
by which numerical dissipation is better controlled
without being too diffusive.
From this experiment, we conclude that
the complete absence of the jet dynamics
in the first two results is related to the excessive
amount of numerical dissipation in the conventional
MOOD approach
\cite{clain2011high,diot2012improved,diot2013multidimensional,diot2012methode},
which suppresses the anticipated jet formation 
at the early stage of the evolution.
Activating the CSD check addresses this issue
and imposes the proper amount of diffusivity in the flow.
Apart from the presence of the jet, its shape is also interesting. 
In all four results with CSD the jet retains
exact symmetry, a key indication that has been
widely adopted to test each method's
capability of maintaining reflective symmetry to 
machine precision
\cite{persson2013shock,stone2020athena++,lee2021single}.
We also observe that the distance the jet head travels
from the origin is related to the amount of numerical
dissipation in the tested methods.
The tip of the mushroom head reaches 
$x=y \approx 0.22$ for GP-MOOD5 and GP-MOOD7,
$x=y \approx 0.18$ for GP-MOOD3, and
$x=y \approx 0.17$ for POL-MOOD3.
Similarly, a long-traveled jet is reported in \cite{stone2020athena++}
using the polynomial-based PPM method.

\subsection{Highly compressible flow tests with strong Mach jets in 2D}\label{sec:Mach_jets}
This section considers two configurations of highly compressible supersonic astrophysical jets
that propagate into dense and light ambient gases.
\subsubsection{Single Mach 100 light jet}\label{sec:single_mach_jet}
Our first selection of a test problem in highly compressible flows is the extragalactic
astrophysical jet evolution studied by Balsara \cite{balsara2012self}. 
This setup is the second type of the two jet configurations 
tested in \cite{balsara2012self},
called the \textit{Mach 100 light adiabatic jet} that enters 
the ten times denser ambient fluid
from a narrow slit in the domain bottom.
Of the two configurations therein, this light jet problem is considered to be 
more stringent than the first jet problem to simulate, 
where the first jet problem considers the Mach 800
dense jet propagation into the ten times less dense ambient fluid 
\cite{ha2008positive,balsara2012self}.
We report our GP-MOOD results on the more challenging 
light dense jet case only, 
although we were able to run the dense jet problem successfully.
Magnetized versions of the astrophysical jet have been also reported
in \cite{wu2019provably,liu2021new}.

The computational domain is a square outflow box, $[0, 1.5] \times [0, 1.5]$,
except for the narrow slit, $[0.7, 0.8]$ at $y=0$, through which
the Mach 100 jet is injected into the domain via the inflow
boundary condition fixed by the jet condition,
\beq\label{eq:jet}
(\rho, u, v, p)_{\mbox{jet}} = (\gamma, 0, 100, 1), \;\;\; 0.7 \le x \le 0.8 \;\; \mbox{and} \;\; y=0,
\eeq
where $\gamma = 1.4$. The jet's Mach number is hence 100.
The ambient fluid is initialized by
\beq\label{eq:ambient}
(\rho, u, v, p)_{\mbox{ambient}} = (10\gamma, 0, 0, 1).
\eeq
There are a couple of important features to monitor
in this problem, including 
(i) the symmetric evolution of the jet,
(ii) positivity preservation throughout the simulation time, and
(iii) the internal structure of the ``cocoon'' that surrounds the
jet. 

\begin{figure}[htpb!]
    \centering
    \includegraphics[height=5cm, trim={0.5cm 0cm 1.8cm 0cm},clip]
    {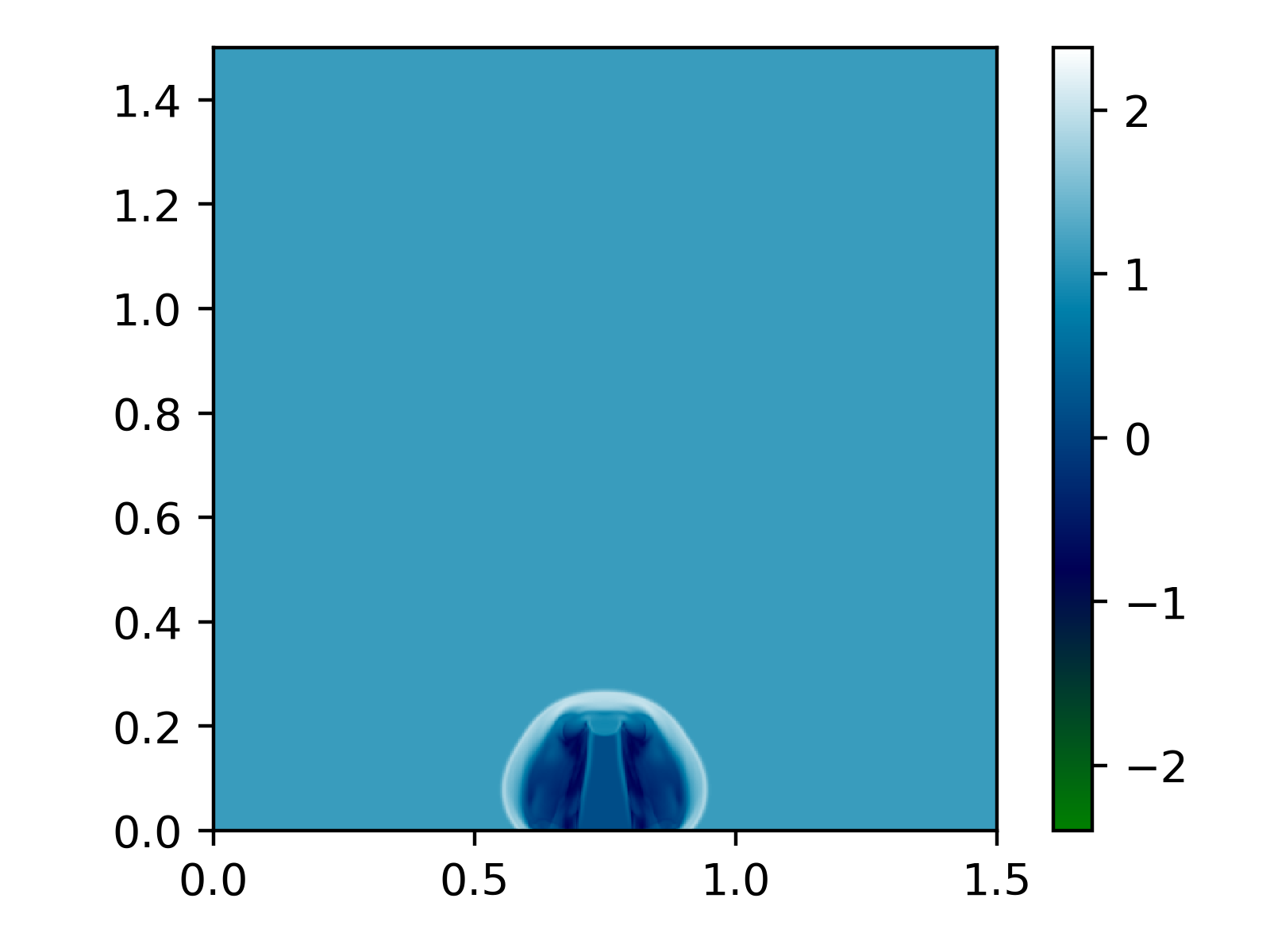}
    \includegraphics[height=5cm, trim={0.5cm 0cm 1.8cm 0cm},clip]
    {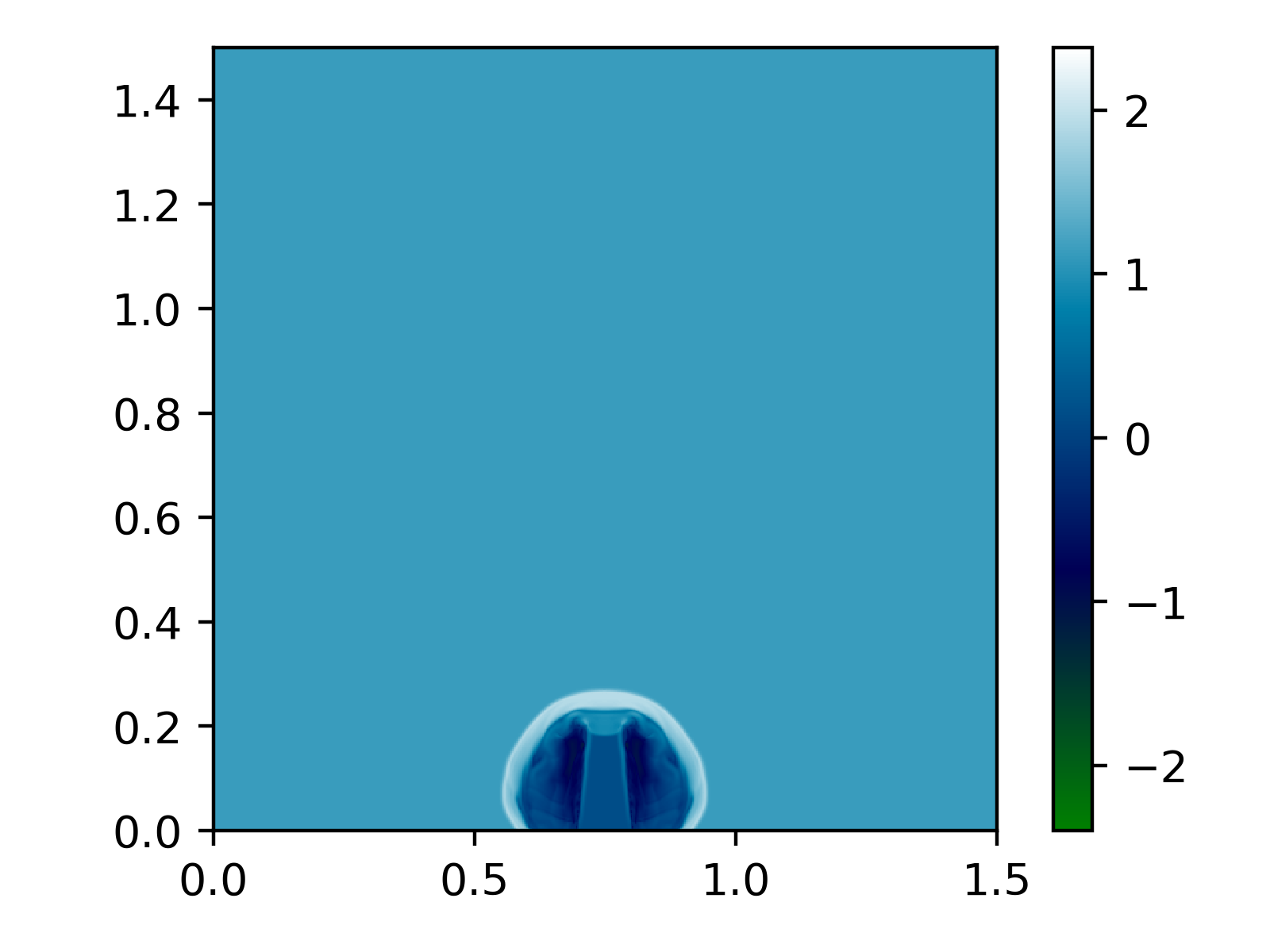}
    \includegraphics[height=5cm, trim={0.5cm 0cm 0.7cm 0cm},clip]
    {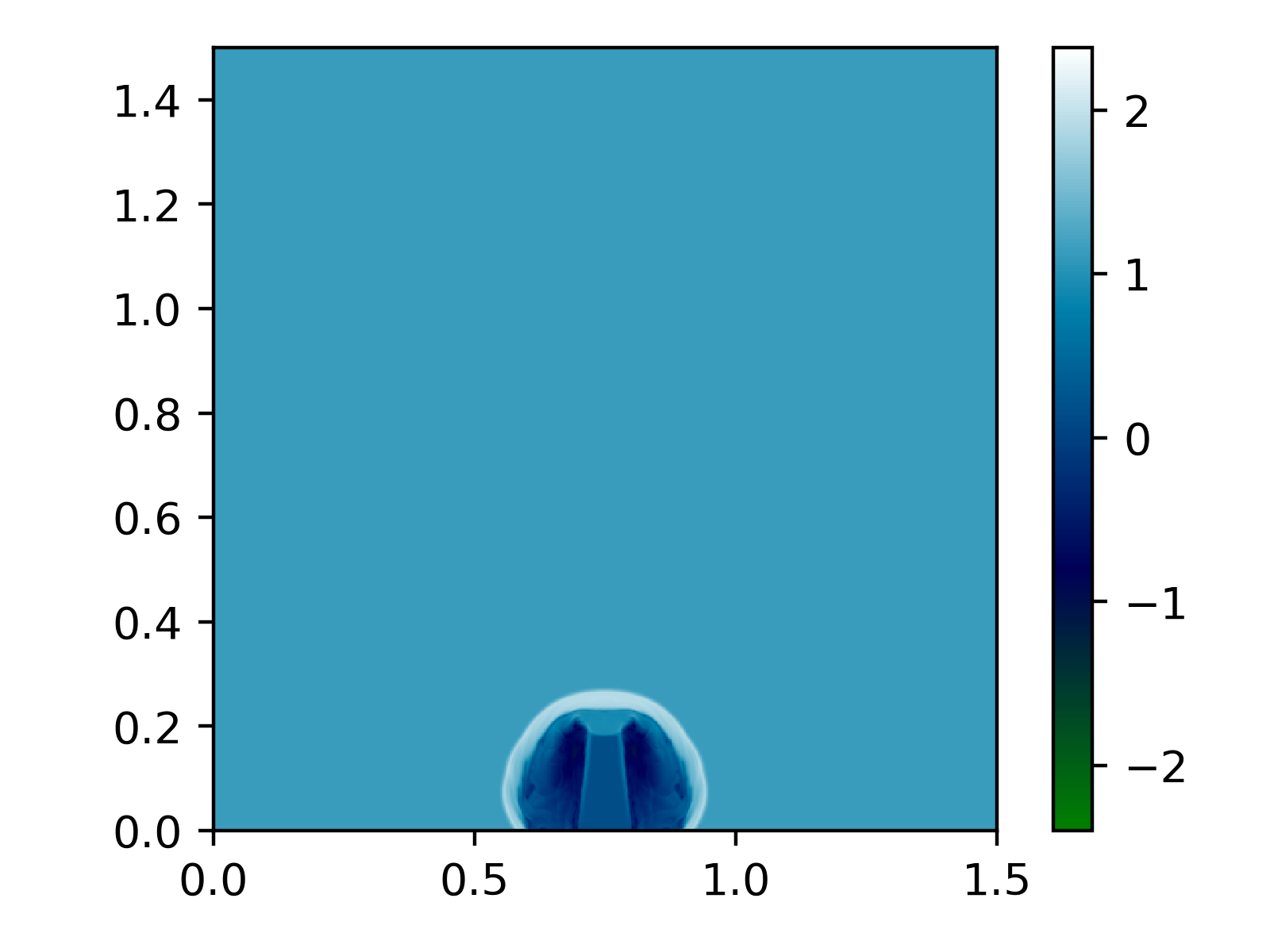}

    \includegraphics[height=5cm, trim={0.5cm 0cm 1.8cm 0cm},clip]
    {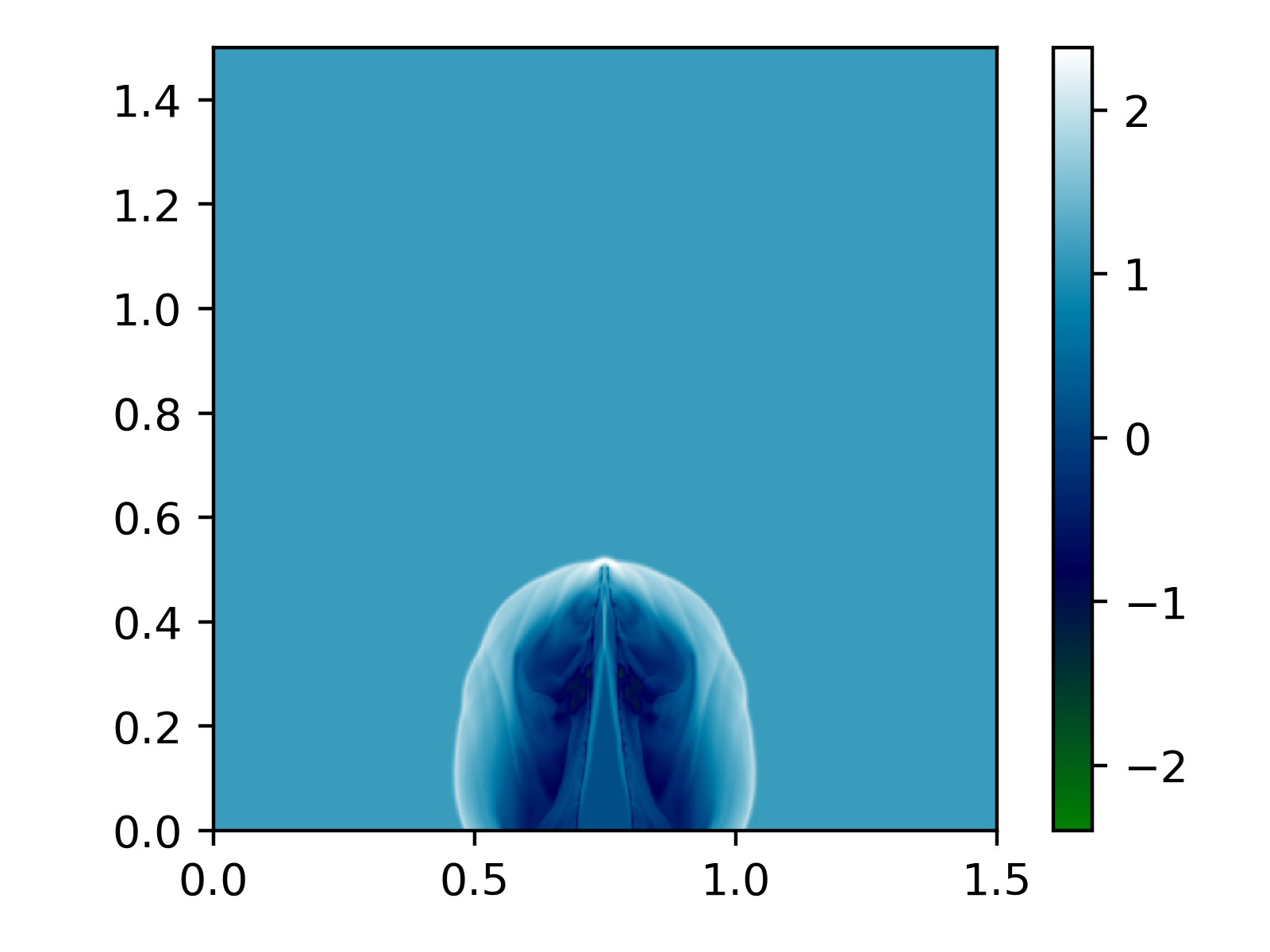}   
     \includegraphics[height=5cm, trim={0.5cm 0cm 1.8cm 0cm},clip]
    {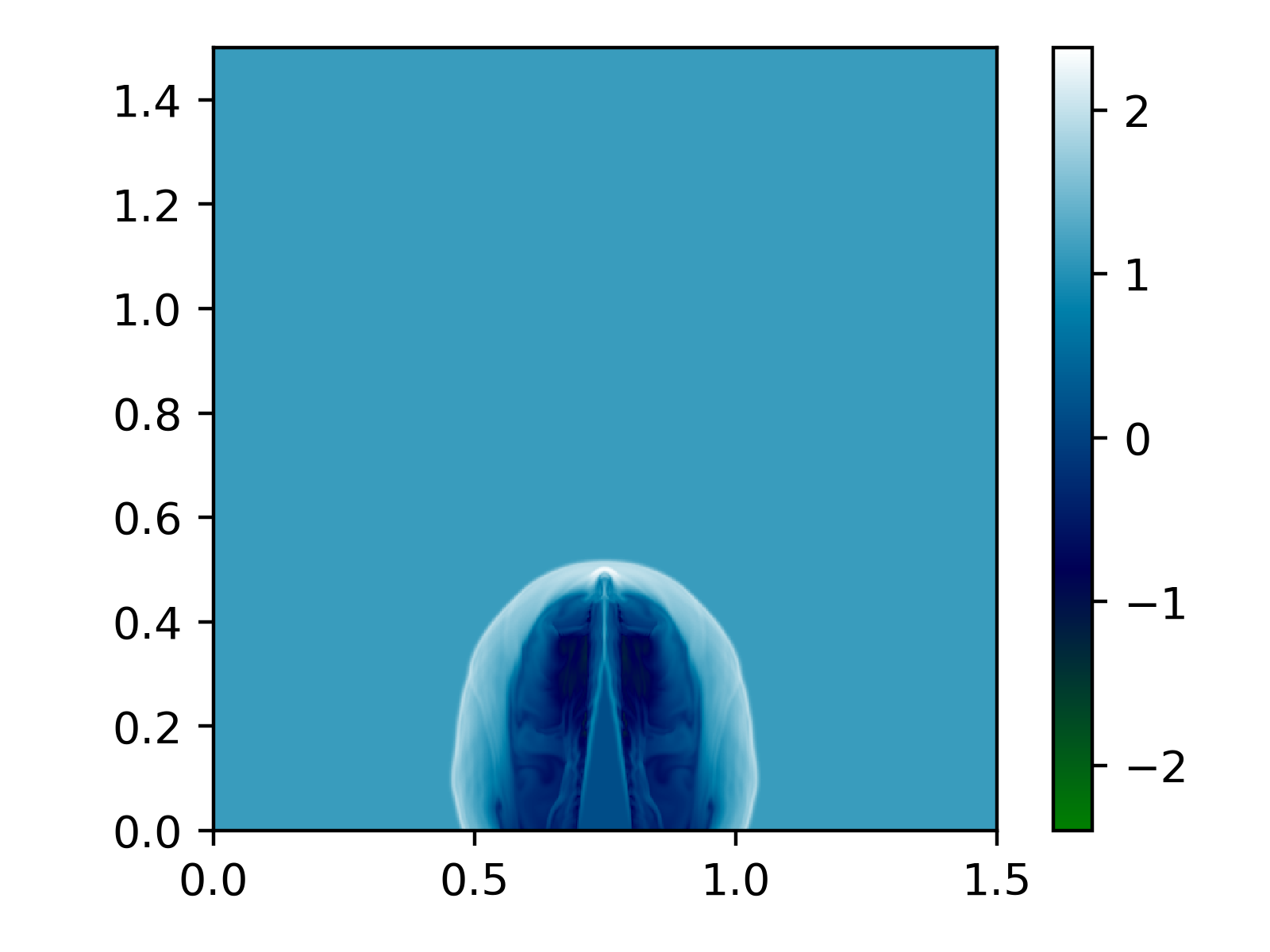} 
   \includegraphics[height=5cm, trim={0.5cm 0cm 0.7cm 0cm},clip]
    {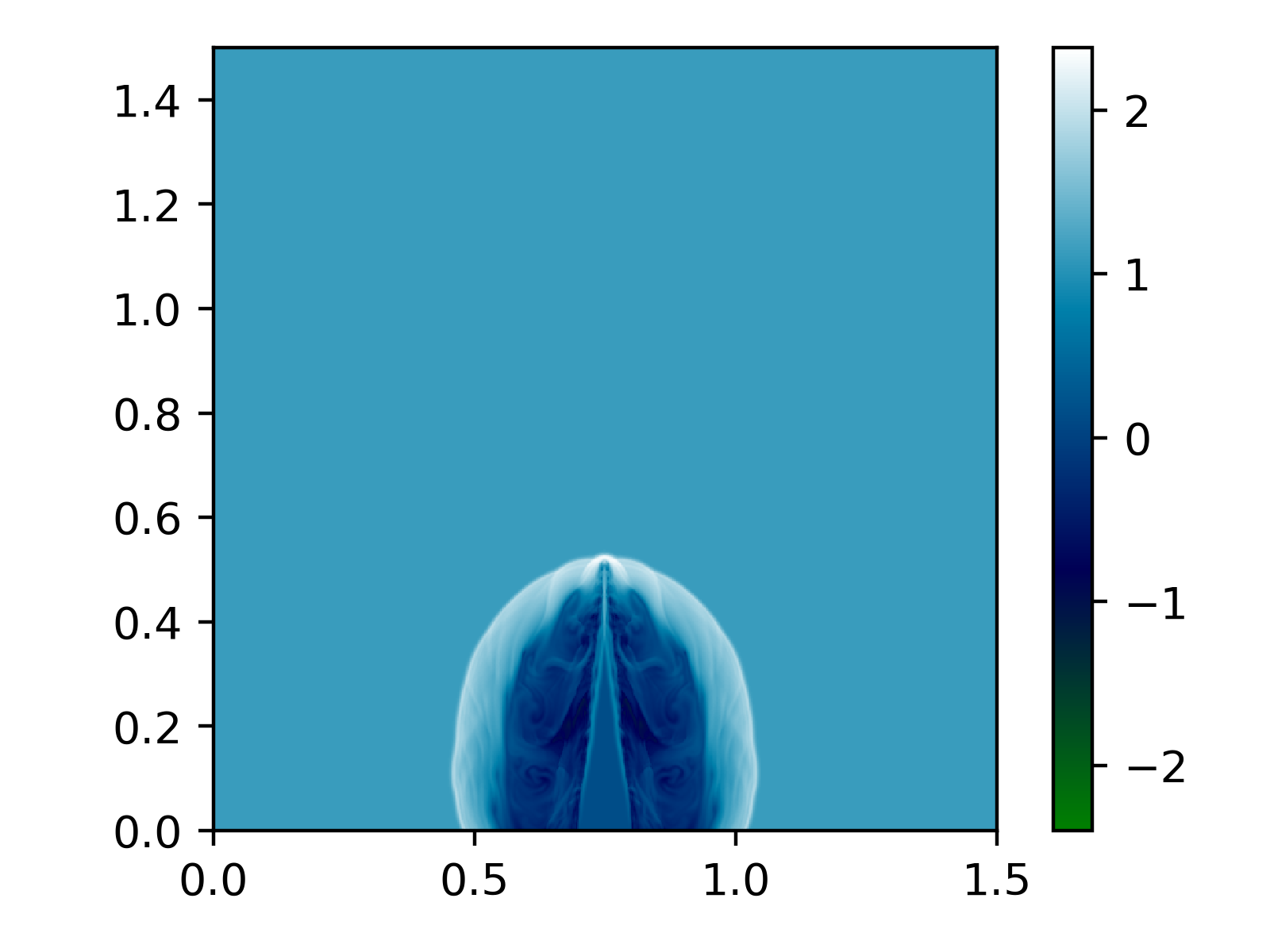}

    \includegraphics[height=5cm, trim={0.5cm 0cm 1.8cm 0cm},clip]
    {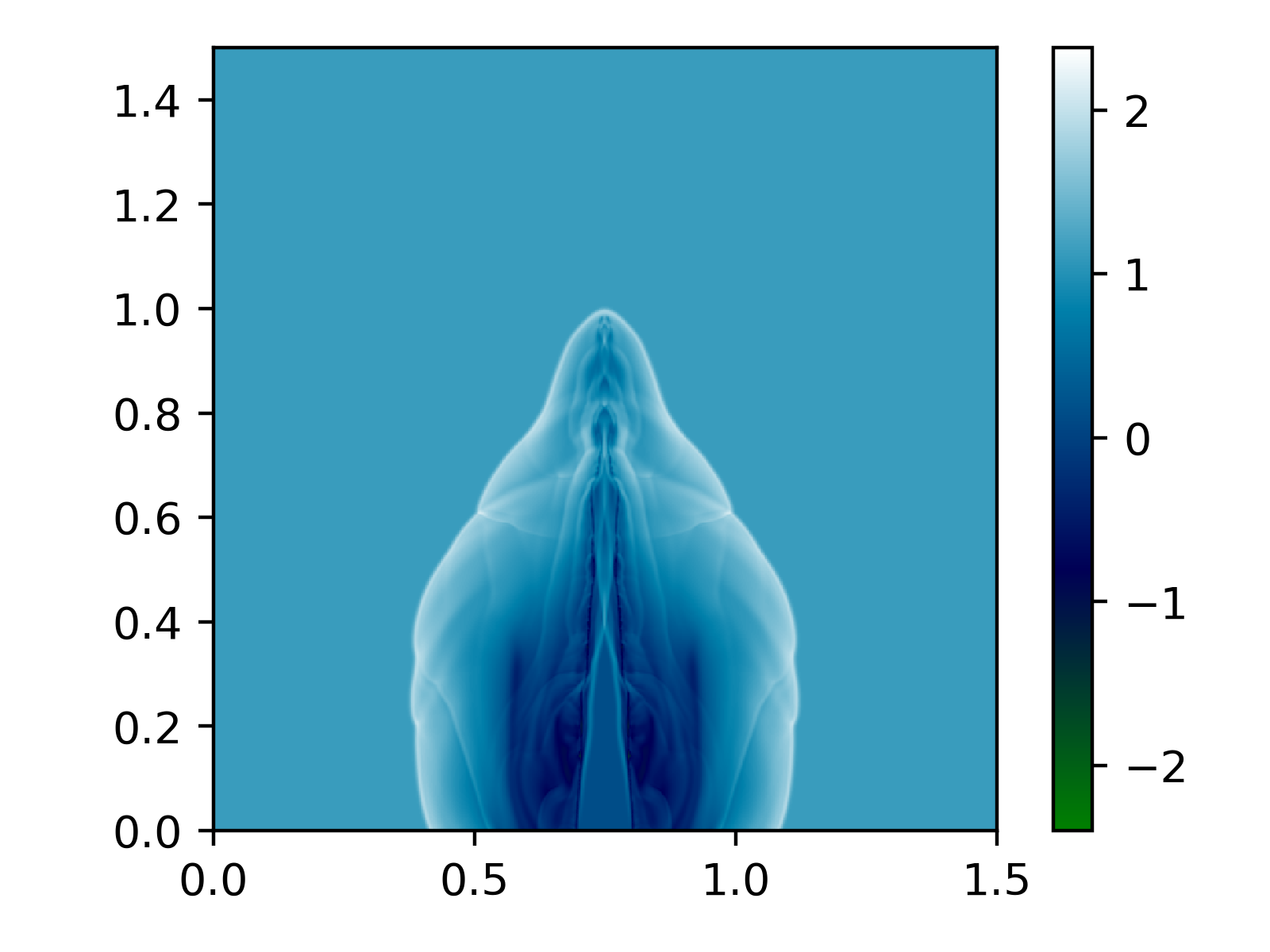}   
     \includegraphics[height=5cm, trim={0.5cm 0cm 1.8cm 0cm},clip]
    {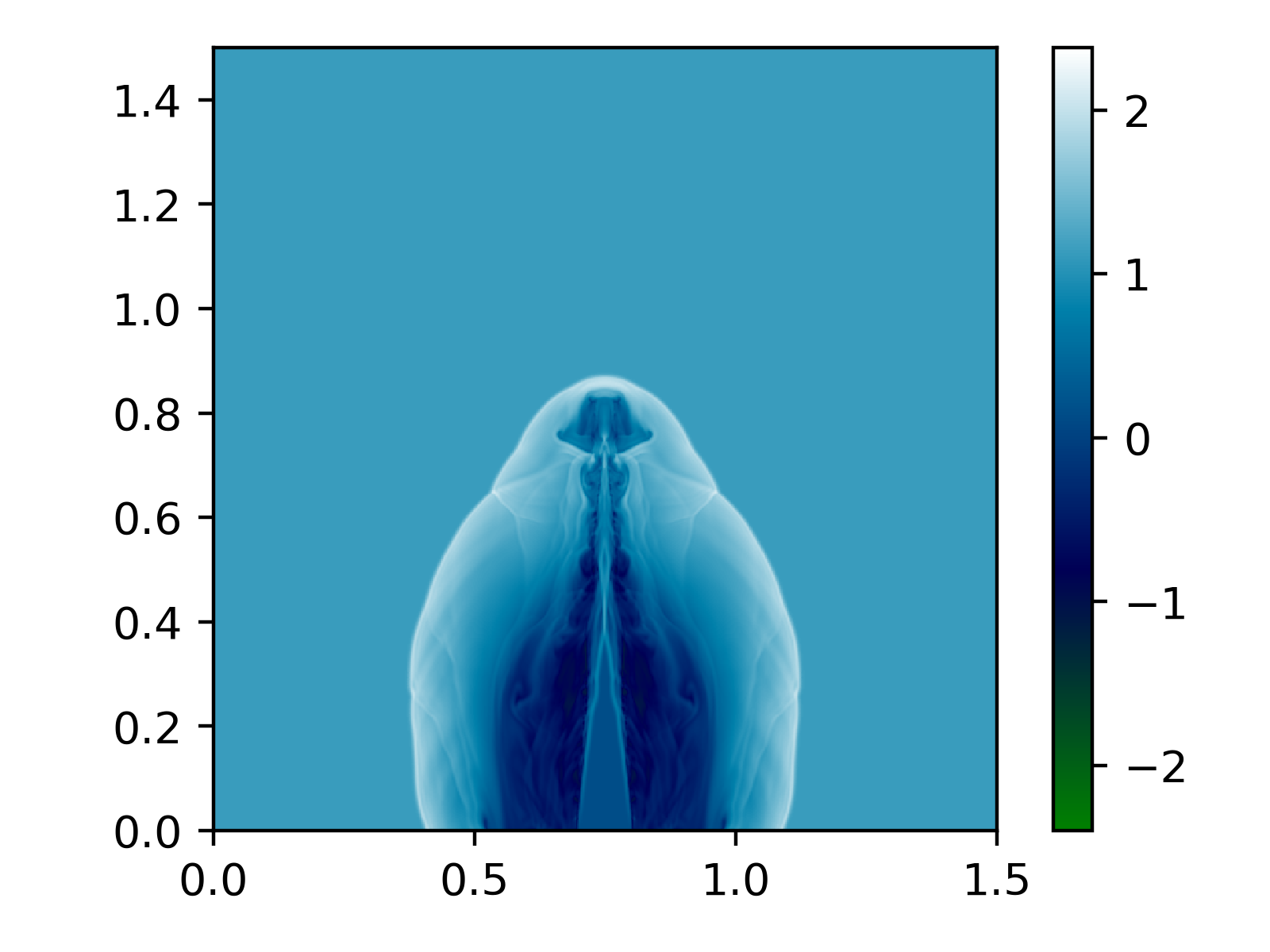} 
   \includegraphics[height=5cm, trim={0.5cm 0cm 0.7cm 0cm},clip]
    {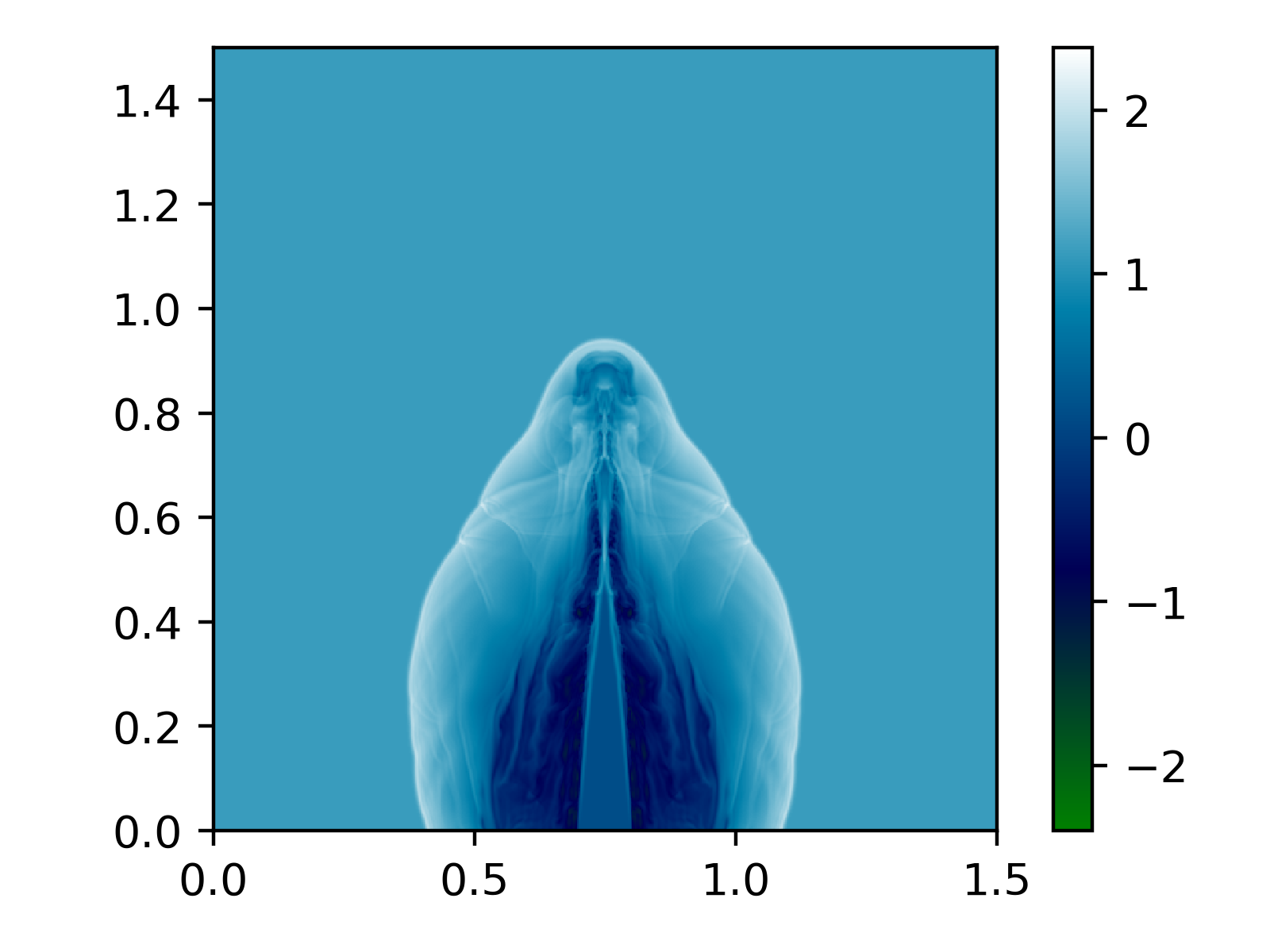}

    \includegraphics[height=5cm, trim={0.5cm 0cm 1.8cm 0cm},clip]
    {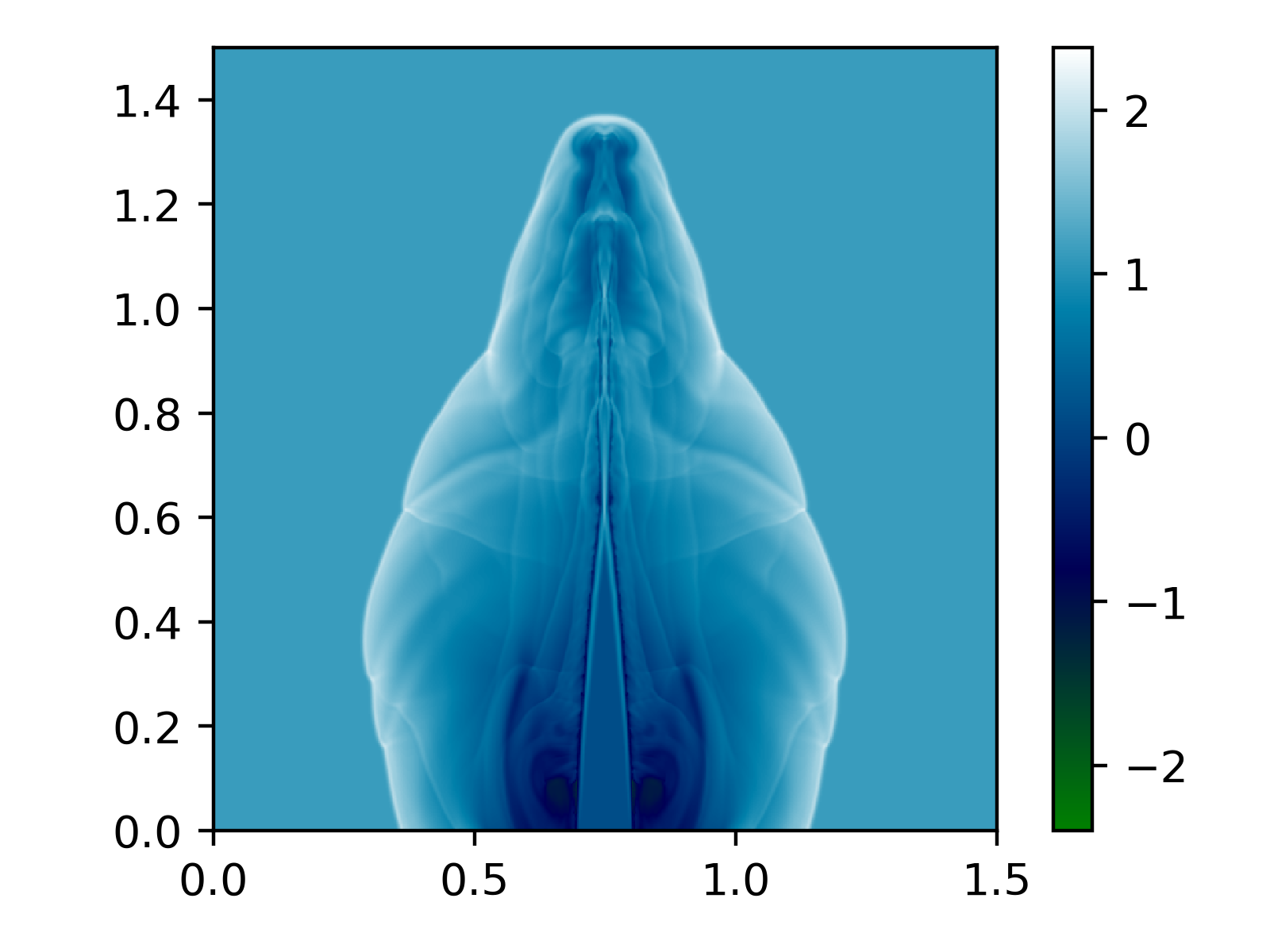}   
     \includegraphics[height=5cm, trim={0.5cm 0cm 1.8cm 0cm},clip]
    {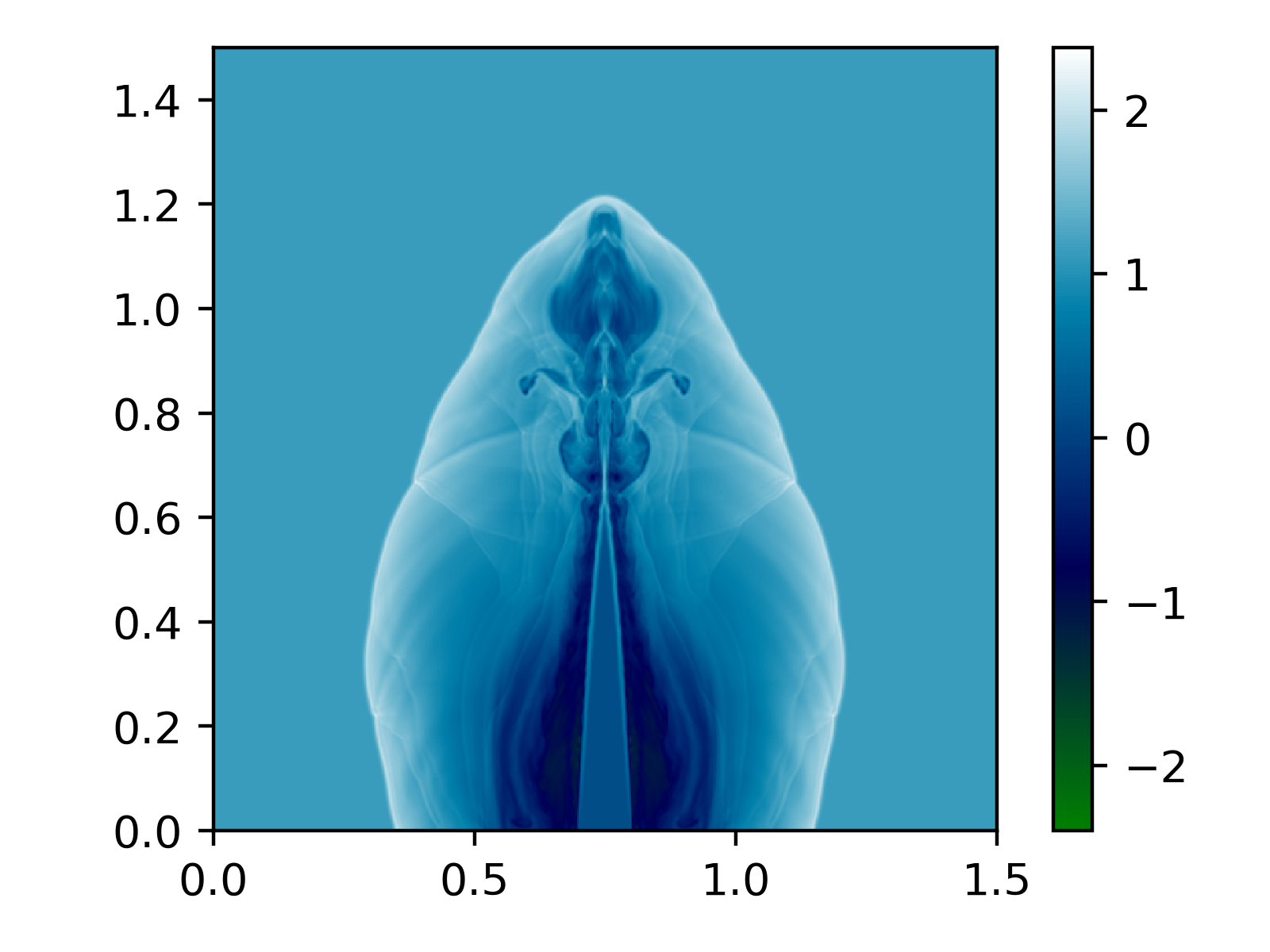} 
   \includegraphics[height=5cm, trim={0.5cm 0cm 0.7cm 0cm},clip]
    {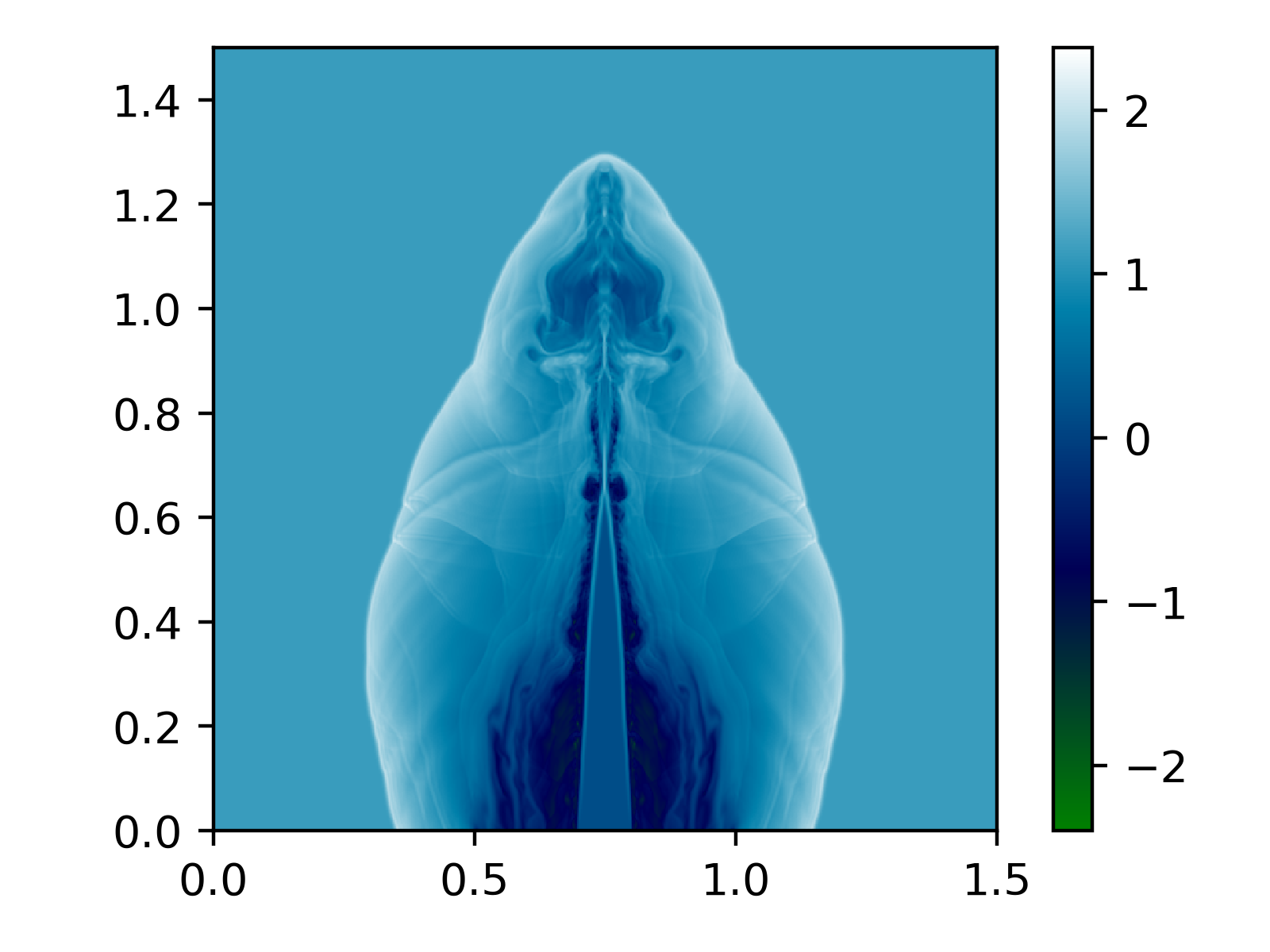} 
    
 \caption{The single Mach 100 light jet problem is resolved on a 
 $600 \times 600$ grid resolution. Density profiles are displayed.
From left to right, the results are computed using
GP-MOOD3, GP-MOOD5, and GP-MOOD7, while
from top to bottom, the densities are plotted at 
$t=0.01, 0.02, 0.03$, and $0.04$.
Except for the use of HLL, the rest of the runtime parameters are the same as before,
including CFL=0.8, $\ell = 12\Delta$ with $\Delta = \dx = \dy$.
 }
\label{fig:single_jet}
\end{figure}

From left to right in \cref{fig:single_jet}, 
each column shows the results using
GP-MOOD3, GP-MOOD5, and GP-MOOD7, which are integrated
until $t=0.04$ on a $600 \times 600$ grid resolution.
Plotted in each row from top to bottom are the density profiles at
$t=0.01, 0.02, 0.03$, and $0.04$, respectively.
As before, all results are computed with 
CFL=0.8 (a relatively very high CFL number for this simulation),
and SSP-RK3 without any time step reduction.
In this simulation and the double jet collision simulation in the next section
we use the HLL solver~\cite{harten1983upstream} 
to mitigate the grid-aligned
carbuncle instability~\cite{quirk1997contribution} 
at the foremost shock front.
The overall jet evolution and the subsequent flow structures
are comparable to the results in Fig. 4 of \cite{balsara2012self}, where 
the reported simulation time was $t=0.03$. 
We doubt that, based on the jet height, the simulation could have
run longer than $t=0.03$, close to $t=0.04$ as we report here.
We were also able to confirm this using FLASH's 
\cite{fryxell2000flash} unsplit hydrodynamics solver, which is
a gas dynamics version of the unsplit magnetohydrodynamics solver
by Lee \& Deane~\cite{lee2009unsplit}; Lee~\cite{lee2013solution}.
%

In the early stage at $t=0.01$, the jet has only propagated
through approximately 15\% of the domain in height, pushing
the denser gas and building up the density and pressure.
Since the ambient gas is denser, the lighter jet's push moving
at Mach 100 creates a highly compressed region at the
nose of the jet, forming a curved shock wave that drapes
down to the bottom. As the jet further pushes
the ambient denser gas, the shock wave expands into
an oval sheath that drapes to the bottom.
The high-pressurized flow just inside the sheath pushes the 
low dense gas between the jet and the sheath 
further inward, squeezing the jet. As a result,
the blunt jet nose turns into a thinner and pointy shape, 
as illustrated in the plots at $t = 0.02$.
This change of the nose shape results in a sudden speed-up
of the nose relative to the shock front, which is seen
as the shock wave penetration of the pointy jet
nose at $t=0.02$. 
Sooner or later, the pointy nose feels the pushbacks
from the denser and highly pressurized shock front
and decelerates, increasing the nose size a bit. This can be
seen in the results at $t=0.03$, where the formerly pointy jet nose at $t=0.02$ 
around 
$y \approx 0.4$ with GP-MOOD3 and GP-MOOD5,
$y \approx 0.5$ with GP-MOOD7,
now has grown to a wider jet at $t=0.03$.

This motion of acceleration and deceleration repeats
and changes the jet into a sausage-like shape.
As a result, the previously oval-shaped shock outer surface at $t\le 0.02$
turns into a more irregular shape, such as an
outer layer of Matryoshka dolls (or Russian dolls),
as plotted at $t=0.03$ and $0.04$.
Another consequence of the 
acceleration-deceleration process
is the shock waves reflected from the kinks at the
shock front, each of which is created when a newly formed pointy
jet nose starts to push the shock front with sudden
acceleration. This push accelerates
the shock front ahead of the jet's nose and stretches
the sheath further, deforming the overall envelope
from a clean oval surface to an elongated Matryoshka doll shape. 
A new pair of kinks is generated at the shock surface, 
symmetrically around the
center of the jet, $x=0.75$.
We see that a rich structure of the shock waves is fully
present in the cocoon region, i.e., the
region between the shock front sheath and the jet.

In contrast, such an acceleration and deceleration process is
not present in the other dense jet configuration \cite{balsara2012self}
since, in this case, the jet is dense enough to retain its width
without experiencing a sausage-like mode. The flow in the 
cocoon region therefore becomes relatively quiescent without
shock waves compared to the light dense jet case considered here.

%

We report in \cref{fig:single_jet} 
that all three GP-MOOD methods successfully 
produce the relevant physics of the flow.
The MOOD approach in all cases strictly enforced the positivity preservation 
without applying any additional positivity preserving mechanism, 
including a fictitious floor value approach 
or the self-adjusting positivity preserving strategy 
proposed in \cite{balsara2012self}.
The acceleration-deceleration process of the jet makes 
the cocoon flow filled with the reflected shock waves.
All three GP-MOOD methods are well capable of capturing
these shock structures without any artifacts, as well as
the sharp profile of the shock front during the run.
All results reveal a great number of fine details 
of the reflected shock waves. There are small-scale 
Kelvin-Helmholtz instabilities driven by 
the upward and the downward flows 
in the low dense valley region between the jet and the 
high dense region where the reflected shocks are present.
The flow attached to the jet experiences 
upward movements as the jet drags the flow.
Trapped between the upward moving jet and the 
the upward moving shock-dominant region, the low-dense
valley feels the downward flow motions, creating
small-scale Kelvin-Helmholtz instabilities. We see that
more abundant small-scale flows are featured as the 
solution accuracy increases from the 3rd-order GP-MOOD3
to the 7th-order GP-MOOD7.

The flow evolution of this problem is truly non-linear, with
continuous interactions between the lighter jet and the denser ambient flows.
Hence, the end results vary with the methods and parameters
used for each run. 
This can be easily seen in different heights of the foremost shock front
and different shapes of the shock front sheath in the three solutions.
Such flows sensitively depend on the sausage-mode evolution of the jet
and the corresponding acceleration-deceleration process.
Interestingly, the tallest height is seen in the GP-MOOD3 solutions,
while the smallest is in the GP-MOOD5.
The mirror symmetry around $x=0.75$ is preserved well in all runs. 

Lastly, we report that the number of cells  that are iterated for
order decrements  in the MOOD loop is less than
3.5\% of the entire cells in each step, implying 96.5\% of the domain
are stably solved by the highest accurate unlimited GP reconstruction method
during the run.

\subsubsection{Double Mach 800 jet collision}\label{sec:double_machJet_collision}

\begin{figure}[htpb!]
    \centering

    \includegraphics[height=5.5cm, trim={0.9cm 0cm 2.35cm 0cm},clip]
    {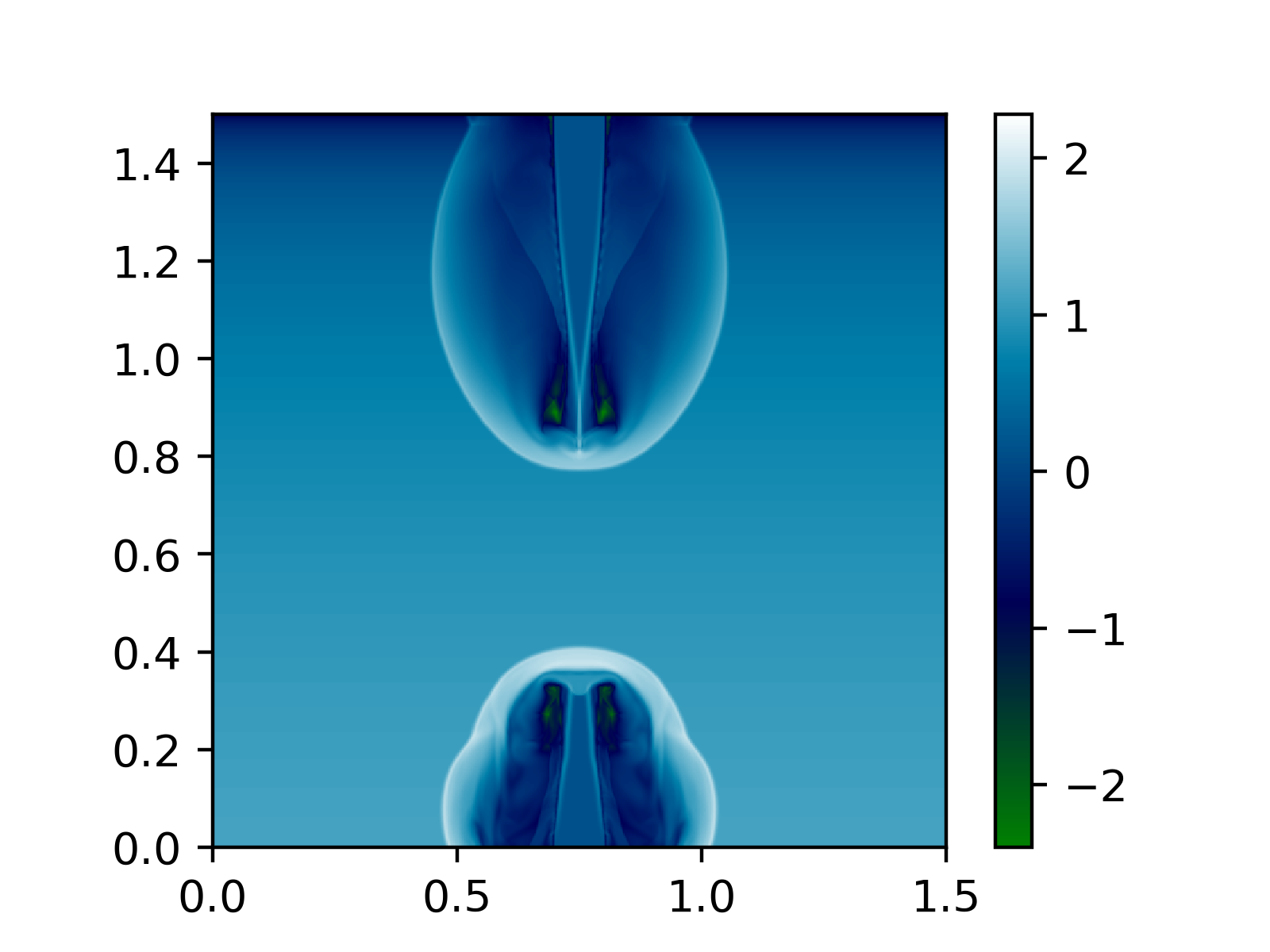}  
     \includegraphics[height=5.5cm, trim={0.9cm 0cm 2.35cm 0cm},clip]
    {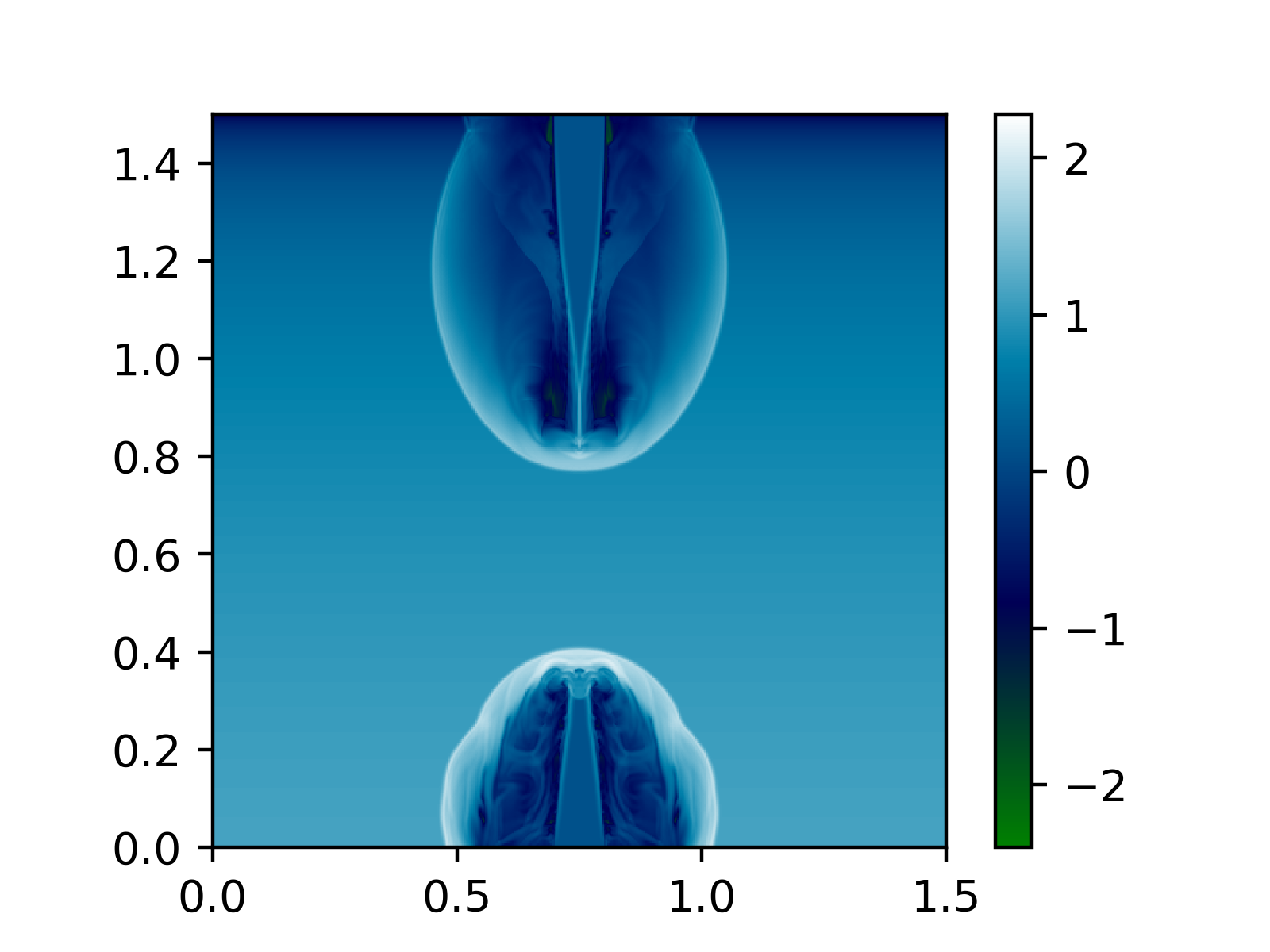}
   \includegraphics[height=5.5cm, trim={0.9cm 0cm 1.15cm 0cm},clip]
    {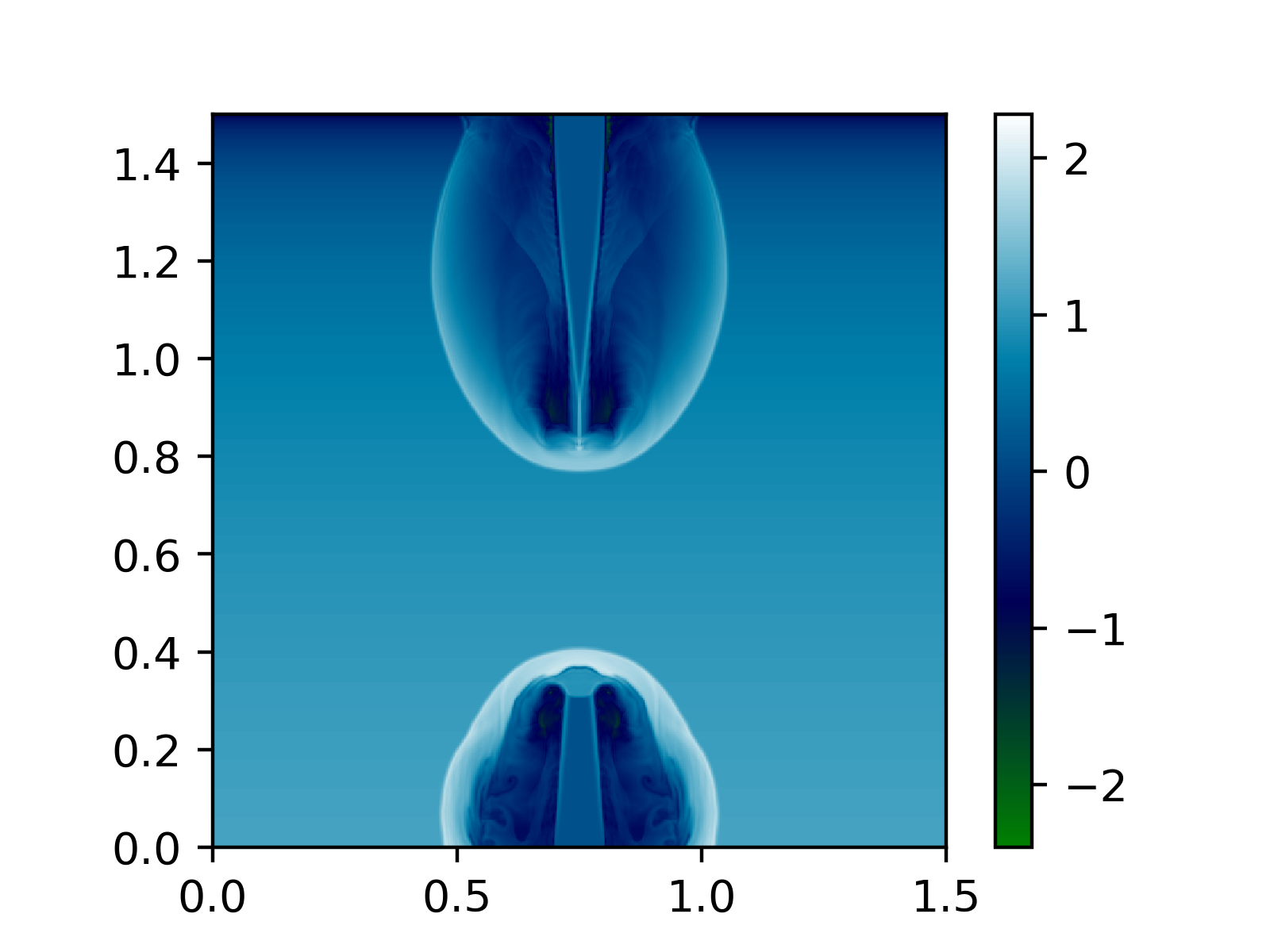}
     \\

    \includegraphics[height=5.5cm, trim={0.9cm 0cm 2.35cm 0cm},clip]
    {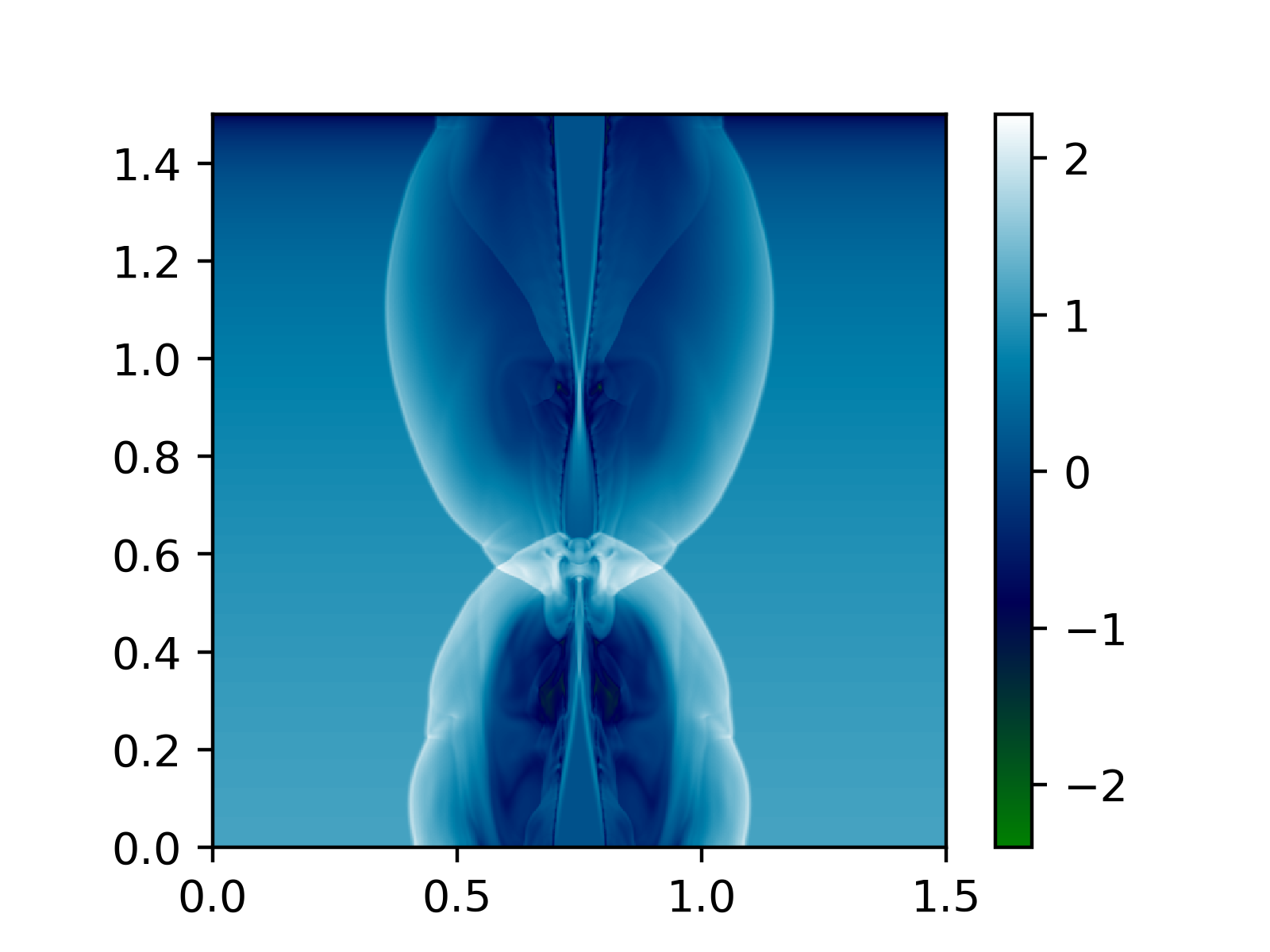}   
     \includegraphics[height=5.5cm, trim={0.9cm 0cm 2.35cm 0cm},clip]
    {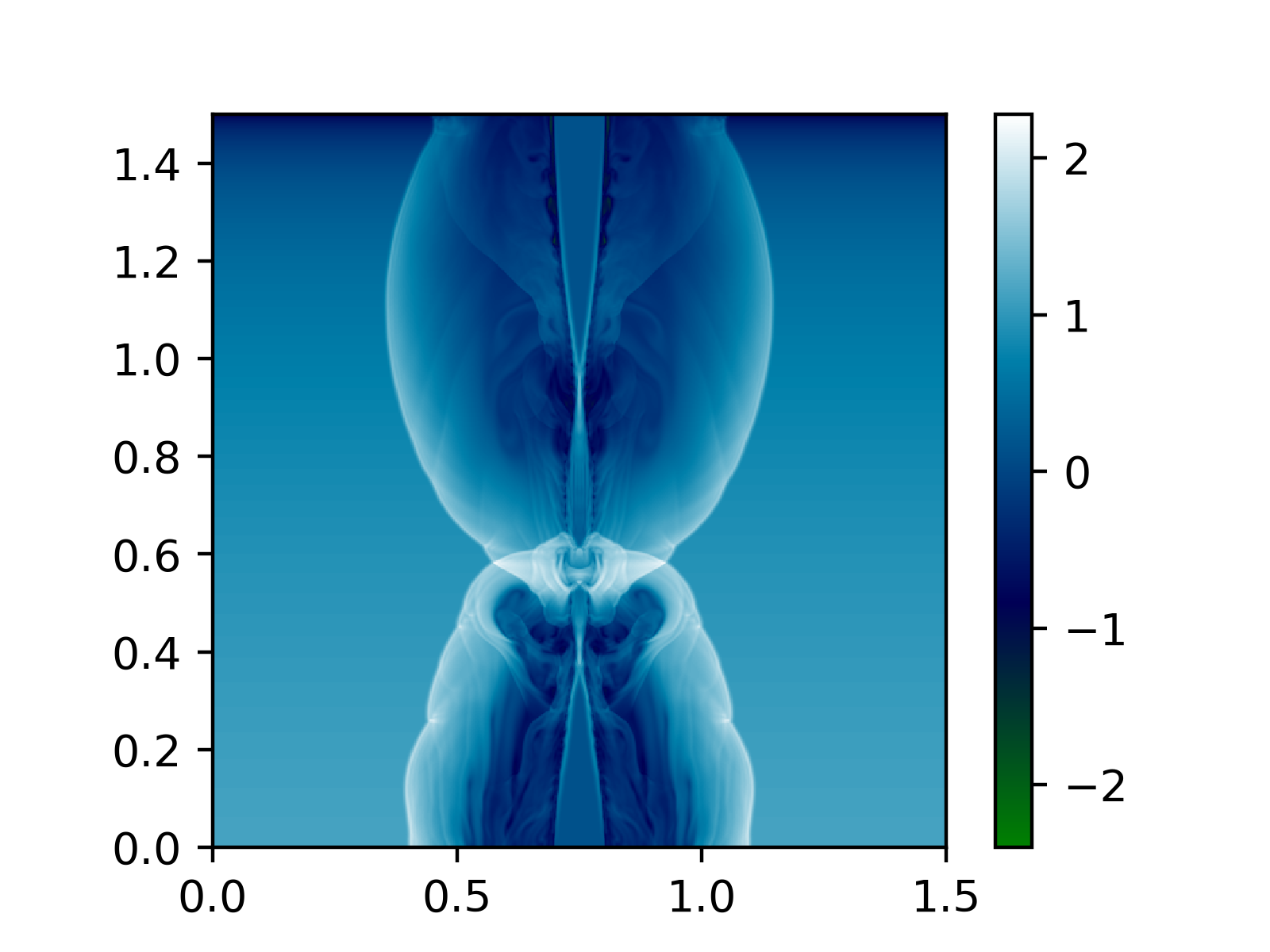} 
   \includegraphics[height=5.5cm, trim={0.9cm 0cm 1.15cm 0cm},clip]
    {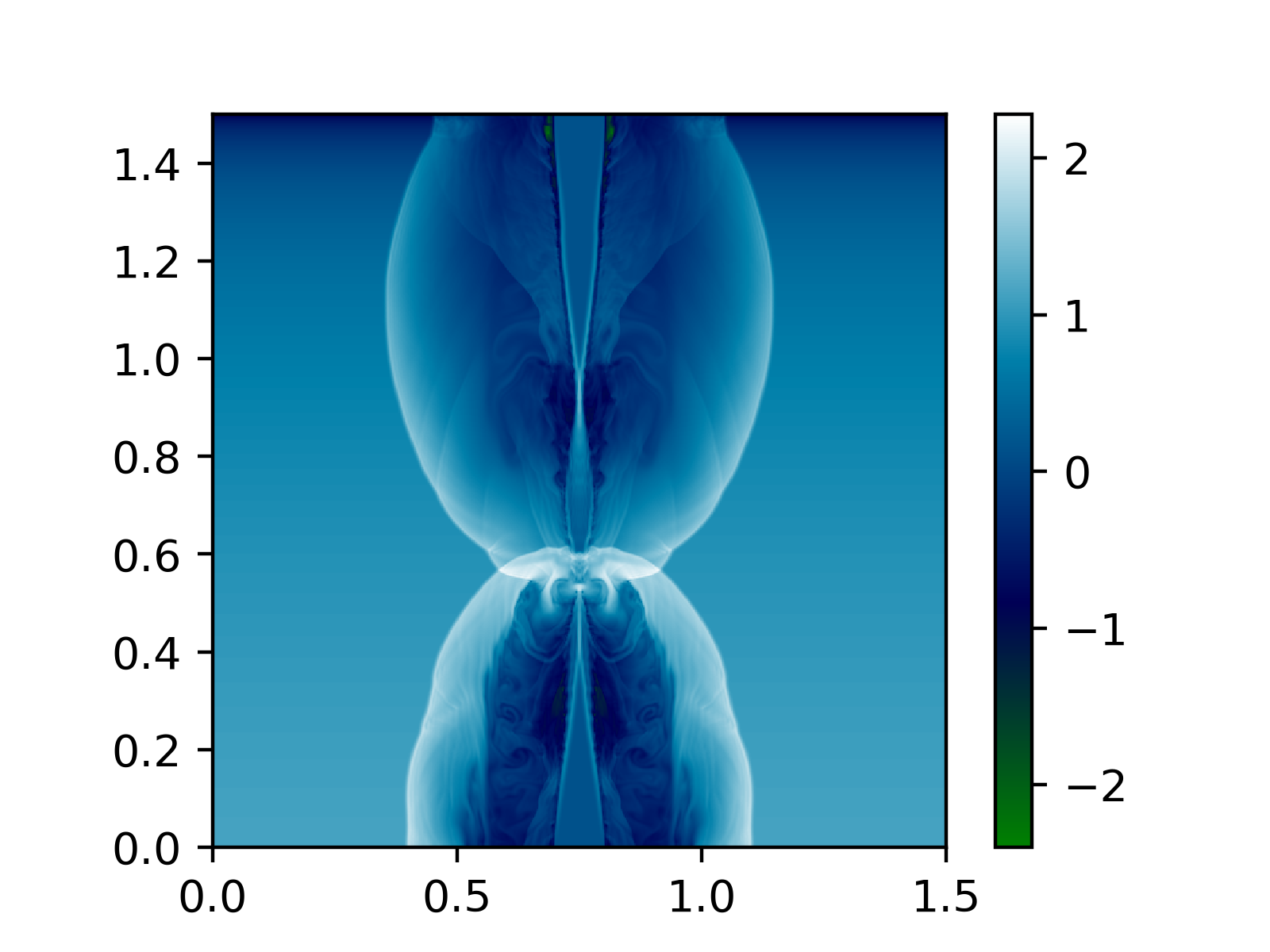}

    \includegraphics[height=5.5cm, trim={0.9cm 0cm 2.35cm 0cm},clip]
    {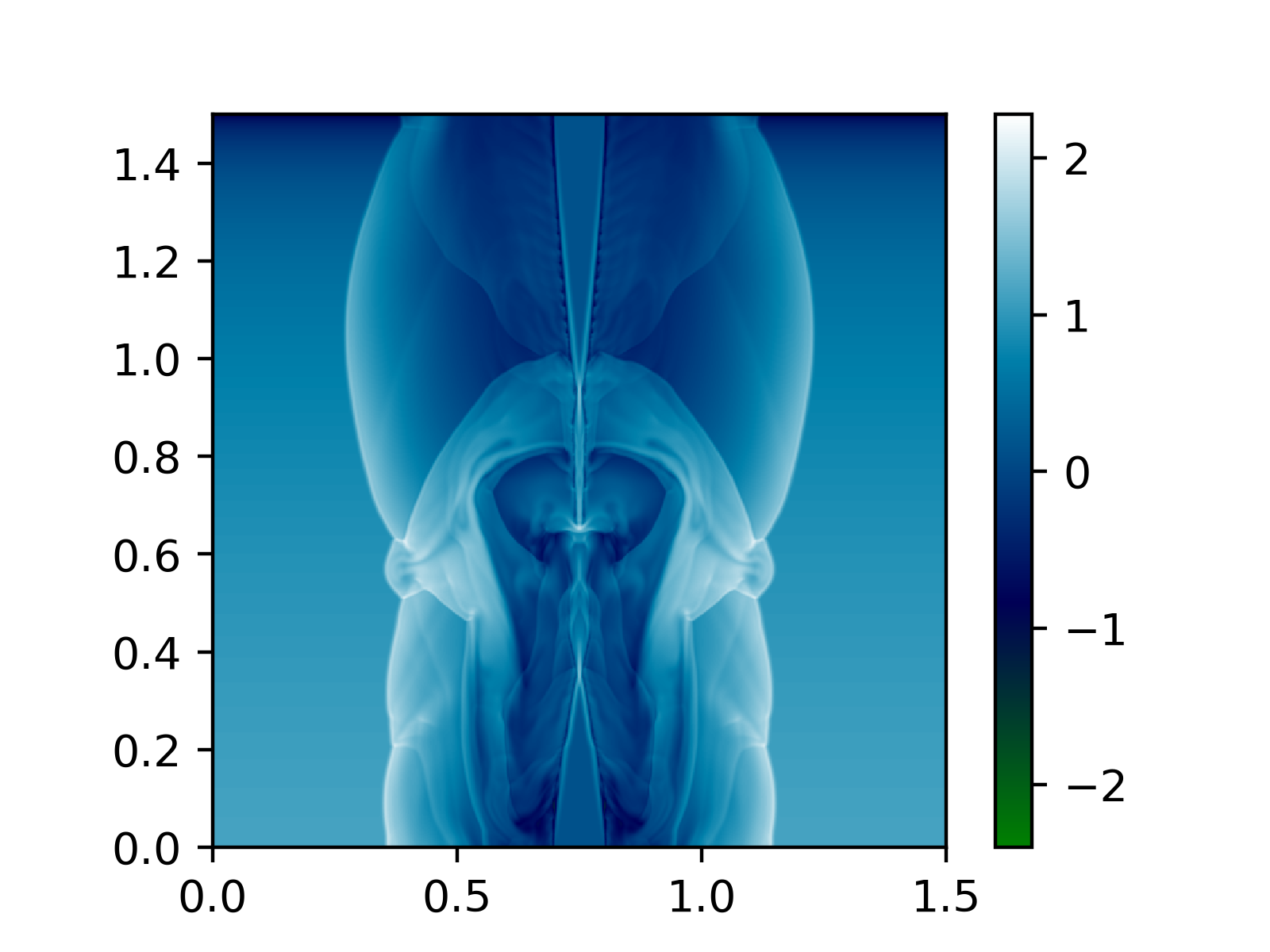}  
     \includegraphics[height=5.5cm, trim={0.9cm 0cm 2.35cm 0cm},clip]
    {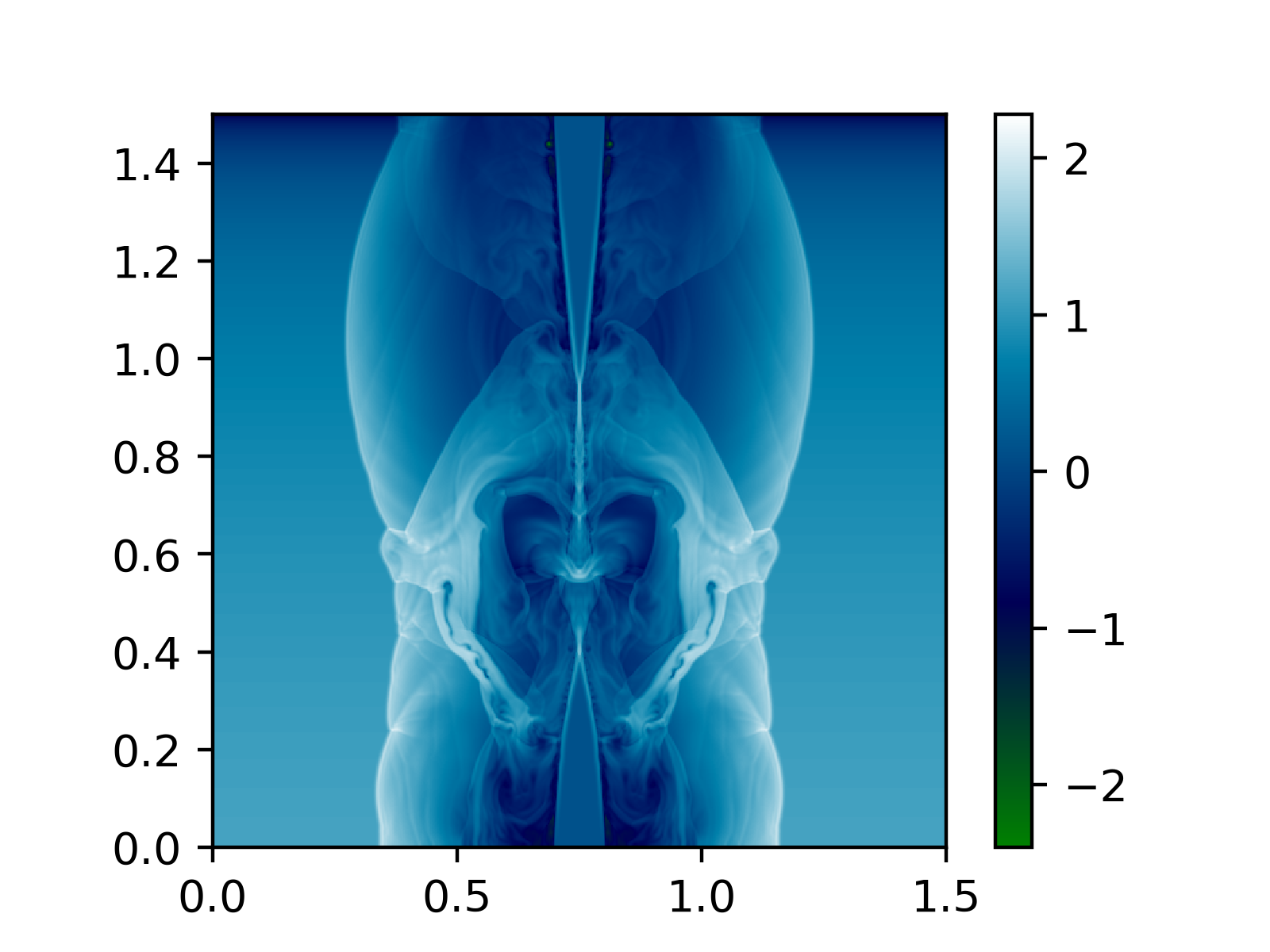} 
   \includegraphics[height=5.5cm, trim={0.9cm 0cm 1.15cm 0cm},clip]
    {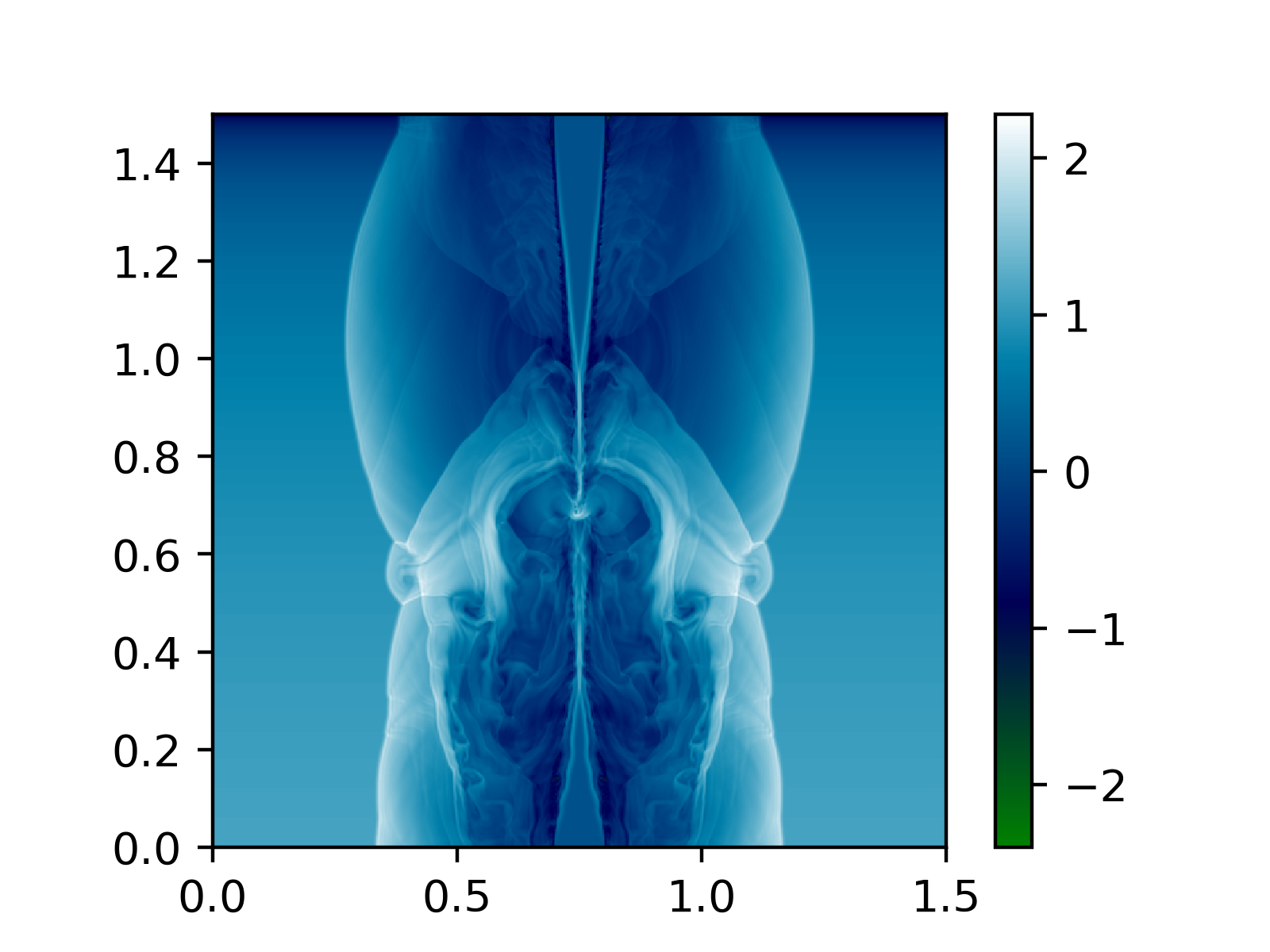}\\

    \includegraphics[height=5.5cm, trim={0.9cm 0cm 2.35cm 0cm},clip]
    {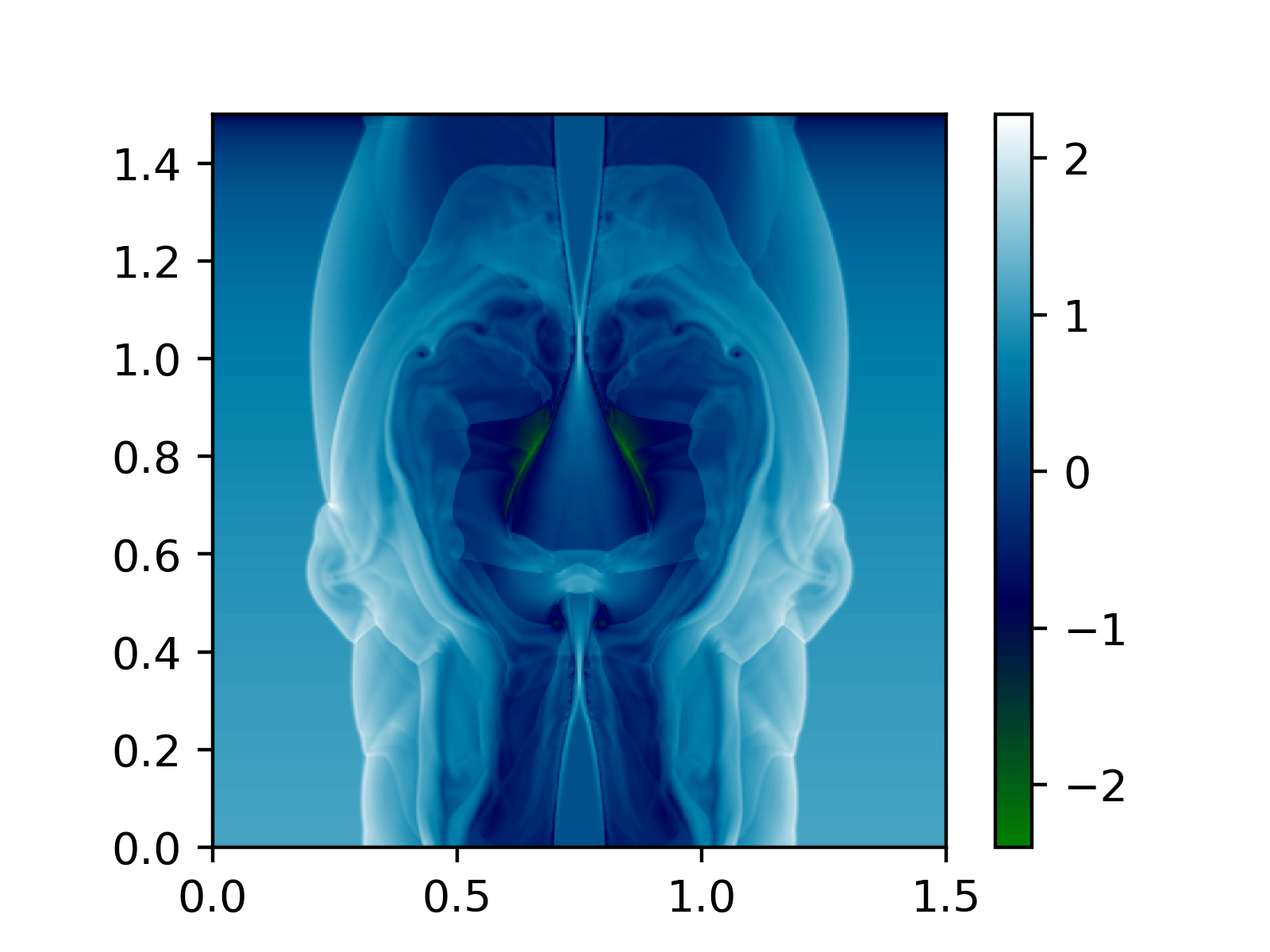}  
     \includegraphics[height=5.5cm, trim={0.9cm 0cm 2.35cm 0cm},clip]
    {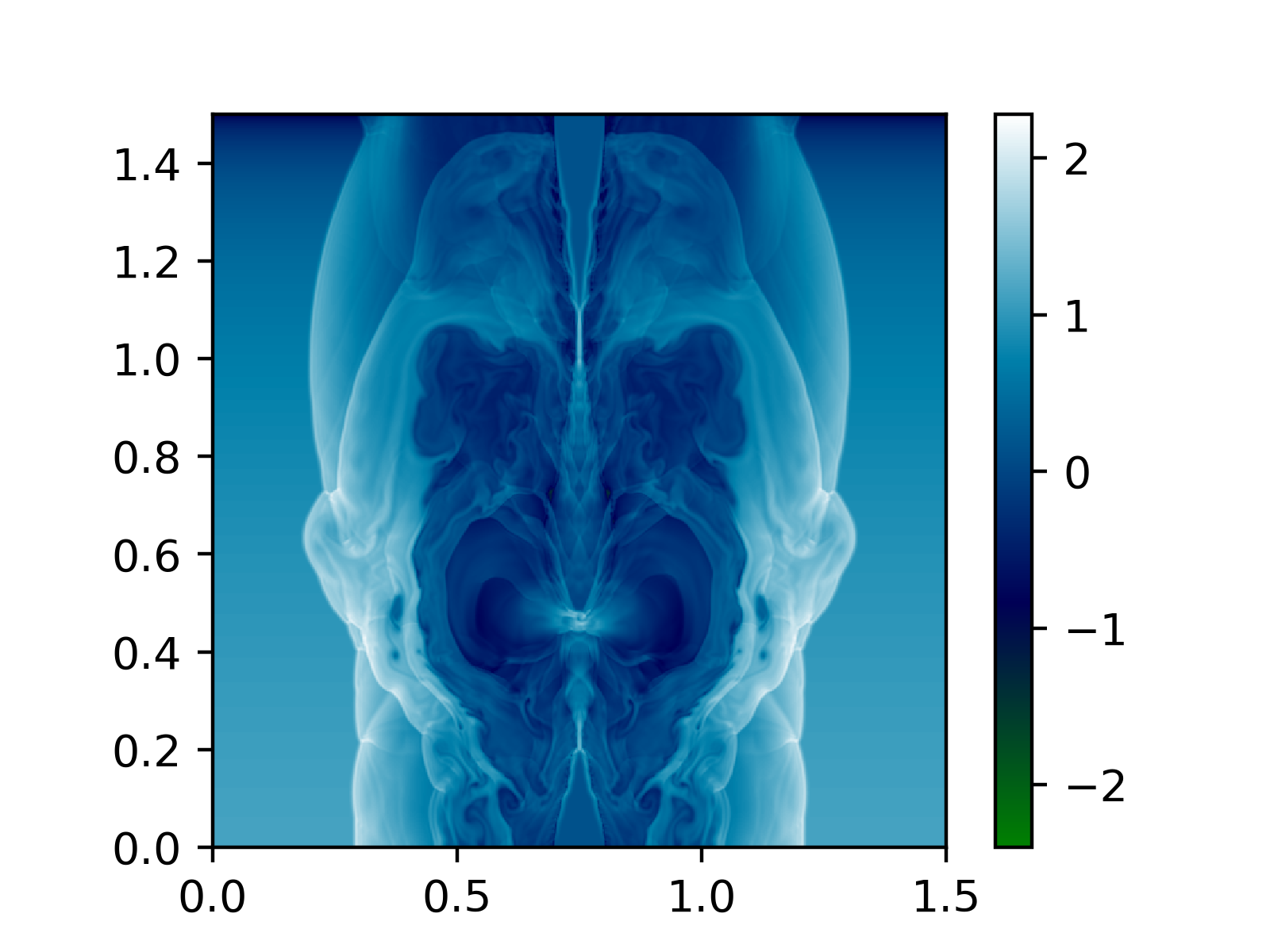} 
   \includegraphics[height=5.5cm, trim={0.9cm 0cm 1.15cm 0cm},clip]
    {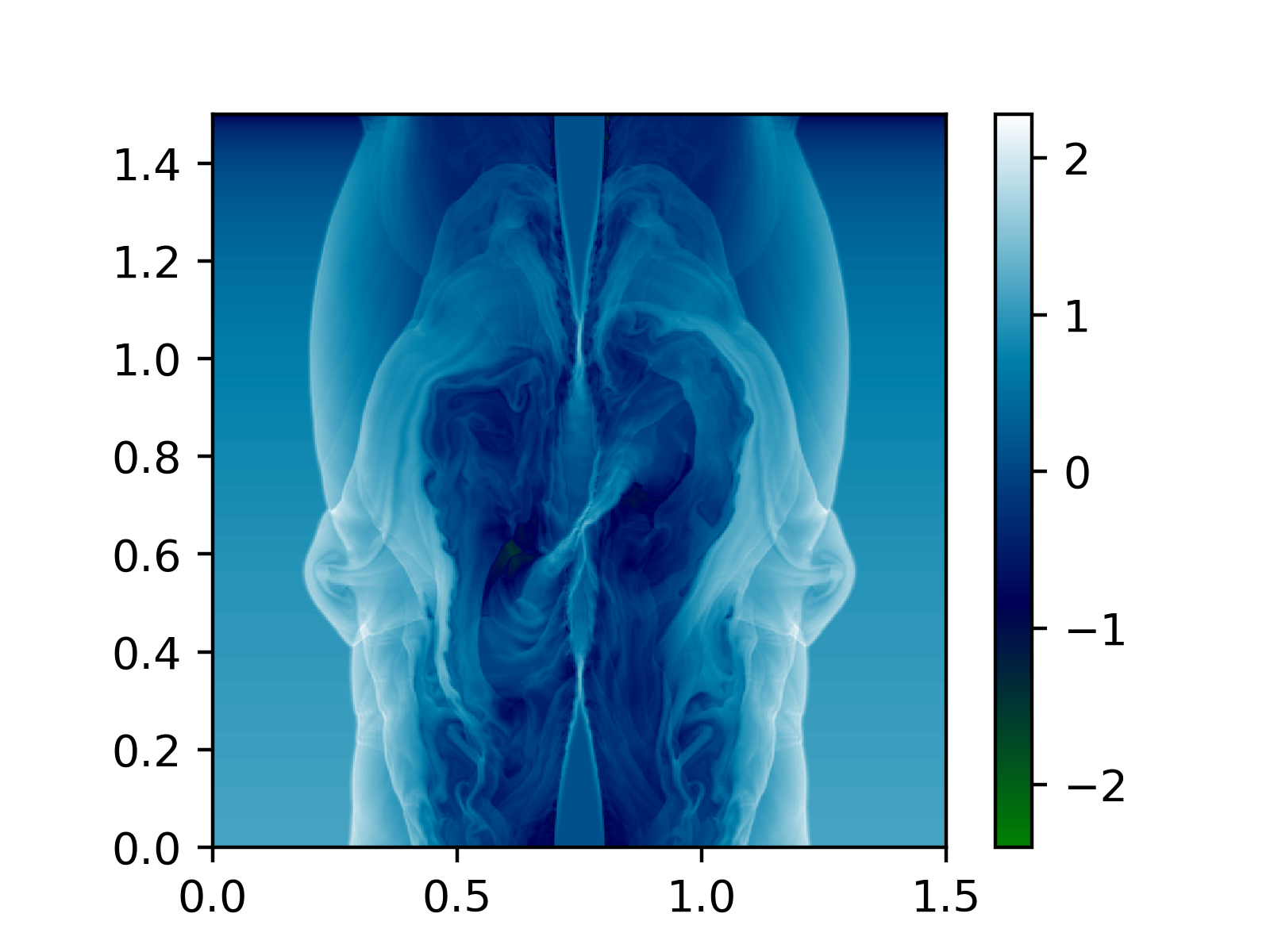} \\

 \caption{The double Mach 800 jet collision problem is resolved on a 
 $600 \times 600$ grid resolution. Density profiles are displayed.
From left to right, the results are computed using
GP-MOOD3, GP-MOOD5, and GP-MOOD7, while
from top to bottom, the densities are plotted at 
$t=0.002, 0.003, 0.004$, and $0.005$.
The runtime parameters are the same as the single jet problem.}
\label{fig:dbl_mach_collision}
\end{figure}

In the next problem, we extend the previous Mach jet configuration to a
double Mach jet collision that leads to a more violent jet-jet head-on collision.
The initial condition of the single Mach 100 light jet in \cref{sec:Mach_jets}
is modified to add two new configurations. 
First, the density of the ambient flow is stratified following
\beq\label{eq:stratified}
\rho_{\mbox{ambient}} = -9.24 y + 14, \;\;\; 0 \le y \le 1.5,
\eeq
which makes the density at the bottom to be 100 times heavier than the
top, i.e., 
$\rho_{\mbox{ambient}} (0) = 14.0$ and 
$\rho_{\mbox{ambient}} (1.5) = 0.14$.
Second, we modify the jet configuration by
introducing a second Mach jet injected from the top,
and at the same time, increasing the two jets' injection velocity by eight times, i.e.,
\beq\label{eq:double_jet}
(\rho, u, v, p)_{\mbox{jet}} = (\gamma, 0, 800, 1), \;\;\; 0.7 \le x \le 0.8,
\eeq
both at $y=0$ and $y=1.5$ with $\gamma = 1.4$.
Since the ambient gas density is ten times lighter than the jet at $y=1.5$,
the top jet becomes a dense Mach 800 jet, while the bottom
jet remains as a light Mach 800 jet similar to the previous single jet setup.

In \cref{fig:dbl_mach_collision}, from the left to right columns,
we display the density evolutions
of GP-MOOD3, GP-MOOD5, and GP-MOOD7.
From top to bottom, the density profiles are plotted at 
$t=0.002, 0.003, 0.004$, and $0.005$.
The top jet propagates through the ambient flow as  \textit{a dense jet}
up to $y\approx 1.3637$, beyond which the ambient density becomes
denser than the jet's density. The jet then becomes a light jet and
starts to undergo the sausage-like mode, as explained
in \cref{sec:Mach_jets}. This is seen in the density results
at $t=0.002$, where the top jet starts to get significantly 
squeezed into a narrower shape for $y \le 1.3637$.
The top jet moves faster 
than the bottom jet towards the collision point since
the top jet propagates through
the ambient density that is lighter in the top part of the domain.
The shape of the outer shock wave sheath still remains 
as a clean oval shape for the top dense jet, 
while the bottom light jet's sheath already exhibits 
a pair of kinks at $y \approx 0.22$.

At $t=0.003$, the two jets have already made a head-on collision.
The collision produces
highly turbulent fluid motions that are progressively amplified
as the two jets continue to make their ways in the opposite directions
for $t \ge 0.003$.
Driven by the turbulence, 
the low-dense gas in the valley pushes the high-dense gas outward,
which further excites the Rayleigh Taylor instabilities in
the already Kelvin-Helmholtz unstable gas in the cocoon.
It is observed that the symmetry is well preserved up to $t=0.003$ 
in all three solutions.
However, developments of asymmetric flows are seen
in the two higher order solutions later times.
For instance, with GP-MOOD7 at $t=0.004$,
a small region of asymmetrical flow is observed 
at $(x,y)\approx (0.75,0.71)$, where the two jets make a strong
collision and keep pushing each other.
With GP-MOOD5, a flow asymmetry is slightly delayed and seen 
at $t=0.005$ around $(x,y)\approx (0.75,0.45)$.
By $t=0.005$, the GP-MOOD7 solution has evolved into
a more asymmetrical flow in the central region.
We relate this asymmetry issue with the recent
study by Fleischmann \textit{et al.}~\cite{fleischmann2019numerical}.
The study demonstrates that symmetry-breaking phenomena,
particularly with shock interactions and turbulent structures,
are closely associated with vanishing numerical viscosity.
As such, high-order methods and highly resolved grid solutions
are more prone to asymmetrical flow developments.
The authors conclude that the main source of these 
asymmetric instabilities originates from floating-point truncation
errors that grow rapidly over time. Such errors result from
algorithmic artifacts (e.g., the lack of associativity)
and do not originate from the nature of high-order schemes.
Given that our GP-MOOD methods do not suffer
from asymmetric instabilities in other problems,
the flow asymmetries in our results with
GP-MOOD5 and GP-MOOD7 may well be addressed
with more rigorous code implementations
that strictly impose the proposed associativity
arithmetics as studied in \cite{fleischmann2019numerical}.
We leave a further study on this topic in our future work.

Lastly, we report that the number of cells that undergo
the order decrement in the MOOD loop remains less than
5.5\% of the entire cells in each step before the jet-jet collision, 
which grows as high as 6.7\% by $t=0.005$.
Although this test is the most stringent
problem in this paper in terms of the 
shock strength and the flow compressibility,
about 93\% of the domain is shown to be
stable enough with the solution
by the high-order unlimited linear GP reconstruction
methods. This finding can support the algorithmic advantage
of MOOD's \textit{a posterior} strategy over the conventional
\textit{a priori} shock-capturing strategy that requires 
actively employing computationally expensive nonlinear limiters all the time. 





\section{Conclusion}\label{sec:conclusion}
We have presented multidimensional, positivity-preserving high-order
(3rd-, 5th-, and 7th-order) finite volume 
GP methods integrated in the MOOD framework.
Extended from the existing polynomial-MOOD methods, the GP-MOOD methods 
replace the unlimited polynomial reconstruction schemes with a new
set of unlimited GP reconstruction algorithms that are genuinely
multidimensional. The GP-MOOD methods have shown that they
can deliver non-oscillatory, highly accurate, and positive-preserving
numerical solutions on a set of stringent benchmark problems,
including two highly compressible supersonic Mach jet test problems.

Compared to the existing 2D polynomial MOOD methods, the GP-MOOD
methods can deliver high-order solutions using smaller local stencil sizes,
e.g., the 5-points GP stencil for the 3rd-order GP-MOOD3, 
the 13-points GP stencil for the 5th-order GP-MOOD5, and
the 25-points GP stencil for the 7th-order GP-MOOD7.
The GP-MOOD methods balance the two demanding numerical requirements, 
solution accuracy and stability, 
via the series of detection algorithms comprising the MOOD loop.
The solution accuracy is decreased at troubled cells
from the highest accurate 
choice of an unlimited GP reconstruction algorithm down to the
least accurate first-order Godunov method, by which 
GP-MOOD gains increasing numerical stability.
This study also introduced a new detection switch called the
``Compressibility-Shock Detection'' (or CSD) that improves
the excessive diffusivity in the existing polynomial MOOD methods.
We have demonstrated that CSD is essential to capture some
crucial flow dynamics that are signature to 
certain benchmark problems.

We demonstrate that the percentage of those troubled cells,
where GP-MOOD's order decrements takes place,
is only a fraction of the entire domain and never exceed
10\%. Our finding is consistent with the former
studies on the polynomial MOOD methods
\cite{diot2013multidimensional,diot2012methode,bourriaud2020priori}.

Our work in this paper designs our GP-MOOD
methods for the hyperbolic systems of the Euler equations in 1D and 2D.
However, our study is not limited to this model; 
the GP-MOOD methods can seamlessly be extended to 
the 3D version of the Euler equations as well as other physics
system such as magnetohydrodynamics, relativistic (magneto)hydrodynamics,
to name a few. 
Lastly, we remark that the accuracy of GP-MOOD is not limited to
7th-order but can be further extended to a higher accuracy of $(2R+1)$ 
with $R>3$. These topics will be investigated in our future work.

\section*{Acknowledgements}\label{sec:ack}
The first author thanks the University of California Exchange Abroad program (UCEAP), 
the Centre de Californie de l'Université de Bordeaux and the international relations 
department of Enseirb-MATMECA for allowing him the opportunity to visit 
UC Santa Cruz and carry out the present work. 
Both authors would also like to thank Emile Josso for his insights on implementation, 
Dr. Youngjun Lee for his help on plotting simulation results, 
Oliver Speed and Ian May for the numerous fruitful discussions.
This work was supported in part by the National Science Foundation 
under grant AST-1908834. We also acknowledge the use of 
the Lux supercomputer at UC Santa Cruz, funded by NSF MRI grant AST-1828315.
\appendix

\section{Tabulated example of $\Delta_{mn,d}$ values}
\label{apdx:delta_kh}

\begin{figure}[h!]
    \centering
    \includegraphics[scale=0.4]{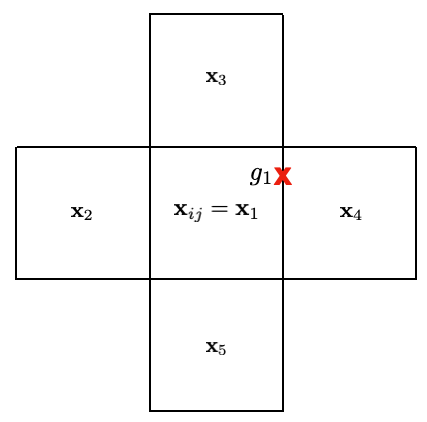}
    \caption{The 5-point GP stencil with a GP radius $R=1$ for the 3rd-order GP-MOOD3 method.
    One of the two quadrature points, $g_1$, is shown with the red cross at the right $x$-directional
    cell interface of the central cell, $\mathbf{x}_1 = \mathbf{x}_{ij}$.}
    \label{fig:example_delta}
\end{figure}

We give the tabulated values of  $\Delta_{mn,d}$ in Eq.~\eqref{eq:delta_mnd}
for the computation of the prediction 
vector at the Gaussian quadrature point, $(x,y)=(x_{i+1/2},g_1)$,
depicted in \cref{fig:example_delta}. See also \cref{fig:multipt-QR}.
We give an example for the 3rd-order GP-MOOD3 with a GP radius $R=1$,
for which the default quadrature rule is the 2-point, 4th-order quadrature rule,
\beq
\int_{y=-\frac{\Delta y}{2}}^{y=\frac{\Delta y}{2}} f(x,y) 
\simeq w_1 f(x,g_1) +  w_2 f(x,g_2),
\eeq
where 
$g_1 = \frac{\Delta y}{2\sqrt{3}}, 
g_2 = -\frac{\Delta y}{2\sqrt{3}},$ and
$w_1 = w_2 =  \frac{1}{2}$.
The values of $\Delta_{mn,d}$ for each $d=x, y$ are tabulated in 
\cref{tab:delta_x} and \cref{tab:delta_y}, respectively. 
The indices $m$ and $n$ represent
the cell index in the 5-point GP stencil in \cref{fig:example_delta}, i.e.,
$m,n = 1,  \dots, 5$ as well as the quadrature point, $g_1$.

\begin{table}[h!]
    \footnotesize
    \centering
    \caption{Tabulated values of $\Delta_{mn,x}$ at each pair of $(\mathbf{x}_m, \mathbf{x}_n)$ points.} 
    \begin{tabular}{@{}ccccccccccccccc@{}}
        \toprule
        {$\Delta_{mn,d}$}  & $m=1$ & $m=2$   & $m=3$  & $m=4$  &  $m=5$ & $m=g_1$     \\                                                                  
        \midrule
        $n=1$                    & $0$    & $1$      &  $0$    & $-1$    &   $0$    & $-1/2$ \\
        $n=2$                    & $-1$   & $0$      &  $-1$   & $-2$    &   $-1$   & $-3/2$ \\
        $n=3$                    & $0$    & $1$      &  $0$    & $-1$     &   $0$    & $-1/2$ \\
        $n=4$                    & $1$    & $2$      &  $1$    & $0$      &   $1$    & $1/2$ \\
        $n=5$                    & $0$    & $1$      & $0$     & $-1$     &   $0$    & $-1/2$ \\
        $n=g_1$                & $1/2$ & $3/2$   & $1/2$  & $-1/2$   &  $1/2$ & $0$ 
    \end{tabular}\label{tab:delta_x}
\end{table}

%
 
 \begin{table}[h!]
    \footnotesize
    \centering
    \caption{Tabulated values of $\Delta_{mn,y}$ at each pair of $(\mathbf{x}_m, \mathbf{x}_n)$ points.} 
    \begin{tabular}{@{}ccccccccccccccc@{}}
        \toprule
        {$\Delta_{mn,d}$}  & $m=1$ & $m=2$   & $m=3$  & $m=4$  &  $m=5$ & $m=g_1$     \\                                                                  
        \midrule
        $n=1$                    & $0$    & $0$      &  $-1$    & $0 $    &   $1$    & $-\frac{1}{2\sqrt{3}}$ \\
        $n=2$                    & $0$   & $0$      &  $-1$   & $0$   &   $1$   & $-\frac{1}{2\sqrt{3}}$ \\
        $n=3$                    & $1$    & $1$      &  $0$    & $1$     &   $2$    & $\frac{1}{2}\left(2-\frac{1}{\sqrt{3}}\right) $ \\
      $n=4$                    & $0$   & $0$      &  $-1$   & $0$   &   $1$   & $-\frac{1}{2\sqrt{3}}$ \\
        $n=5$                    & $-1$    & $-1$      & $-2$     & $-1$     &   $0$    & $-\frac{1}{2}\left(2+\frac{1}{\sqrt{3}}\right)$ \\
        $n=g_1$                & $\frac{1}{2\sqrt{3}}$ & $\frac{1}{2\sqrt{3}}$   
                                      & $-\frac{1}{2}\left(2-\frac{1}{\sqrt{3}}\right) $  & $\frac{1}{2\sqrt{3}}$   
                                      &  $\frac{1}{2}\left(2+\frac{1}{\sqrt{3}}\right)$ & $0$ 
    \end{tabular}\label{tab:delta_y}
\end{table}

%

\section{GP second derivatives formulas}\label{apdx:gp_2ndDer}
A GP-based second derivative formula can be derived to 
compute the second derivatives in \cref{eq:u2_curvature}.
%
This can be done by differentiating $\mathbf{t}_*$ twice in 
\cref{eq:gp_reconstructor} in each direction, $d=x, y$, since
only the prediction vector $\mathbf{t}_*$ is affected by the derivative,
\begin{equation}
\frac{\partial^2\tilde{m}_*} 
{\partial (\mathbf{e}_d \cdot \mathbf{x}_{*})^2 } 
= 
\frac{\partial^2 \mathbf{t}_*^{T}}
{\partial (\mathbf{e}_d \cdot \mathbf{x}_{*})^2 } 
\mathbf{C}^{-1}\mathbf{q}.
\label{eq:der}   
\end{equation}

The final expression of the GP-based second derivative is written as
\bea
    \frac{\partial^2 \mathbf{t}_{*,m}^{T}}
    {\partial (\mathbf{e}_d \cdot \mathbf{x}_{*})^2 } 
    &= &
    \prod_{\delta \neq d} \sqrt{\frac{\pi}{2}}\frac{\ell}{\Delta_{\delta}} \left \{
     \erf{\frac{\Delta_{m*,{\delta}}+1/2}{\sqrt{2}\ell/\Delta_{\delta}}}
     - \erf{\frac{\Delta_{m*,{\delta}}-1/2}{\sqrt{2}\ell/\Delta_{\delta}}}
     \right \} \\
     &\times& \frac{\partial^2}
     {\partial (\mathbf{e}_d \cdot \mathbf{x}_{*})^2 } 
     \left[\sqrt{\frac{\pi}{2}}\frac{\ell}{\Delta_d} \left \{
     \erf{\frac{\Delta_{m*,d}+1/2}{\sqrt{2}\ell/\Delta_d}}
     - \erf{\frac{\Delta_{m*,d}-1/2}{\sqrt{2}\ell/\Delta_d}}
     \right \}\right]\\
      &=& \prod_{\delta \neq d} \sqrt{\frac{\pi}{2}}\frac{\ell}{\Delta_{\delta}} \left \{
     \erf{\frac{\Delta_{m*,{\delta}}+1/2}{\sqrt{2}\ell/\Delta_{\delta}}}
     - \erf{\frac{\Delta_{m*,{\delta}}-1/2}{\sqrt{2}\ell/\Delta_{\delta}}}
     \right \} \\
     &\times& \frac{1}{\ell^2}
     \left[ 
     \left(\Delta_{m*,d} - \frac{1}{2} \right) \expo{ - \frac {(\Delta_{m*,d} - \frac{1}{2})^2}{2(\ell/\Delta_d)^2}} - 
     \left(\Delta_{m*,d} + \frac{1}{2} \right) \expo{ - \frac {(\Delta_{m*,d} + \frac{1}{2})^2}{2(\ell/\Delta_d)^2}}
     \right].
\eea
%
Here, we choose $\mathbf{x}_{*}$ as the center of the cell to obtain the cell-centered second derivative.
%
%

\section{The full set of coefficients of the 3rd-order polynomial, $p_3(x,y)$}
\label{apdx:p3_coeffs}
The five coefficients, $a_0, \dots, a_4$, of $p_3(x,y)$ are obtained by solving
\cref{eq:p3_coeffs}, written in terms of $\overline{\bq}_{i_m j_m}$, $1 \le i_mj_m \le 5$, 
as follows:
\bea
a_0 &=&\frac{7}{6}\overline{\bq}_{1}
            -\frac{1}{24}\overline{\bq}_{2}
            -\frac{1}{24}\overline{\bq}_{3}
            -\frac{1}{24}\overline{\bq}_{4}
            -\frac{1}{24}\overline{\bq}_{5},\\
a_1 &=&
            -\frac{1}{2}\overline{\bq}_{2}
            +\frac{1}{2}\overline{\bq}_{4},\\
a_2 &=&-\overline{\bq}_{1}
            +\frac{1}{2}\overline{\bq}_{2}
            +\frac{1}{2}\overline{\bq}_{4},\\
a_3 &=&
             \frac{1}{2}\overline{\bq}_{3}
            -\frac{1}{2}\overline{\bq}_{5},\\     
a_4 &=&-\overline{\bq}_{1}
            +\frac{1}{2}\overline{\bq}_{3}
            +\frac{1}{2}\overline{\bq}_{5}.
\eea

\bibliography{refs_merged_new}

\begin{thebibliography}{10}

\bibitem{van1974towards}
Bram Van~Leer.
\newblock Towards the ultimate conservative difference scheme. ii. monotonicity
  and conservation combined in a second-order scheme.
\newblock {\em Journal of computational physics}, 14(4):361--370, 1974.

\bibitem{tadmor1988convenient}
Eitan Tadmor.
\newblock Convenient total variation diminishing conditions for nonlinear
  difference schemes.
\newblock {\em SIAM journal on numerical analysis}, 25(5):1002--1014, 1988.

\bibitem{hubbard1999multidimensional}
ME~Hubbard.
\newblock Multidimensional slope limiters for muscl-type finite volume schemes
  on unstructured grids.
\newblock {\em Journal of Computational Physics}, 155(1):54--74, 1999.

\bibitem{harten1997high}
Ami Harten.
\newblock High resolution schemes for hyperbolic conservation laws.
\newblock {\em Journal of computational physics}, 135(2):260--278, 1997.

\bibitem{colella1984piecewise}
Phillip Colella and Paul~R Woodward.
\newblock The piecewise parabolic method ({PPM}) for gas-dynamical simulations.
\newblock {\em Journal of computational physics}, 54(1):174--201, 1984.

\bibitem{mccorquodale2011high}
Peter McCorquodale and Phillip Colella.
\newblock A high-order finite-volume method for conservation laws on locally
  refined grids.
\newblock {\em Communications in Applied Mathematics and Computational
  Science}, 6(1):1--25, 2011.

\bibitem{harten1997uniformly}
Ami Harten and Stanley Osher.
\newblock Uniformly high-order accurate nonoscillatory schemes. i.
\newblock In {\em Upwind and High-Resolution Schemes}, pages 187--217.
  Springer, 1997.

\bibitem{shu1988efficient}
Chi-Wang Shu and Stanley Osher.
\newblock Efficient implementation of essentially non-oscillatory
  shock-capturing schemes.
\newblock {\em Journal of computational physics}, 77(2):439--471, 1988.

\bibitem{liu1994weighted}
Xu-Dong Liu, Stanley Osher, Tony Chan, et~al.
\newblock Weighted essentially non-oscillatory schemes.
\newblock {\em Journal of computational physics}, 115(1):200--212, 1994.

\bibitem{jiang1996efficient}
Guang-Shan Jiang and Chi-Wang Shu.
\newblock Efficient implementation of weighted {ENO} schemes.
\newblock {\em Journal of computational physics}, 126(1):202--228, 1996.

\bibitem{balsara2000monotonicity}
Dinshaw~S Balsara and Chi-Wang Shu.
\newblock Monotonicity preserving weighted essentially non-oscillatory schemes
  with increasingly high order of accuracy.
\newblock {\em Journal of Computational Physics}, 160(2):405--452, 2000.

\bibitem{gerolymos2009very}
GA~Gerolymos, D~S{\'e}n{\'e}chal, and I~Vallet.
\newblock Very-high-order weno schemes.
\newblock {\em Journal of Computational Physics}, 228(23):8481--8524, 2009.

\bibitem{levy2000compact}
Doron Levy, Gabriella Puppo, and Giovanni Russo.
\newblock Compact central weno schemes for multidimensional conservation laws.
\newblock {\em SIAM Journal on Scientific Computing}, 22(2):656--672, 2000.

\bibitem{ivan2014high}
Lucian Ivan and Clinton~PT Groth.
\newblock High-order solution-adaptive central essentially non-oscillatory
  (ceno) method for viscous flows.
\newblock {\em Journal of Computational Physics}, 257:830--862, 2014.

\bibitem{semplice2016adaptive}
Matteo Semplice, Armando Coco, and Giovanni Russo.
\newblock Adaptive mesh refinement for hyperbolic systems based on third-order
  compact weno reconstruction.
\newblock {\em Journal of Scientific Computing}, 66(2):692--724, 2016.

\bibitem{dumbser2017central}
Michael Dumbser, Walter Boscheri, Matteo Semplice, and Giovanni Russo.
\newblock Central weighted eno schemes for hyperbolic conservation laws on
  fixed and moving unstructured meshes.
\newblock {\em SIAM Journal on Scientific Computing}, 39(6):A2564--A2591, 2017.

\bibitem{qiu2004hermite}
Jianxian Qiu and Chi-Wang Shu.
\newblock Hermite weno schemes and their application as limiters for
  runge--kutta discontinuous galerkin method: one-dimensional case.
\newblock {\em Journal of Computational Physics}, 193(1):115--135, 2004.

\bibitem{balsara2007sub}
Dinshaw~S Balsara, Christoph Altmann, Claus-Dieter Munz, and Michael Dumbser.
\newblock A sub-cell based indicator for troubled zones in rkdg schemes and a
  novel class of hybrid rkdg+ hweno schemes.
\newblock {\em Journal of Computational Physics}, 226(1):586--620, 2007.

\bibitem{balsara2016efficient}
Dinshaw~S Balsara, Sudip Garain, and Chi-Wang Shu.
\newblock An efficient class of weno schemes with adaptive order.
\newblock {\em Journal of Computational Physics}, 326:780--804, 2016.

\bibitem{reyes2018new}
Adam Reyes, Dongwook Lee, Carlo Graziani, and Petros Tzeferacos.
\newblock A new class of high-order methods for fluid dynamics simulations
  using gaussian process modeling: One-dimensional case.
\newblock {\em Journal of Scientific Computing}, 76(1):443--480, 2018.

\bibitem{reyes2019variable}
Adam Reyes, Dongwook Lee, Carlo Graziani, and Petros Tzeferacos.
\newblock A variable high-order shock-capturing finite difference method with
  {GP-WENO}.
\newblock {\em Journal of Computational Physics}, 381:189--217, 2019.

\bibitem{reeves2020application}
Steve Reeves, Dongwook Lee, Adam Reyes, Carlo Graziani, and Petros Tzeferacos.
\newblock An application of gaussian process modeling for high-order accurate
  adaptive mesh refinement prolongation.
\newblock {\em arXiv preprint arXiv:2003.08508; Accepted for publication in
  CAMCoS}, 2021.

\bibitem{kent2014determining_part2}
James Kent, Christiane Jablonowski, Jared~P Whitehead, and Richard~B Rood.
\newblock Determining the effective resolution of advection schemes. part {II}:
  Numerical testing.
\newblock {\em Journal of Computational Physics}, 278:497--508, 2014.

\bibitem{clain2011high}
St{\'e}phane Clain, Steven Diot, and Rapha{\"e}l Loub{\`e}re.
\newblock A high-order finite volume method for systems of conservation
  laws---multi-dimensional optimal order detection ({MOOD}).
\newblock {\em Journal of computational Physics}, 230(10):4028--4050, 2011.

\bibitem{diot2012improved}
Steven Diot, St{\'e}phane Clain, and Rapha{\"e}l Loub{\`e}re.
\newblock Improved detection criteria for the multi-dimensional optimal order
  detection ({MOOD}) on unstructured meshes with very high-order polynomials.
\newblock {\em Computers \& Fluids}, 64:43--63, 2012.

\bibitem{diot2013multidimensional}
Steven Diot, Rapha{\"e}l Loub{\`e}re, and Stephane Clain.
\newblock The multidimensional optimal order detection method in the
  three-dimensional case: very high-order finite volume method for hyperbolic
  systems.
\newblock {\em International Journal for Numerical Methods in Fluids},
  73(4):362--392, 2013.

\bibitem{diot2012methode}
Steven Diot.
\newblock {\em La m{\'e}thode MOOD Multi-dimensional Optimal Order Detection:
  la premi{\`e}re approche a posteriori aux m{\'e}thodes volumes finis d'ordre
  tr{\`e}s {\'e}lev{\'e}}.
\newblock PhD thesis, Universit{\'e} de Toulouse, Universit{\'e} Toulouse
  III-Paul Sabatier, 2012.

\bibitem{semplice2018adaptive}
Matteo Semplice and Rapha{\"e}l Loub{\`e}re.
\newblock Adaptive-mesh-refinement for hyperbolic systems of conservation laws
  based on a posteriori stabilized high order polynomial reconstructions.
\newblock {\em Journal of Computational Physics}, 354:86--110, 2018.

\bibitem{toro2001towards}
EF~Toro, RC~Millington, and LAM Nejad.
\newblock Towards very high order {Godunov} schemes.
\newblock In {\em Godunov methods}, pages 907--940. Springer, 2001.

\bibitem{titarev2002ader}
Vladimir~A Titarev and Eleuterio~F Toro.
\newblock {ADER}: Arbitrary high order {Godunov} approach.
\newblock {\em Journal of Scientific Computing}, 17(1-4):609--618, 2002.

\bibitem{titarev2005ader}
Vladimir~A Titarev and Eleuterio~F Toro.
\newblock {ADER} schemes for three-dimensional non-linear hyperbolic systems.
\newblock {\em Journal of Computational Physics}, 204(2):715--736, 2005.

\bibitem{dumbser2014posteriori}
Michael Dumbser, Olindo Zanotti, Rapha{\"e}l Loub{\`e}re, and Steven Diot.
\newblock A posteriori subcell limiting of the discontinuous galerkin finite
  element method for hyperbolic conservation laws.
\newblock {\em Journal of Computational Physics}, 278:47--75, 2014.

\bibitem{dumbser2016simple}
Michael Dumbser and Rapha{\"e}l Loub{\`e}re.
\newblock A simple robust and accurate a posteriori sub-cell finite volume
  limiter for the discontinuous galerkin method on unstructured meshes.
\newblock {\em Journal of Computational Physics}, 319:163--199, 2016.

\bibitem{bourriaud2020priori}
Alexandre Bourriaud, Rapha{\"e}l Loub{\`e}re, and Rodolphe Turpault.
\newblock A priori neural networks versus a posteriori mood loop: A high
  accurate 1d fv scheme testing bed.
\newblock {\em Journal of Scientific Computing}, 84(2):1--36, 2020.

\bibitem{shi_technique_2002}
Jing Shi, Changqing Hu, and Chi-Wang Shu.
\newblock A {Technique} of {Treating} {Negative} {Weights} in {WENO} {Schemes}.
\newblock {\em Journal of Computational Physics}, 175(1):108--127, January
  2002.

\bibitem{shu_high-order_2003}
Chi-Wang Shu.
\newblock High-order finite difference and finite volume {WENO} schemes and
  discontinuous {{G}alerkin} methods for {CFD}.
\newblock {\em International Journal of Computational Fluid Dynamics},
  17(2):107--118, 2003.

\bibitem{balsara2009divergence}
Dinshaw~S Balsara.
\newblock Divergence-free reconstruction of magnetic fields and weno schemes
  for magnetohydrodynamics.
\newblock {\em Journal of Computational Physics}, 228(14):5040--5056, 2009.

\bibitem{mccorquodale2015adaptive}
P.~McCorquodale, P.~Ullrich, H.~Johansen, and P.~Colella.
\newblock An adaptive multiblock high-order finite-volume method for solving
  the shallow-water equations on the sphere.
\newblock {\em Communications in Applied Mathematics and Computational
  Science}, 10(2):121--162, 2015.

\bibitem{zhang2012fourth}
Q.~Zhang, H.~Johansen, and P.~Colella.
\newblock A fourth-order accurate finite-volume method with structured adaptive
  mesh refinement for solving the advection-diffusion equation.
\newblock {\em SIAM Journal on Scientific Computing}, 34(2):B179--B201, 2012.

\bibitem{amrex}
W.~Zhang, A.~Almgren, V.~Beckner, J.~Bell, J.~Blaschke, C.~Chan, M.~Day,
  B.~Friesen, K.~Gott, D.~Graves, M.~Katz, A.~Myers, T.~Nguyen, A.~Nonaka,
  M.~Rosso, S.~Williams, and M.~Zingale.
\newblock Amrex: a framework for block-structured adaptive mesh refinement.
\newblock {\em Journal of Open Source Software}, 4(37):1370, 5 2019.

\bibitem{gottlieb1998total}
Sigal Gottlieb and Chi-Wang Shu.
\newblock Total variation diminishing {Runge-Kutta} schemes.
\newblock {\em Mathematics of computation of the American Mathematical
  Society}, 67(221):73--85, 1998.

\bibitem{spiteri2002new}
Raymond~J Spiteri and Steven~J Ruuth.
\newblock A new class of optimal high-order strong-stability-preserving time
  discretization methods.
\newblock {\em SIAM Journal on Numerical Analysis}, 40(2):469--491, 2002.

\bibitem{lee2017piecewise}
Dongwook Lee, Hugues Faller, and Adam Reyes.
\newblock The piecewise cubic method ({PCM}) for computational fluid dynamics.
\newblock {\em Journal of Computational Physics}, 341:230--257, 2017.

\bibitem{rasmussen2005}
C.E. Rasmussen and C.K.I. Williams.
\newblock {\em Gaussian {Processes} for Machine Learning}.
\newblock Adaptive Computation And Machine Learning. MIT Press, 2005.

\bibitem{bishop2007pattern}
C.~Bishop.
\newblock Pattern recognition and machine learning (information science and
  statistics), 1st edn. 2006. corr. 2nd printing edn.
\newblock {\em Springer, New York}, 2007.

\bibitem{fornberg2004stable}
Bengt Fornberg and Grady Wright.
\newblock Stable computation of multiquadric interpolants for all values of the
  shape parameter.
\newblock {\em Computers \& Mathematics with Applications}, 48(5-6):853--867,
  2004.

\bibitem{fornberg2011stable}
Bengt Fornberg, Elisabeth Larsson, and Natasha Flyer.
\newblock Stable computations with gaussian radial basis functions.
\newblock {\em SIAM Journal on Scientific Computing}, 33(2):869--892, 2011.

\bibitem{fornberg2008stable}
Bengt Fornberg and C{\'e}cile Piret.
\newblock A stable algorithm for flat radial basis functions on a sphere.
\newblock {\em SIAM Journal on Scientific Computing}, 30(1):60--80, 2008.

\bibitem{fasshauer2012stable}
Gregory~E Fasshauer and Michael~J McCourt.
\newblock Stable evaluation of gaussian radial basis function interpolants.
\newblock {\em SIAM Journal on Scientific Computing}, 34(2):A737--A762, 2012.

\bibitem{fornberg2013stable}
Bengt Fornberg, Erik Lehto, and Collin Powell.
\newblock Stable calculation of gaussian-based rbf-fd stencils.
\newblock {\em Computers \& Mathematics with Applications}, 65(4):627--637,
  2013.

\bibitem{wright2003radial}
Grady~Barrett Wright.
\newblock {\em Radial basis function interpolation: numerical and analytical
  developments}.
\newblock University of Colorado at Boulder, 2003.

\bibitem{wright2017stable}
Grady~B Wright and Bengt Fornberg.
\newblock Stable computations with flat radial basis functions using
  vector-valued rational approximations.
\newblock {\em Journal of Computational Physics}, 331:137--156, 2017.

\bibitem{godunov1959difference}
Sergei~Konstantinovich Godunov.
\newblock A difference method for numerical calculation of discontinuous
  solutions of the equations of hydrodynamics.
\newblock {\em Matematicheskii Sbornik}, 89(3):271--306, 1959.

\bibitem{lee2017new}
Dongwook Lee, Adam Reyes, Carlo Graziani, and Petros Tzeferacos.
\newblock New high-order methods using gaussian processes for computational
  fluid dynamics simulations.
\newblock In {\em Journal of Physics: Conference Series}, volume 837, page
  012018. IOP Publishing, 2017.

\bibitem{balsara1999staggered}
Dinshaw~S Balsara and Daniel~S Spicer.
\newblock A staggered mesh algorithm using high order godunov fluxes to ensure
  solenoidal magnetic fields in magnetohydrodynamic simulations.
\newblock {\em Journal of Computational Physics}, 149(2):270--292, 1999.

\bibitem{mignone2011pluto}
Andrea Mignone, C~Zanni, Petros Tzeferacos, B~Van~Straalen, P~Colella, and
  G~Bodo.
\newblock The pluto code for adaptive mesh computations in astrophysical fluid
  dynamics.
\newblock {\em The Astrophysical Journal Supplement Series}, 198(1):7, 2011.

\bibitem{padioleau:tel-03130146}
Thomas Padioleau.
\newblock {\em {Development of ''all-r{\'e}gime'' AMR simulation methods for
  fluid dynamics, application in astrophysics and two-phase flows}}.
\newblock Theses, {Universit{\'e} Paris-Saclay}, December 2020.

\bibitem{lee2021single}
Youngjun Lee and Dongwook Lee.
\newblock A single-step third-order temporal discretization with
  {Jacobian}-free and {Hessian}-free formulations for finite difference
  methods.
\newblock {\em Journal of Computational Physics}, 427:110063, 2021.

\bibitem{LEE2021100098}
Youngjun Lee, Dongwook Lee, and Adam Reyes.
\newblock A recursive system-free single-step temporal discretization method
  for finite difference methods.
\newblock {\em Journal of Computational Physics: X}, 12:100098, 2021.

\bibitem{gottlieb2011strong}
Sigal Gottlieb, David~I Ketcheson, and Chi-Wang Shu.
\newblock {\em Strong stability preserving {Runge-Kutta} and multistep time
  discretizations}.
\newblock World Scientific, 2011.

\bibitem{toro1994restoration}
Eleuterio~F Toro, M~Spruce, and W~Speares.
\newblock Restoration of the contact surface in the {HLL}-{R}iemann solver.
\newblock {\em Shock Waves}, 4(1):25--34, 1994.

\bibitem{harten1983upstream}
Amiram Harten, Peter~D Lax, and Bram~van Leer.
\newblock On upstream differencing and godunov-type schemes for hyperbolic
  conservation laws.
\newblock {\em SIAM review}, 25(1):35--61, 1983.

\bibitem{quirk1997contribution}
James~J Quirk.
\newblock A contribution to the great riemann solver debate.
\newblock In {\em Upwind and High-Resolution Schemes}, pages 550--569.
  Springer, 1997.

\bibitem{shu1998essentially}
Chi-Wang Shu.
\newblock Essentially non-oscillatory and weighted essentially non-oscillatory
  schemes for hyperbolic conservation laws.
\newblock In {\em Advanced numerical approximation of nonlinear hyperbolic
  equations}, pages 325--432. Springer, 1998.

\bibitem{shu1989efficient}
Chi-Wang Shu and Stanley Osher.
\newblock Efficient implementation of essentially non-oscillatory
  shock-capturing schemes, {II}.
\newblock In {\em Upwind and High-Resolution Schemes}, pages 328--374.
  Springer, 1989.

\bibitem{sedov1993similarity}
Leonid~Ivanovich Sedov.
\newblock {\em Similarity and dimensional methods in mechanics}.
\newblock CRC press, 1993.

\bibitem{fryxell2000flash}
Bruce Fryxell, Kevin Olson, Paul Ricker, FX~Timmes, Michael Zingale, DQ~Lamb,
  Peter MacNeice, Robert Rosner, JW~Truran, and H~Tufo.
\newblock {FLASH}: An adaptive mesh hydrodynamics code for modeling
  astrophysical thermonuclear flashes.
\newblock {\em The Astrophysical Journal Supplement Series}, 131(1):273, 2000.

\bibitem{woodward1984numerical}
Paul Woodward and Phillip Colella.
\newblock The numerical simulation of two-dimensional fluid flow with strong
  shocks.
\newblock {\em Journal of computational physics}, 54(1):115--173, 1984.

\bibitem{hui1999unified}
WH~Hui, PY~Li, and ZW~Li.
\newblock A unified coordinate system for solving the two-dimensional {Euler}
  equations.
\newblock {\em Journal of Computational Physics}, 153(2):596--637, 1999.

\bibitem{liska2003comparison}
Richard Liska and Burton Wendroff.
\newblock Comparison of several difference schemes on {1D} and {2D} test
  problems for the {Euler} equations.
\newblock {\em SIAM Journal on Scientific Computing}, 25(3):995--1017, 2003.

\bibitem{persson2013shock}
Per-Olof Persson.
\newblock Shock capturing for high-order discontinuous galerkin simulation of
  transient flow problems.
\newblock In {\em 21st AIAA computational fluid dynamics conference}, page
  3061, 2013.

\bibitem{stone2020athena++}
James~M Stone, Kengo Tomida, Christopher~J White, and Kyle~G Felker.
\newblock The athena++ adaptive mesh refinement framework: Design and
  magnetohydrodynamic solvers.
\newblock {\em The Astrophysical Journal Supplement Series}, 249(1):4, 2020.

\bibitem{balsara2012self}
Dinshaw~S Balsara.
\newblock Self-adjusting, positivity preserving high order schemes for
  hydrodynamics and magnetohydrodynamics.
\newblock {\em Journal of Computational Physics}, 231(22):7504--7517, 2012.

\bibitem{ha2008positive}
Youngsoo Ha and Carl~L Gardner.
\newblock Positive scheme numerical simulation of high mach number
  astrophysical jets.
\newblock {\em Journal of Scientific Computing}, 34(3):247--259, 2008.

\bibitem{wu2019provably}
Kailiang Wu and Chi-Wang Shu.
\newblock Provably positive high-order schemes for ideal magnetohydrodynamics:
  analysis on general meshes.
\newblock {\em Numerische Mathematik}, 142(4):995--1047, 2019.

\bibitem{liu2021new}
Mengqing Liu, Man Zhang, Caixia Li, and Fang Shen.
\newblock A new locally divergence-free wls-eno scheme based on the
  positivity-preserving finite volume method for ideal mhd equations.
\newblock {\em Journal of Computational Physics}, page 110694, 2021.

\bibitem{lee2009unsplit}
Dongwook Lee and Anil~E Deane.
\newblock An unsplit staggered mesh scheme for multidimensional
  magnetohydrodynamics.
\newblock {\em Journal of Computational Physics}, 228(4):952--975, 2009.

\bibitem{lee2013solution}
Dongwook Lee.
\newblock A solution accurate, efficient and stable unsplit staggered mesh
  scheme for three dimensional magnetohydrodynamics.
\newblock {\em Journal of Computational Physics}, 243:269--292, 2013.

\bibitem{fleischmann2019numerical}
Nico Fleischmann, Stefan Adami, and Nikolaus~A Adams.
\newblock Numerical symmetry-preserving techniques for low-dissipation
  shock-capturing schemes.
\newblock {\em Computers \& Fluids}, 189:94--107, 2019.

\end{thebibliography}
\end{document}